\documentstyle[epsf,12pt,fullpage,subfigure] {article}
\font\tenmsb=msbm10
\font\sevenmsb=msbm7
\font\fivemsb=msbm5
 
\newfam\msbfam 
\textfont\msbfam=\tenmsb
\scriptfont\msbfam=\sevenmsb
\scriptscriptfont\msbfam=\fivemsb

\def\Bbb#1{\fam\msbfam\relax#1}
 
\newtheorem{thm}{Theorem}[subsection]
\newtheorem{prop}[thm]{Proposition}
\newtheorem{cor}[thm]{Corollary}
\newtheorem{lem}[thm]{Lemma}
\newtheorem{conj}[thm]{Conjecture}
\newtheorem{exa}[thm]{Example}
\newtheorem{defn}[thm]{Definition}
\newtheorem{clm}[thm]{Claim}

\newcommand{\ben}{\begin{enumerate}}
\newcommand{\een}{\end{enumerate}}
\newcommand{\ble}{\begin{lem}}
\newcommand{\ele}{\end{lem}}
\newcommand{\bcl}{\begin{clm}}
\newcommand{\ecl}{\end{clm}}
\newcommand{\bthm}{\begin{thm}}
\newcommand{\ethm}{\end{thm}}
\newcommand{\bpr}{\begin{prop}}
\newcommand{\epr}{\end{prop}}
\newcommand{\bco}{\begin{cor}}
\newcommand{\eco}{\end{cor}}
\newcommand{\bcon}{\begin{conj}}
\newcommand{\econ}{\end{conj}}
\newcommand{\bde}{\begin{defn}}
\newcommand{\ede}{\end{defn}}
\newcommand{\bex}{\begin{exa}}
\newcommand{\eex}{\end{exa}}
\newcommand{\barr}{\begin{array}}
\newcommand{\earr}{\end{array}}
\newcommand{\btab}{\begin{tabular}}
\newcommand{\etab}{\end{tabular}}
\newcommand{\beq}{\begin{equation}}
\newcommand{\eeq}{\end{equation}}
\newcommand{\bea}{\begin{eqnarray*}}
\newcommand{\eea}{\end{eqnarray*}}
\newcommand{\bce}{\begin{center}}
\newcommand{\ece}{\end{center}}
\newcommand{\bpi}{\begin{picture}}
\newcommand{\epi}{\end{picture}}
\newcommand{\bfi}{\begin{figure} \begin{center}}
\newcommand{\efi}{\end{center} \end{figure}}

\newcommand{\bsl}{\begin{slide}{}}
\newcommand{\esl}{\end{slide}}

\newcommand{\bib}{thebibliography}

\newcommand{\qed}{\rule{1ex}{1ex}}

\newcommand{\sbs}{\subset}
\newcommand{\sbe}{\subseteq}

\long\def\forget#1\forgotten{}

\newcommand{\al}{\alpha}
\newcommand{\be}{\beta}
\newcommand{\ga}{\gamma}
\newcommand{\de}{\delta}
\newcommand{\ep}{\epsilon}

\newcommand{\la}{\lambda}

\newcommand{\si}{\sigma}

\newcommand{\Ga}{\Gamma}
\newcommand{\De}{\Delta}

\newcommand{\Si}{\Sigma}

\newcommand{\Z}{{\Bbb Z}}
\newcommand{\R}{{\Bbb R}}
\newcommand{\C}{{\Bbb C}}
\newcommand{\F}{{\Bbb F}} 
\newcommand{\PP}{{\Bbb P}}

\newcommand{\cB}{{\cal B}}

\newcommand{\cL}{{\cal L}}

\newcommand{\cR}{{\cal R}}

\begin{document}

\title{\textsc{The Fundamental Group's Structure of the Complement of 
           Some Configurations of Real Line Arrangements}}
\author{David Garber and Mina Teicher \\ Department of Mathematics and Computer Sciences \\ Bar-Ilan University }

\date{\today \\[1in]
               }
\maketitle

\begin{abstract}
In this paper, we give a fully detailed exposition of computing 
fundamental groups of complements of line arrangements using 
the Moishezon-Teicher technique for computing the braid monodromy  of a curve 
and the Van-Kampen theorem which induces a presentation of the 
fundamental group of the complement from the braid monodromy of the curve.
  For example, we treated the cases 
where the arrangement has $t$ multiple intersection points  and  the 
rest are simple intersection points.
In this case, the fundamental group of the complement
is a direct sum of infinite cyclic groups and $t$ free groups. 
Hence, the fundamental groups in these cases is ``big''.
These calculations will be useful
in computing the fundamental group of Hirzebruch covering surfaces.

\end{abstract}
\newpage

\tableofcontents

\newpage

\section{Introduction}

In this paper, we give a fully detailed exposition of calculations
of fundamental groups of
the complements of certain configurations of real line arrangements
using the Moishezon-Teicher algorithm (which calculates the braid monodromy 
of curves),  the Van-Kampen theorem
(which induces a finite presentation, in terms of generators and relations, 
of the fundamental group of curves' complements, 
from its braid monodromy), and some group computations.

\medskip

In particular, we got:
\begin{enumerate}
\item Let $\cL$ be  a real line arrangement which is a union of 
 $t$ subsets of lines  each of which consists of
 $k_i+1$ lines meeting in a single point, and any two lines belonging 
to different subsets meet in a simple point.  Then:
$$\pi _1 (\C ^2 -\cL, u_0 ) \cong ({\bigoplus _{i=1} ^t} \F ^{k_i}) \oplus \Z ^t$$
and
$$\pi _1 (\C\PP ^2 -\cL, u_0 ) \cong ({\bigoplus _{i=1} ^t} \F ^{k_i}) \oplus \Z^{t-1} $$
\item  Let $\cL$ be a real line arrangement
which consists of $t$ subsets of lines each of which consists of
 $k_i+1$ lines meeting in a single point and all the $t$ multiple points
lie on the same line $L \in \cL$. Then:
$$\pi _1 (\C ^2 -\cL, u_0 ) \cong ({\bigoplus _{i=1} ^t} \F ^{k_i}) \oplus \Z$$
and
$$\pi _1 (\C\PP ^2 -\cL, u_0 ) \cong {\bigoplus _{i=1} ^t} \F ^{k_i}$$
\item Generalizations:
Let $\cL$ be a real line arrangement in $\C\PP ^2$ consists of $n$ lines. We
choose the line at infinity such that all the lines are intersected in $\C ^2$.
Assume that there are $k$ multiple intersection points $p_1, \cdots ,p_k$
with multiplicities $m_1, \cdots, m_k$ respectively. 
Assume also that all the multiple 
intersection points in every equivalence class (of multiple points) 
are collinear, i.e. in every 
equivalence class (of multiple points) 
there is a unique line of $\cL$ which all the multiple
points of that class lie on it. 
Then: 
$$\pi _1 (\C ^2 -\cL, u_0 ) \cong {\bigoplus_{i=1}^k} \F ^{m_i -1} 
\oplus \Z ^{n-({\sum_{i=1}^k} (m_i -1))}$$ 
and     
$$\pi _1 (\C\PP ^2 -\cL, u_0 ) \cong {\bigoplus_{i=1}^k} \F ^{m_i -1} 
\oplus \Z ^{n-1-({\sum_{i=1}^k} (m_i -1))}$$ 
The number of infinite cyclic groups in the affine case is a sum of 
two numbers: the number of 
equivalence classes (see definitions in section 5) 
and the number of lines which have only simple intersection
points.
\item Therefore, in all the above cases, the fundamental group 
is ``big''.
\end{enumerate}

\bigskip

We will organize the paper as follows: in  section 2, 
we introduce the needed background for the techniques which will be used, 
and we give a detailed description of the Moishezon-Teicher algorithm
for the case of line arrangements and the Van-Kampen theorem. 

In section 3, we compute 
 the structure of the fundamental group 
of the complement of a line arrangement which consists of $t$ subsets 
of lines  and the multiple points are not collinear. 

In section 4, we compute
the structure of  the fundamental group of the complement 
of a line arrangement which consists of $t$ subsets of lines  and 
the multiple points are collinear. 

In section 5, we generalize the results of the calculations of sections 
3 and 4. In section 6, we discuss the bigness of the groups which 
have been treated.

\section{Preliminaries}

\subsection{Some background}

This topic starts with Zariski, who proved in [Z, p. 317] that:

\medskip

\noindent
{\bf Proposition (Zariski)}\\
{\it The fundamental group of the complement of $n$ lines in 
general position is abelian.}

\medskip

Among the modern works on this topic, one can mention [Fa1], [Fa2], [OS], [Sa],
[Ra] and more.

Moishezon and Teicher developed an algorithm for computing
fundamental groups of complements of branch curves of generic projection of 
surfaces of general type (see [MoTe1],[MoTe2]). 
This algorithm can be used also for computing 
fundamental groups of complement of line arrangements. In this paper we
give a detailed exposition of this technique in some configurations
of line arrangements.  

\medskip

Simultaneously and independently, by entirely different methods, 
Fan proved in [Fa1],[Fa2] the following results for the projective case:

\medskip

\noindent
{\bf Proposition (Fan)} \\
{\it Let $\Si=\bigcup l_i$ be a line arrangement in $\C \PP^2$ and assume that there 
is 
a line $L$ of $\Si$ such that for any singular point $S$ of $\Si$ 
with multiplicity $\geq$ 3, we have $S \in L$. Then: $\pi _1(\C \PP ^2 - \Si)$ 
is isomorphic to a direct product of free groups. }

\medskip

\noindent
{\bf Proposition (Fan)} \\
{\it Let $\Si$ be an arrangement of $n$ lines and $S=\{a_1, \cdots, a_k\} $
be the set of all singularities of $\Si$ with multiplicity $\geq 3$. 
Suppose that 
$\be (\Si)=0$, where $\be (\Si)$ is the first Betti number of the subgraph of $\Si$ 
which contains only the higher singularities 
(i.e. with multiplicity $\geq 3 $) and their 
edges. Then:
$$\pi _1 (\C \PP ^2 - \Si) \cong \Z ^r \oplus \F ^{m(a_1)-1} \oplus \cdots \oplus \F ^{m(a_k)-1}$$
where $r=n+k-1-m(a_1)- \cdots - m(a_k)$.}

\medskip

It has to be noted that the assumption $\be (\Si)=0$ is equivalent to the
assumption that $\Si$ is a union of
trees. The $r$ in the last proposition is actually a sum of two 
combinatorial ingredients: the number of the 
trees in $\Si$ minus 1 and the number of lines which are intersected only in 
simple intersection points.

\medskip

 Oka and Sakamoto proved in [OS] the following theorem, which will be a useful
tool in some of our calculations:

\medskip

\noindent
{\bf Theorem  (Oka-Sakamoto)} \\
{\it Let $C_1$ and $C_2$ be algebraic plane curves in $\C ^2$. 
Assume that the intersection $C_1 \cap C_2$ 
consists of distinct $d_1 \cdot  d_2$ points, where $d_i \ (i=1,2)$ are the 
respective degrees of $C_1$ and $C_2$.\\
Then:
$$\pi _1 (\C ^2 - (C_1 \cup C_2)) \cong \pi _1 (\C ^2 -C_1) \oplus \pi _1 (\C ^2 -C_2)$$}

\medskip

Our computations on the fundamental groups of complements of 
line arrangements have 
applications to the fundamental groups of complements of branch curves, 
which is an important invariant of surfaces [Te2] (when we degenerate 
a surface to a union of planes, the branch curve 
degenerates to a union of lines). Moreover, the methods of this 
paper are important tools
 in the computations of  the fundamental groups of Hirzebruch covering
surfaces.

\subsection{Definition of g-base}

Here, we will present the required definitions and results for the
presentation of the algorithm of Moishezon-Teicher. 
We  follow the presentation of [MoTe1].

\medskip

In this section, we will define the notion of {\it g-base (good geometric
base)} for $\pi _1(D-K,*)$, where $K$ is a finite set in a disk $D$.
For this definition, we have to define:

\bde \label{l_si} {\bf $l( \ga )$}  \\
Let $D$ be a disk. Let $w_i, \ i=1, \cdots , n$, be small disks in ${\rm Int}(D)$ such
that:
$$ \ w_i \cap w_j = \emptyset , \forall i \not = j.$$ 
Let $u \in \partial D$.
Let $\ga$ be a simple path connecting $u$ with one of the $w_i$'s, 
say $w_{i_0}$, which does not meet any other $w_j,\ j \not = i_0$. \\
We assign to $\ga$ a loop $l( \ga )$ (actually an element of $\pi _1(D-K,u)$)
as follows: let $c$ be a simple loop equal to the (oriented) boundary of a 
small neighbourhood $V$ of $w_{i_0}$ chosen such that $\ga ' = \ga - V \cap \ga$ is
a simple path. \\
 Then: $l( \ga ) = \ga ' \cup c \cup (\ga ')^{-1}$ (we will not 
distinguish between $l( \ga )$ and its representative in $\pi _1(D-K,u)$).
\ede

\medskip

\begin{center}
\epsfysize=4cm
\epsfbox{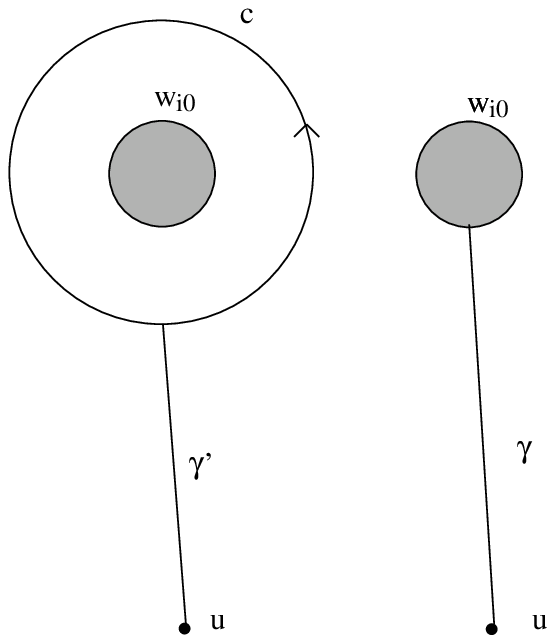}
\end{center}

\medskip

\bde \label{g_base} {\bf Bush, g-base (good geometric base)} \\
Let $D$ be a disk, $K \sbs D$, $\# K < \infty$. Let $u \in D-K$.  
A set of simple paths $\{ \ga _i \}$ is a {\bf bush} in $(D,K,u)$, if 
$\forall i,j,\ \ga _i \cap \ga _j = u; \ \forall i,\ \ga _i \cap K$ 
= one point, and $\ga _i$ are ordered counterclockwise around u. 
Let $\Ga _i = l(\ga _i) \in \pi _1(D-K,u)$ be a loop around $K \cap \ga _i$ 
determined by $\ga _i$. $\{ \Ga _i \}$ is called a {\bf g-base} of $\pi _1(D-K,u)$. 
\ede

\medskip

\begin{center}
\epsfysize=5cm
\epsfbox{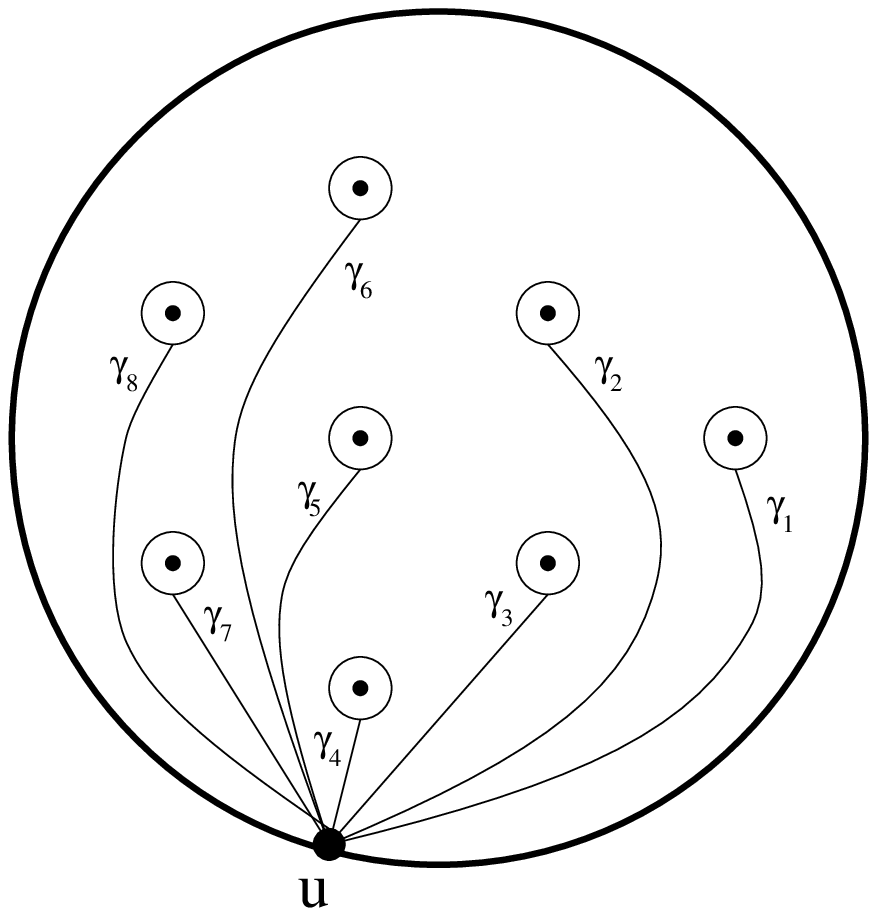}  
\end{center}

\medskip

\subsection{Braid group and braid monodromy}

Let $D$ be a closed disk in $\R^2$, $K \sbs D$ a finite set, $u \in\partial D$.
In such a case, we can define the {\it braid group} $B_n[D,K]\ (n=\# K)$:

\bde {\bf Braid group - $B_n[D,K]$} \\
Let $\cB$ be the group of all diffeomorphisms $\be$ of $D$ such that $\be(K)=K$,
$\be|_{\partial D} = {\rm Id}|_{\partial D}$. Such diffeomorphism acts naturally on
$\pi_1(D-K,u)$. We say that two such diffeomorphisms
are equivalent if they define the same automorphism on $\pi_1(D-K,u)$.
The quotient of $\cB$ by this equivalence relation is called the 
{\bf braid group $B_n[D,K]$}. An element of $B_n [D,K]$ is called a {\bf braid}.
A composition of braids is from {\bf left to right}.
\ede

\medskip

Let us now define the concept of a {\it half-twist braid}. After fixing an 
orientation on $\R^2$, we can define a simple path $\si$ such that 
$[ \si ] \sbe (D-\partial D-K) \cup \{ a,b \}$, $\si$ connects $a$ with $b$
($a,b \in K$).
Choose now a small regular neighbourhood U of $\si$, and an orientation 
preserving diffeomorphism $f:\R^2 \to \C$ ($\C$ is taken with the usual 
``complex'' orientation) such that 
$f(\si)=[-1,1],\ f(U)= \{ z \in \C \ | \ |z|<2 \}$. 
Let $\al (x)$ be any real smooth monotone function such that
$$
\al (x) =  \left\{ \matrix{
            1  &  x \in [0, {3 \over 2} ] \cr
            0  &  x \ge 2
} \right.
$$
 
\noindent
With this function, we define a diffeomorphism $ h: \C \to \C $ as follows:
for any $z=re^{i \varphi} \in \C$, we define: 
$h(z)=re^{i ( \varphi + \al (r) \pi )}$. It is clear that 
$\forall z, |z| \le {3 \over 2}$, $h(z)$ is a positive rotation on $180^\circ$ 
and $h(z)={\rm Id} \ \forall z, |z| \ge 2$. After these preparations, we can define:
 
\bde {\bf $H(\si)$ - (positive) half-twist defined by $\si$} \\
$H(\si)$ is the braid defined by $(f^{-1} \cdot h \cdot f) |_D$. 
\ede

\medskip

We have also another way to look at braids - via {\it motions} of K.
 
\bde {\bf Motion of K' to K} \\
Let $K'= \{ a'_1, a'_2, \cdots ,a'_n \}, \ K= \{ a_1,a_2, \cdots ,a_n \}$. 
{\bf A motion of K' to K in D} is n continuous functions 
$m_i:[0,1] \to D, \ i=1, \cdots ,n$, such that:  \\
(a) $\forall i,\  m_i(0) = a'_i, \ m_i(1) = a_i$.  \\
(b) $\forall i \not = j, \ m_i(t) \not = m_j(t) \ \forall t \in [0,1]$. 
\ede
 
According to the following proposition, we can define a family of
diffeomorphisms induced from the motion (under the condition that K=K').
  
\bpr
Given a motion $\cR$, there exists a continuous family of diffeomorphisms 
$D_{\cR ,t} : D \to D, \ t \in [0,1]$, such that: \\
(a) $D_{\cR ,t}|_{\partial D} = {\rm Id}|_{\partial D}$.  \\
(b) $\forall t,i,\  D_{\cR ,t}(a'_i)=m_i(t)$.
\epr

\bde \label{b_m} {\bf $b_\cR$ (braid induced from a motion $\cR$)} \\
When K=K', $b_\cR$ is the braid defined by the diffeomorphism $D_{\cR ,1}$.
\ede

\medskip

We define another important notion: 

\bde {\bf Skeleton in $(D,K,K'')$} \\
Let $K'' \sbs K, K'' = \{ b_1, \cdots , b_m \}$. A {\bf skeleton} in $(D,K,K'')$ 
is represented by a consecutive sequence of simple paths 
$( p_1, \cdots , p_{m-1})$ in $D- \partial D$ such that  each $p_i$ connects $b_i$ to
$b_{i+1}$.
We say that two such sequences, say 
$( p_1, \cdots , p_{m-1}), ( \tilde p_1, \cdots , \tilde p_{m-1})$, represent the
same skeleton, if $H( p_i ) = H ( \tilde p_i ),  i =1, \cdots, m-1$.
\ede  

Before introducing the definition of {\it braid monodromy}, we have to make
some more constructions. From now, we will work in $\C^2$. Let $E$ (resp. $D$) 
be a closed disk on $x$-axis (resp. $y$-axis), 
and let $C$ be a part of an algebraic
curve in $\C^2$ located in $E \times D$. Let $\pi _1 : E \times D \to E$ and
$\pi _2 : E \times D \to D$ be the canonical projections, and let 
$\pi = \pi _1 |_C: C \to E$. Assume $\pi$ is a proper map, and $\deg \pi = n$. 
Let $ N = \{ x \in E \ |\ \# \pi ^ {-1}(x) < n \} $, and assume 
$ N \cap \partial E= \emptyset $. Now choose $ M \in \partial E$ 
and let $K = K(M) = \pi ^ {-1} (M)$. By the assumption that 
$\deg \pi = n \ \ (\Rightarrow \#  K=n)$,  we can write: 
$K = \{ a_1,a_2, \cdots ,a_n \} $.
Under these constructions, from each loop in $E-N$, we can define a braid in
$B_n [M \times D, K]$ in the following way: 
\begin{itemize}
\item[(1)] Because $ \deg \pi = n$, we can lift any loop in $E-N$ with 
a base point $M$ to a system of $n$ paths in $ (E-N) \times D $ which 
start and finish at $ \{ a_1,a_2, \cdots ,a_n \} $.
\item[(2)] Project this system into $D$ (by $\pi _2$), to get $n$ paths in $D$ 
which start and end at the image of $K$ in $D$ (under $\pi_2$). These paths 
actually form a motion.
\item[(3)] Induce a braid from this motion, as we did in definition \ref{b_m}.
\end{itemize}

\medskip

To conclude, we can match a braid to each loop. 
Therefore, we get a map $ \varphi : \pi _1(E-N,M) \to B_n [M \times D, K]$,
which is also a group homomorphism which is called the
{\bf braid monodromy of $C$ with respect to $E \times D,\pi _1,M$}.      

\medskip

For the next definitions, let us assume $M_0,M_1 \in E-N$ and $T:[0,1] \to E-N$
be a path which connects $M_0$ with $M_1$. We know that there exists a 
continuous family of diffeomorphisms 
$\psi _{(t)}: M_0 \times D \to T(t) \times D, \ \forall t \in [0,1]$, such that: \\
(a) $\psi _{(0)} = {\rm Id}|_{M_0 \times D}$. \\
(b) $\forall t \in [0,1], \ \psi _{(t)}(\pi^{-1} _1(M_0) \cap C) =
                            \pi^{-1} _1(T(t)) \cap C$.\\
(c) $\forall y \in \partial D, \ \psi _{(t)}(M_0,y) = (T(t),y)$.

\bigskip

In this situation, we can define the {\it Lefschetz diffeomorphism induced by T}:

\bde {\bf $\psi _T$, Lefschetz diffeomorphism induced by T}
$$\psi _T = \psi _{(1)}: M_0 \times D\  \tilde\to \ M_1 \times D $$
\ede

Let $s=(x(s),y(s)) \in C$ be a singular point of $\pi$ 
(i.e. $x(s) \in N$). Let $D'(s)$ be such a small disk on $y$-axis centered at $y(s)$ 
that $(x(s) \times D'(s)) \cap C = s$, i.e. there are no other branches of $C$ which 
intersect $D'(s)$. Therefore, for any sufficiently small neighbourhood $U$
of $x(s)$ on the $x$-axis centered at $x(s)$ such that $\forall x \in U -x(s)$,
$\# (x \times {\rm Int}(D'(s))) \cap C$ is independent of $x$ (we call this number the
{\it local degree of $\pi$ at $s$} and denote it by $\deg _s \pi$). Let $k=\deg _s \pi$
and $E'$ be a small closed disk on the $x$-axis centered at $x(s)$, such that
$\forall x  \in E' -x(s)$, $\# (x \times {\rm Int}(D'(s))) \cap C=k$. Choose a point $a(s) \in \partial E'$
and let $T: [0,1] \to \C$ be a path in $E-N-{\rm Int}(E')$ connecting $a(s)$ to a point 
$M' \in E-N$. Let $K_{a,s} =(a(s) \times D'(s)) \cap C$.

\bde {\bf $\tilde\psi _T$, Lefschetz embedding induced by T}\\
Let $\psi _T$ be the Lefschetz diffeomorphism as defined above. 
Let $T$ be as above, $a=a(s)$, $D'=D'(s)$.
Then:
$$\tilde\psi _T = \psi _T | _{a \times D'}: a \times D' \to M' \times D$$
\ede   

\noindent
{\bf Remark:} Take $k$ liftings of $T$ to $C$ starting at the 
different points of\\ 
$K_{a,s} =(a \times D') \cap C$. These liftings are real curves in $T \times D$. We can
think of $\tilde\psi _T$ as ``pulling'' of $a \times D'$ in $T \times D$ along 
these real curves.

\bde {\bf $\cL _{T,s}$, Lefschetz injection induced by T}\\
Consider $\tilde\psi _T : a \times D' \to M' \times D$, Lefschetz embedding induced by T.\\
Let $K(M') = (M' \times D) \cap C$. We have
$$\tilde\psi _T (K_{a,s}) \sbs K(M'), (K(M') - \tilde\psi _T (K_{a,s})) \cap \tilde\psi _T ({\rm Int} (D')) = \emptyset$$
Therefore, 
the following canonical injection is well defined:
$$\cL _{T,s} = \psi _T ^{\vee} : B _k [a \times D' , K _{a,s}] \hookrightarrow
B _n [M' \times D ,K]$$
\ede

In order to define the {\it Lefschetz vanishing cycle}, we need the following 
definition:

\bde {\bf Linear frame of a braid group $B_n [D,K]$}\\
Let $K= \{ a_1, a_2, \cdots ,a_n \}$. Let $\{ \xi _1,\xi _2, \cdots ,\xi _{n-1} \}$ be
a system of straight line segments in $D- \partial D$ such that 
each $\xi _i$ connects
$a_i$ with $a_{i+1}$ (and does not intersect any other $\xi _j$ except of
end points). Let $H_i=H(\si _i)$. The ordered system of positive half-twists
$( H_1, H_2, \cdots ,H_{n-1} )$ is called {\bf a linear frame of $B_n[D,K]$
defined by $\{ \xi _1,\xi _2, \cdots ,\xi _{n-1} \}$}.
\ede

Now, we come to one of the most important definitions:

\bde {\bf $\cL$.V.C.$(T,H')$, Lefschetz vanishing cycle induced by $T$}\\
We call $\cL$.V.C.$(T,H')$ a skeleton $<\xi _1, \cdots, \xi _{k-1}>$ in 
$(M' \times D, K ,  \tilde\psi _T (K_{a,s}))$ corresponding $\cL _{T,s}$ and a linear frame 
$(H') = (H _1', \cdots, H'_{k-1})$ of $B _k [a \times D' , K _{a,s}]$, that is 
$\cL _{T,s} (H' _i)=H(\xi _i), \ i=1, \cdots, k-1$.
\ede   

Because of the fact that such a linear frame is unique only when 
all the points of $K$ are 
on a straight line in $D \sbs \R ^2$, $\cL$.V.C.$(T,H')$ 
will be well defined if
all the points of $K _{a,s}$ are on a straight line in $a \times \C$. If all the
points of $K _{a,s}$ are real,  we will
choose the unique linear frame $(H' _1, \cdots, H'_{k-1})$ determined by an
increasing sequence of consecutive real segments on the real axis
of $a \times \C$.

\subsection{The braid monodromy of a real line arrangement}

\bde {\bf Line arrangement in $\C \PP ^2$} \\
A {\bf Line arrangement in $\C \PP ^2$} is an algebraic curve in $\C \PP ^2$
which is a union of projective lines.
\ede

If the lines are given by the linear forms $l_1,l_2, \cdots, l_k$, the union
of the lines is the reducible curve defined by
$$l_1 l_2 \cdots l_k = 0$$

We say that the arrangement is {\it real} if each line
can be defined by an equation with real coefficients (i.e. each linear form
$l_i$ has real coefficients).

Let $\C ^2 = \C \PP ^2 -({\rm projective \ line})$ be an affine part of $\C \PP ^2$.
Let $E$ (resp. $D$) be a closed disk on $x$-axis (resp. on $y$-axis) with
the center on the real part of $x$-axis (resp. $y$-axis).
Let $\pi _1 : E \times D \to E, \pi _2 : E \times D \to D$ be the canonical projections.

\bde {\bf Real line arrangement in a polydisk $E \times D$}\\
We say that $C$ is a {\bf real line arrangement in a polydisk $E \times D$} 
(as above),
if there exists a real line arrangement $\hat C$ in $\C \PP ^2$, such that: \\
(a) $C=\hat C \cap (E \times D)$.\\
(b) $\forall x \in E, \ \pi _1 ^{-1} (x) \cap C \sbs x \times {\rm Int}(D)$.
\ede

Let $\pi = \pi _1 | _C, \ n=\deg \pi$ (=number of lines in $C$),
$N= \{ x \in E\ |\ \# \pi ^{-1} (x) < n \}, \ K _x = \pi ^{-1} (x)$. Therefore,
for any real $x \not\in N$, we have $n$ distinct real points
$(x,y_i(x)), \ 1 \leq i \leq n,$ in $K_x$. We choose a numeration in
$\{ y_1(x), \cdots , y_n(x) \}$, such that $y_1 (x) < y_2(x) < \cdots < y_n (x)$.

Let $\tilde D =\{ z \in \C \ | \ |z - {n+1 \over 2} | \leq {n+1 \over 2} \}, \
\tilde K = \{ 1, 2, \cdots , n \} \sbs \tilde D$ ($\tilde D$ is a model which simplifies
the treatment with the theoretic calculations of the braid monodromy).
Let $\tilde H = (\tilde H_1, \tilde H_2, \cdots, \tilde H _{n-1})$ be the linear frame
of $B_n [\tilde D,\tilde K]$ defined by the sequence of real segments 
$\tilde\xi = ([1,2],[2,3], \cdots, [n-1,n])$, i.e. $\tilde H_j = H([j,j+1])$.

\medskip

For the set $E' _{\R} = \{ x \in E-N \ | \ x\ {\rm real} \}$, we can construct a set of
diffeomorphisms $\{ \be _x \ | \ x \times D \tilde\to \tilde D \}$ with 
the following properties:
\begin{itemize}
\item[(a)] $\be _x (K_x) =\tilde K$.
\item[(b)] $\be _x (x \times {\rm real\ part\ of}\ D) ={\rm  \ real\ part\ of}\ \tilde D$ 
(order preserved). 
\item[(c)] $\forall x,x' \in E' _{\R} , y \in \partial D, \ \be _x (x,y)=\be _{x'} (x',y)$.
\item[(d)] On each connected component $\tilde\cL$ of $E' _{\R}$, 
$\{ \be _x \ | \ x \in \tilde\cL \}$ is a continuous family of diffeomorphisms.
\end{itemize}   

Let $\xi _x = \{ \xi _{x,1}, \xi _{x,2}, \cdots , \xi _{x,n-1} \} \ (x \in E' _{\R})$ be the sequence
of real segments $[y_i (x), y _{i+1} (x)], \ 1 \leq i \leq n-1,$ 
in $x \times D$ and let 
$H_x = (H _{x,1}, H _{x,2}, \cdots , H _{x,n-1})$ be the linear frame of 
$B_n [x \times D, K_x]$ defined by $\xi _x$. 

Now, we assume that $\forall x_j \in N$, there is only one singular point 
of $C$ over $x_j$.

Let $x_j \in N$. Choose $x' _j = x_j +\ep$, $\ep > 0$ a very small number.
Let $A_j$ be the singularity of $C$ over $x_j$ (i.e. $x(A_j)=x_j$), and let 
$Y_j$ be the union of irreducible components of $C$ containing $A_j$.
In $\{ y_1(x' _j), \cdots , y_n(x' _j) \}$, there is a subsequence with consecutive
indices $\{ y_{k_j} (x' _j), y_{k_j+1} (x' _j), \cdots , y_{l_j}(x' _j) \}$ which is equal to 
$K'_{x'_j} = Y_j \cap (x' _j \times D)$.

In this situation, we can define the following notions:

\bde {\bf Local $\cL$.V.C. of $A_j$ (``Local Lefschetz vanishing cycle of $A_j$'')}\\
A skeleton in $(x' _j \times D , K _{x' _j} , K' _{x' _j})$ represented by the sequence 
of real segments
$$[y_{k_j+r-1} (x' _j), y_{k_j+r} (x' _j)], \ 1 \leq r \leq l_j-k_j$$
is called a {\bf local $\cL$.V.C. of $A_j$}. 
\ede

\bde {\bf $(k_j,l_j)$, Lefschetz pair of $A_j$}\\
The smallest and biggest indices $k_j,l_j$ in the sequence considered above form 
a pair $(k_j,l_j)$, which is called the {\bf Lefschetz pair of $A_j$}.
\ede

Obviously, the local $\cL$.V.C. of $A_j$ is uniquely defined by
the Lefschetz pair $(k_j,l_j)$.

\bde {\bf $<k_j,l_j>$, skeleton representing local $\cL$.V.C. of $A_j$}\\
Denote by $<k_j,l_j>$ the skeleton in $(\tilde D, \tilde K, (k_j, k_j +1, \cdots , l_j))$
represented by consecutive real segments connecting points of 
$(k_j, k_j +1, \cdots , l_j)$.
\ede

\ble\label{bm_lvc2}
Let $\ga$ be a simple path in $E-N$ connecting $x_j$ with $M (\in \partial E)$,
$[x_j, x' _j] \sbs \ga$. Let $\ga '$ be the part of $\ga$ from $x'_j$ to $M$.
Let 
$$\varphi : \pi _1 (E-N, M) \to B_n [M \times D, K_M]$$ 
be the braid monodromy of $C$ 
w.r.t. $E \times D, \pi _1, M$. Let $\Ga$ be the element represented by $l (\ga)$.\\ 
Then:
$$\varphi (\Ga) = \De ^2 <\cL{\rm .V.C.} (\ga ', H(<\xi _x>))>$$
(where, intuitively, $\De < {\rm skeleton}>$ is a generalized half-twist
which is defined according to the skeleton, and $\De ^2 < {\rm skeleton}>$
is applying this half-twist twice).
\ele

\subsection{The algorithm of Moishezon-Teicher}\label{MT_sec}

Following lemma \ref{bm_lvc2}, in order to calculate the braid monodromy,
we have to find the appropriate Lefschetz vanishing cycles. This is given by
the following theorem [MoTe1]:
   
\begin{thm}\label{MT} {\bf (Moishezon-Teicher)}\\
Let $N=\{ x_1, x_2, \cdots, x_q \}$ with 
$x_q < x_{q-1} < \cdots < x_2 < x_1$, 
$M \in \partial E \cap {\rm (real\ axis)}$, with $M>x_1$, 
and $\ep >0$ a very small number. Let $T_j (1 \leq j \leq q)$ 
be the path from $x_j -\ep$ to
$x_j + \ep$ along the semicircle below real axis centered at $x_j$.

\noindent
Let $\ga _j$ be the path from $x_j$ to $M$ defined by 
$$\ga _j = [x_j,x_{j-1}- \ep] \cdot T _{j-1} \cdot [x_{j-1}+ \ep ,x_{j-2}- \ep] \cdot T _{j-2} \cdots T_1 \cdot [x_1,M]$$
\[ (\ga _j = [x_j,x_{j-1}- \ep] \cdot T _{j-1} \cdot ({\prod _{r=j-1} ^2}  [x_r+ \ep ,x_{r-1}- \ep] \cdot T _{r-1}) \cdot [x_1,M] )\]

\medskip

\begin{center}
\epsfysize=2cm
\epsfbox{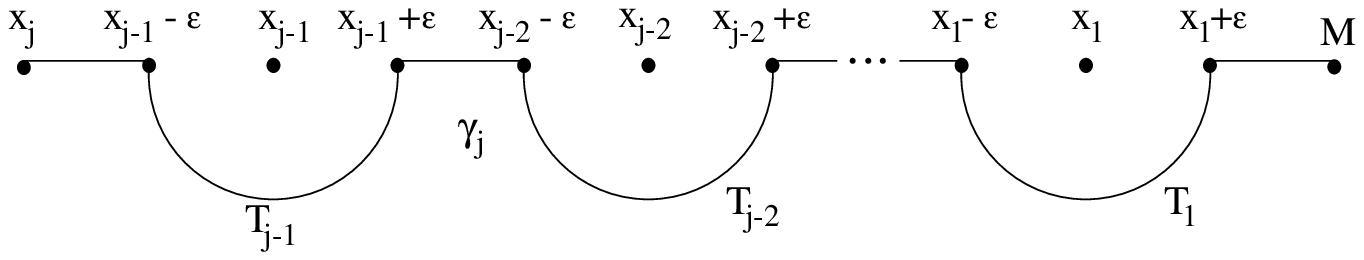}  
\end{center}

\medskip

\noindent
Considering $l (\ga _j)$'s, we get a g-base $\{ \de _1, \de _2, \cdots, \de _q \}$ in 
$\pi _1 (E-N, M)$. 

Assume that for all $x_j$, $1 \leq j \leq q$, there is only one singular point $A_j$
with $x(A_j)=x_j$. Let $(k_j,l_j)$ be the Lefschetz pair of $A_j$, and $<k_j,l_j>$ 
be the skeleton in $(\tilde D, \tilde K , (k_j, k_j+1,\cdots, l_j-1,l_j))$ representing
local $\cL$.V.C. of $A_j$. Let $\ga '_j$  be the part of $\ga _j$
from $x' _j=x_j+\ep$ to $M$.

\noindent
Then:
\[ \cL{\rm .V.C.}(\ga ' _j) = \be _M ^{-1} (<k_j,l_j> \cdot {\prod _{m=j-1} ^1} \De <k_m,l_m>)\]   
$ \displaystyle ( {\rm where \ } {\prod _{m=j-1} ^1} \De <k_m,l_m> =$
$$ \De <k_{j-1},l_{j-1}> \cdot \De <k_{j-2},l_{j-2}> \cdots \De <k_1,l_1> \in B_n[\tilde D,\tilde K])$$
and
$$ \cL{\rm .V.C.}(\ga ' _1) = \be _M ^{-1} (<k_1,l_1>)$$
\end{thm}

\medskip

According to this theorem, in order to compute the braid monodromy of
a line arrangement, we have to do the following steps:
\begin{enumerate}
\item Check that the line arrangement fulfills the assumption that there are no
more than one  intersection point with the same $x$-coordinate (so we can apply
the theorem).
\item Find the Lefschetz pairs of all the intersection points.
\item Calculate the Lefschetz vanishing cycle of every intersection point
according to the last theorem (\ref{MT}).
\item The braid monodromy is the $\De ^2$ of this $\cL$.V.C.
\end{enumerate}

\subsection{The Van-Kampen theorem}\label{VK-section}

The Van-Kampen theorem induces a finite presentation of the fundamental group
of complements of curves by meaning of generators and relations.
From this finite presentation, we will
calculate the structure of the group in our cases
(the original theorem is  in [VK], other versions can be found
at [Mo, pp. 127-130], [MoTe3], [MoTe4, ch. 13], [Te1]. 
The theorems presented here are from [MoTe3],[MoTe4] and [Te1]).

\medskip

Let $S$ be an algebraic curve in $\C^2$ ($p=\deg S$). 
Let $\pi = \pi_1 : \C ^2 \to \C$
be the canonical projection on the first coordinate.   
Let $\C _x = \pi ^{-1}(x)$, and now define: $K_x=\C _x \cap S$ (By assumption
${\rm deg}S=p$, we know $\# K_x \le p$). 

Let $N=\{x\ |\ \# K_x < p \} $.
Choose now $u \in \C$, $u$ real, such that $x \ll u, \ \forall x \in N$, and define:
$B_p = B_p [ \C _u, \C _u \cap S ]$. Let $\varphi _u : \pi _1(\C -N, u) \to B_p$
be the braid monodromy of S w.r.t $\pi, u$.
Also choose  $u_0 \in \C _u, \ u_0 \not\in S$, $u_0$ below real line far
enough such that $B_p$ does not move $u_0$.
It is known that the group $\pi _1(\C _u - S, u_0)$ is free.
There exists an epimorphism
$\pi _1(\C _u - S, u_0) \to \pi _1(\C ^2 - S, u_0)$, so a set of generators 
for $\pi _1(\C _u - S, u_0)$ determines a set of generators for      
$\pi _1(\C ^2 - S, u_0)$. 

\medskip

In this situation, Van-Kampen's theorem says:

\bthm\label{VK1} {\bf Van-Kampen's Theorem - classic version} \\
Let $S$ be an algebraic curve, $u,u_0,\varphi _u$ 
defined as above. Let $\{ \de _i \}$ be a g-base of $\pi _1(\C -N,u)$. 
Let $\{ \Ga _j \ | \ 1 \le j \le p \}\ (p=\deg S)$ be a g-base for
$\pi _1 (\C _u -S,u_0)$. \\
Then, $ \pi _1(\C ^2 -S, u_0)$ is generated
by the images of $\Ga _j$ in $\pi _1(\C ^2 -S, u_0)$ and we get a 
complete set of relations
from those induced from
$$(\varphi _u(\de _i))(\Ga _j) = \Ga _j; \forall i \forall j$$
\ethm

\medskip

Here we present also the classic Van-Kampen theorem for the 
projective case. The only 
difference between the affine case and the projective case is that there 
is one additional relation in the projective case  - the multiplication of
all the generators is equal to the identity of the group.

\bthm\label{VKP} {\bf Van-Kampen's Theorem for projective case - classic version} \\
Let $S$ be an algebraic curve, $u,u_0,\varphi _u$ 
defined as above. Let $\{ \de _i \}$ be a g-base of $\pi _1(\C -N,u)$. 
Let $\{ \Ga _j \ | \ 1 \le j \le p \}\ (p=\deg S)$ be a g-base for
$\pi _1 (\C _u -S,u_0)$. \\
Then, $ \pi _1(\C\PP ^2 -S, u_0)$ is generated
by the images of $\Ga _j$ in $\pi _1(\C ^2 -S, u_0)$ and we get a 
complete set of relations
from those induced from
$$(\varphi _u(\de _i))(\Ga _j) = \Ga _j; \forall i \forall j$$
with one additional relation:
$$\Ga _p \Ga _{p-1} \cdots \Ga _1 = 1$$
\ethm

Oka [O] proved the following connection between the fundamental group 
of the affine case and the fundamental group of the projective case:

\bthm\label{oka} {\bf (Oka)} \\
Let $C$ be a curve in $\C\PP ^2$ and let $L$ be a general line to $C$.
Then, we have a central extension:
$$1 \to \Z \to \pi _1(\C\PP ^2 -(C \cup L)) \to \pi _1(\C\PP^2-C) \to 1$$ 
\ethm

Due to the fact that $L$ is in a general position to $C$, we can say:
$$\pi _1(\C\PP ^2 - (C \cup L)) \cong \pi _1 ((\C\PP ^2-L)-C) \cong \pi _1(\C ^2 -C)$$
(by choosing $L$ as the line at infinity).
Therefore, we get the following short exact sequence (see also [OS]):
$$1 \to \Z \to \pi _1(\C^2 -C) \to \pi _1(\C\PP^2-C) \to 1$$ 

We will show that in the cases which we treat, we get:
$$\pi _1 (\C ^2-C) \cong \pi _1 (\C\PP ^2-C ) \oplus \Z$$ 
and therefore, this short exact sequence splits. 

\medskip

Now we return to the affine case.
In order to give a more precise version of Van-Kampen's theorem for 
cuspidal curves, i.e. for curves with only nodes and cusps as singularities,
we need the following two lemmas.  

\ble\label{lem_av_bv}
Let $V$ be a half-twist in $B_p[D,K]$, $u_0 \not\in K$.\\
Then: there exists $A_V,B_V \in \pi _1(D-K,u_0)$, such that: \\
{\rm (a)} $\{A_V,B_V\}$ can be extended to a g-base of 
    $\pi _1(D-K,u_0)$. \\
{\rm (b)} $V(A_V)=B_V$.

\vspace{5mm}

\epsfysize=8cm
\epsfbox{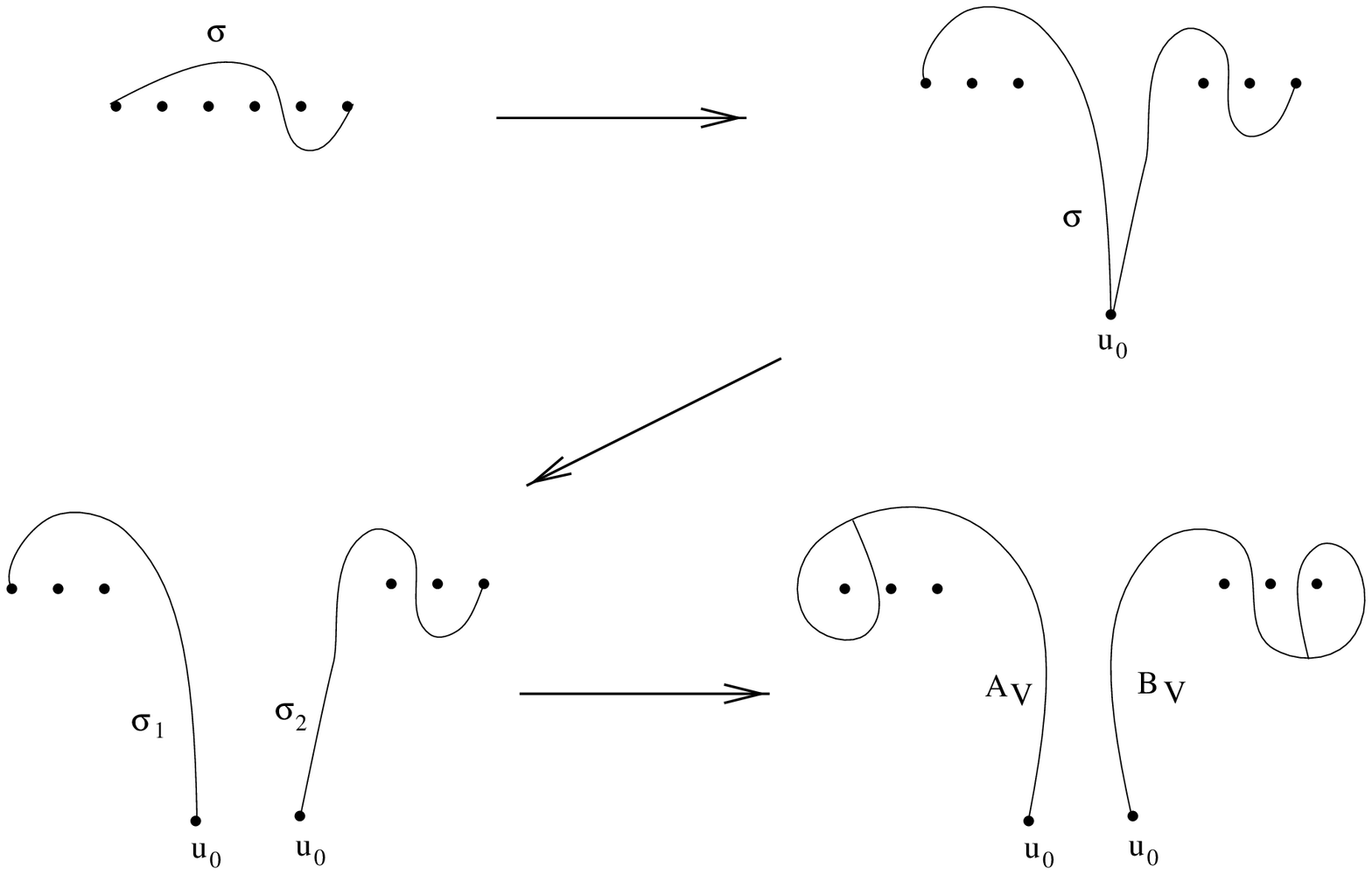}  
\label{av_bv}

\ele
 
Let $S$ be a cuspidal curve in $\C ^2$ ($p=\deg S$). 
We assume that for every $x \in N$ ($N$ as above), 
there is only one singular point
over it (in $\C ^2$). Thus, for every $x \in N$, let $x'$ be the singular point 
over $x$. Because $S$ is a cuspidal curve, the point $x'$ is 
either a branch point, a node or a cusp.

\medskip

\ble \label{lem_VK}
Let $\{ \de _i \}$ be a g-base for $\pi _1(\C - N,u)$.
For every $\de _i$, there exists $V_i$ and $\nu _i$, where $V_i$
is a half-twist and $\nu _i$ is a number such that
$\varphi _u(\de _i)=V_i ^{\nu _i}$. Moreover, $\nu _i =1,2,3$ if $c'_i$ (the
singular point) = a  branch point, a node or a  cusp respectively.
\ele   

We denote:
$$[A,B]=ABA^{-1}B^{-1}$$
$$<A,B>=ABAB^{-1}A^{-1}B^{-1}$$  

\medskip

Now, we can give the precise version of the Van-Kampen theorem for cuspidal
curves:

\bthm\label{VK2} {\bf Van-Kampen's theorem for cuspidal curves} \\
Let $S$ be a cuspidal curve, $u,u_0,\varphi _u,A_{V_i},B_{V_i}$ 
defined as above. Let $\{ \de _i \}$ be a g-base of $\pi _1(\C -N,u)$. 
Let $\varphi _u(\de _i)=V_i ^{\nu _i}$, $V_i$ is a half-twist, 
$\nu _i = 1,2,3$ (as above).

Let $\{ \Ga _j \ | \ 1 \le j \le p \}\ (p=\deg S)$ be a g-base for
$\pi _1 (\C _u -S,u_0)$.\\
Then: $ \pi _1(\C ^2 -S, u_0)$ is generated
by the images of $\Ga _j$ in $\pi _1(\C ^2 -S, u_0)$ and we get a 
complete set of relations
from those induced from $\varphi _u(\de _i) = V_i ^{\nu _i}$, as follows
(when $A_{V_i},B_{V_i}$ are expressed in terms of $\{ \Ga _j \}$): 
\begin{itemize}
\item[{\rm (a)}] $A_{V_i}=B_{V_i}$, when $\nu _i = 1$.
\item[{\rm (b)}] $[ A_{V_i} , B_{V_i} ] =1$, when $\nu _i = 2$. 
\item[{\rm (c)}] $< A_{V_i} , B_{V_i} > =1$, when $\nu _i = 3$. 
\end{itemize}
\ethm

\medskip

What do we get from this theorem? After we calculate the appropriate braid 
monodromy, we can get a finite presentation of the  desired 
fundamental group.

Note that  it is easy to see that the relation,
which is induced from the braid
monodromy, is uniquely determined by the half-twist $V$, and is independent of
the choice of $A_V,B_V$.

\bigskip

Now, we will present the version of Van-Kampen's theorem for 
an arrangement with a single multiple point, i.e. an arrangement where 
all the lines meet in one point (the proof is easy, and can be found, for
example, in [Ga, p. 25]):

\ble\label{vk-cyc} {\bf Van-Kampen's theorem for a single multiple point}\\
Let $l_1, \cdots, l_k$ be $k$ real lines in $\C\PP ^2$ meeting in a single
point $p$.
Let $ \de $ be a loop in $\pi _1 (E-N,u_0)$ around $x(p)$.
Let $\{ \Ga _1, \cdots ,\Ga _k \}$ be a g-base of
$\displaystyle \pi _1 (\C _{u_0} - {\bigcup _{i=1} ^k} l_i)$.\\
Then, the relations which are induced from this intersection point are:
$$\Ga _k \Ga _{k-1} \cdots  \Ga _1 =
 \Ga _1 \Ga _k \cdots  \Ga _3 \Ga _2 = \cdots =
 \Ga _{k-1} \Ga _{k-2} \cdots \Ga _1 \Ga _k$$
\ele

\subsection{An application of the Van-Kampen theorem}

Here, we will prove a simple proposition, which will help
us in the future.
We denote $[x,y]=xyx^{-1} y^{-1}$ for $x,y$ in a group $G$.

\bpr\label{lml}
Let $p$ be an intersection point of $k$ real lines $l_{j_1}, \cdots, l_{j_k}$ 
in $\C\PP ^2$.
Let $ \de  $ be a loop in $\pi _1 (E-N,u_0)$ around $x(p)$. \\
Let $\{ \Ga _{j_1}, \cdots ,\Ga _{j_k} \}$ be a g-base of
$\displaystyle \pi _1 (\C _{u_0} - {\bigcup _{i=1} ^k} l_{j_i})$.\\
Then: the relations which are induced from this intersection point are:
$$[\Ga _{j_k} \Ga _{j_{k-1}} \cdots  \Ga _{j_1},\Ga _{j_i}] = 1;\ 1 \leq i \leq k$$
\epr

\noindent
{\it Proof:}  By the Van-Kampen version for a multiple point
(\ref{vk-cyc}), the following set of relations is induced from
the intersection point $p$:
$$\Ga _{j_k} \Ga _{j_{k-1}} \cdots  \Ga _{j_1} =
 \Ga _{j_{k-1}} \cdots  \Ga _{j_1} \Ga _{j_k} = \cdots =
 \Ga _{j_1} \Ga _{j_k} \cdots \Ga _{j_2}$$

We will prove now that this set of relations is equivalent to the set of
relations in the formulation of the proposition.

\medskip

\noindent
($\Rightarrow$) Let $1 \leq i \leq k$. We have to show that
$$\Ga _{j_k} \Ga _{j_{k-1}} \cdots  \Ga _{j_1} \Ga _{j_i} =
  \Ga _{j_i} \Ga _{j_k} \Ga _{j_{k-1}} \cdots  \Ga _{j_1}$$
We know (from the first set of relations) that
$$(*)\ \Ga _{j_k} \Ga _{j_{k-1}} \cdots  \Ga _{j_1} =
 \Ga _{j_i} \Ga _{j_{i-1}} \cdots \Ga _{j_1} \Ga _{j_k} \cdots \Ga _{j_{i+1}}$$ 
$$(**)\ \Ga _{j_k} \Ga _{j_{k-1}} \cdots  \Ga _{j_1} = 
   \Ga _{j_{i-1}} \Ga _{j_{i-2}} \cdots \Ga _{j_1} \Ga _{j_k} \cdots \Ga _{j_i}$$  
Now: 
$$(\Ga _{j_k} \Ga _{j_{k-1}} \cdots  \Ga _{j_1}) \Ga _{j_i} \stackrel{(*)}{=}
 (\Ga _{j_i} \Ga _{j_{i-1}} \cdots \Ga _{j_1} \Ga _{j_k} \cdots \Ga _{j_{i+1}})
 \Ga _{j_i} = $$
 $$=\Ga _{j_i} (\Ga _{j_{i-1}} \cdots \Ga _{j_1} \Ga _{j_k} \cdots \Ga _{j_{i+1}}
 \Ga _{j_i}) \stackrel{(**)}{=}
 \Ga _{j_i} (\Ga _{j_k} \Ga _{j_{k-1}} \cdots  \Ga _{j_1})$$

\medskip

\noindent
($\Leftarrow$)
From the first relation we have:
$$[\Ga _{j_k} \Ga _{j_{k-1}} \cdots  \Ga _{j_1},\Ga _{j_1}] = 1$$
 i.e.
$\Ga _{j_k} \Ga _{j_{k-1}} \cdots  \Ga _{j_1} \Ga _{j_1} =
\Ga _{j_1} \Ga _{j_k} \Ga _{j_{k-1}} \cdots  \Ga _{j_1}$.
Now, multiply it by $\Ga _{j_1} ^{-1}$ from the right to get:
 $$(***)\ \Ga _{j_k} \Ga _{j_{k-1}} \cdots  \Ga _{j_1} = \Ga _{j_1}   \Ga _{j_k} \cdots \Ga _{j_2}$$
From the second relation we have:
$(\Ga _{j_k} \Ga _{j_{k-1}} \cdots  \Ga _{j_1}) \Ga _{j_2} =
\Ga _{j_2} (\Ga _{j_k} \Ga _{j_{k-1}} \cdots  \Ga _{j_1})$,
but from $(***)$ we get:
$(\Ga _{j_1} \Ga _{j_k}\cdots \Ga _{j_2}) \Ga _{j_2} =
\Ga _{j_2} (\Ga _{j_1}  \Ga _{j_k} \cdots \Ga _{j_2})$.
Now, multiply it  by $\Ga _{j_2} ^{-1}$ from the right to get:
$$\Ga _{j_1}  \Ga _{j_k} \cdots \Ga _{j_2} = \Ga _{j_2} \Ga _{j_1} \Ga _{j_k} \cdots \Ga _{j_3}.$$
Applying the same argument together with the rest of the commutative relations
give us the requested cyclic relations. \hfill $\qed$

\bigskip

\subsection{Outline of the computation of the fundamental group of the complement of line arrangements}\label{outline}

Let us summarize the steps we have
to follow in order to compute the fundamental group of the complement 
of a given real line arrangement $\cL$:
\begin{itemize}
\item[(1)] Calculation of the braid monodromy of $\cL$:
   \begin{itemize}
   \item[-] Check that the line arrangement fulfills the assumption that
   there are no
   more than one  intersection point with the same $x$-coordinate
   (so we can apply the theorem).
   \item[-] Find the Lefschetz pairs of all the intersection points.
   \item[-] Calculate the Lefschetz vanishing cycle of every intersection point
   according to the Moishezon-Teicher theorem.
   \end{itemize}
\item[(2)] Calculation of the relations induced on $\pi _1(\C ^2 -\cL)$ from the
     braid monodromy:
   \begin{itemize}
   \item[-] Choose $u$ as in section \ref{VK-section}.
   \item[-] Choose a g-base for $\pi _1 (\C _u -\cL)$: $\{ \Ga _1, \cdots , \Ga _n\}$.
   \item[-] Calculate the $A_{V_i},B_{V_i}$ from the $\cL$.V.C.
            for every singular point in terms of $\Ga _i, i=1,\cdots, n$.
   \item[-] Find the induced relations according to the Van-Kampen theorem.
   \end{itemize}
\item[(3)] Computing the structure of $\pi _1(\C ^2 -\cL)$ from the relations 
        in (2). This step contains some group calculations and combinatorics.
\end{itemize}

\section{Arrangements with $t$ non-collinear multiple points}

In this section, we are going to calculate the fundamental group
of the complement of line arrangements where there is no line on which there
are two multiple points. Thus, we can divide the arrangement into $t$ subsets 
of lines where all the lines in each subset intersect at a single 
(multiple) point and any two such subsets intersect in  simple points only. 
We define:
\bde {\bf Simple point, multiple point, multiplicity of a point} \\
A {\bf simple point} in a line arrangement is a point where two lines meet.
A {\bf multiple point} in a line arrangement is a point where more than two lines
meet. The {\bf multiplicity of a point} is the number of lines which meet in the 
point.
\ede

\bde {\bf An arrangement with $t$ non-collinear multiple points} \\
An {\bf arrangement with $t$ non-collinear multiple points} is an arrangement where
there is no line on which there are two multiple points and we can divide it 
into $t$ subsets of lines where  all the lines in each subset intersect in
a single multiple point.
\ede

We denote by $\F ^k$ the free group with $k$ generators.

\subsection{The affine case}

We calculate  the affine case:

\bthm\label{thm_not_cnntd}
Let $\cL$ be a real line arrangement in $\C\PP ^2$ with $t$ 
non-collinear multiple points.
Let $k_i +1$ be the multiplicity of the multiple point $P_i$, $1 \leq i \leq t$. \\
Then:
$$\pi _1 (\C ^2 -\cL) \cong ({\bigoplus _{i=1} ^t} \F ^{k_i}) \oplus \Z^t $$
\ethm

\noindent
{\it Proof:} Randell [Ra] showed that the fundamental group of the complement
of a real line arrangement which consists of $n$ lines meet in a single point
is $\F ^{n-1} \oplus \Z$. 

We can observe $\cL$ as a union of $t$ subsets of lines  $\cL _i, 1 \leq i \leq t$,
 where every such subset $\cL _i,\ 1 \leq i \leq t$, consists of 
$k_i +1$ lines which are passing through the multiple point $P_i$ 
(there is no $l \in \cL_i \cap \cL_j$, because then $l$ connects $P_i$ and $P_j$, 
a contradiction to the assumption). The degree
of each $\cL _i$ is exactly $k_i +1$, because there are $k_i +1$ lines which pass
through the point $P_i$. Moreover, $\cL _i \cap \cL _j =  (k_i +1)(k_j+1)$ points, because
every line in $\cL _i$ meets every line in $\cL _j$.

Every $\cL _i,\ 1 \leq i \leq t$, consists of $k_i +1$ lines which  pass 
through the multiple point $P_i$. This is the configuration of Randell. 
Therefore:
$$\pi _1 (\C ^2 -\cL _i) = \F ^{k_i} \oplus \Z$$

Now we can use the Oka-Sakamoto theorem (see section 2.1), 
in order to compute the fundamental group of the complement of $\cL$:
$$\displaystyle \pi _1 (\C ^2 -\cL) = \pi _1 (\C ^2 -{\bigcup _{i=1} ^t} \cL _i) \stackrel{\rm (O-S)}{\cong} {\bigoplus _{i=1}^t} (\pi _1 (\C ^2 -\cL _i)) =$$
$$ = {\bigoplus _{i=1}^t}( \F ^{k_i} \oplus \Z) =({\bigoplus _{i=1}^t} \F ^{k_i}) \oplus \Z ^t  $$
\hfill $\qed$
  
\bigskip

The Oka-Sakamoto theorem gives us  a new inductive approach to prove
Zariski's proposition:
\bpr {\bf (Zariski)} \\
The fundamental group of the complement of $n$ lines in 
general position is abelian.
\epr

\noindent
{\it Proof:} It is known that for a line $L$:
$$\pi _1(\C ^2 -L) \cong \Z$$
Due to the general position of the lines in the arrangement, we
can use the Oka-Sakamoto theorem (see section 2.1) inductively in the following way:
$$\displaystyle \pi _1(\C ^2 -\cL) = \pi _1 (\C ^2 - {\bigcup _{i=1} ^n} l_i) \stackrel{\rm (O-S)}{\cong} {\bigoplus_{i=1} ^n} (\pi _1 (\C ^2 -l_i)) \cong {\bigoplus_{i=1} ^n} \Z \cong \Z ^n$$
And $\Z^n$ is an abelian group (see [O] too). \hfill $\qed$

\subsection{The projective case}

Now, we will investigate the projective case.

\bthm\label{proj_non_cnntd} 
Let $\cL$ be a real line arrangement in $\C\PP ^2$ with $t$ 
non-collinear multiple points.
Let $k_i +1$ be the multiplicity of the multiple point $P_i$, $1 \leq i \leq t$. \\
Then:
$$\pi _1 (\C\PP ^2 -\cL) \cong ({\bigoplus _{i=1} ^t} \F ^{k_i}) \oplus \Z^{t-1} $$
\ethm

\noindent
{\it Proof:} First, we will prove this theorem for $t=1$, i.e. if $\cL$ 
is a real line arrangement in $\C\PP^2$ which consists of  $k+1$ lines meeting 
in one point $P$, then $\pi _1 (\C\PP ^2 -\cL) \cong \F ^k$.   

Let $\{ \Ga _1, \cdots, \Ga _{k+1} \}$ be a g-base of $\pi _1(\C _u-\cL)$ 
(see section \ref{outline}).
In this line arrangement, we have only one 
singular point - $P$. Therefore, according to lemma \ref{vk-cyc} and
proposition \ref{lml}, this singular point induced the following set 
of relations:
$$[\Ga _{k+1} \Ga _k \cdots \Ga_1, \Ga _i]=1, i=1,\cdots ,k+1$$
Hence, the fundamental group of its affine complement  
has the following presentation:
$$\pi _1(\C ^2 -\cL)\  =\  <\Ga _1, \cdots , \Ga _{k+1} \ | \ [\Ga _{k+1} \Ga _k \cdots \Ga_1, \Ga _i]=1, i=1,\cdots ,k+1>$$ 

We will compute now another presentation for this group.

Let us modify the set of generators $g=\{ \Ga _1, \cdots , \Ga _{k+1} \}$ 
by replacing the generator $\Ga _1$ by the generator
$$\Ga '  =\Ga _{k+1} \Ga _k \cdots \Ga_1$$

Then, we have to check that after the modifications 
we get an equivalent set of generators, 
and we have to calculate the new set of relations. 

\medskip

\bcl
 After replacing $\Ga _1$ by $\Ga '$ 
 (which was defined above) in g, 
we again get  a set of generators. We denote 
this set of generators by $\tilde g$.
\ecl

\noindent
{\it Proof:} we have to show that $\Ga _1 \in <\tilde g>$. But this is
obvious, because: 
$$\Ga _1 = \Ga _2 ^{-1} \Ga _3 ^{-1} \cdots \Ga _{k+1}^{-1} \Ga '$$
\hfill $\qed$
   
\medskip

The next step is the calculation of the new set of relations for  $\tilde g$.

\bcl
The set of relations:
$$\{ [\Ga' ,\Ga]=1 \  |\  \forall \Ga \in \tilde g \}$$
is a complete set of relations for $\tilde g$.
\ecl

\noindent
{\it Proof:} We have to show that 
$$(*) \ \ \{ [\Ga' ,\Ga]=1 \  |\  \forall \Ga \in \tilde g \}$$
is an equivalent set of relations to
$$(**) \ \ \{ [\Ga _{k+1} \Ga _k \cdots \Ga_1, \Ga _i]=1 \ | \  1 \leq i \leq k+1 \}$$
under the assignment: $\Ga ' = \Ga _{k+1} \cdots \Ga _1$.

Let us assume $(*)$. All the relations are equal except the first one. 
We have to show that:
$$[\Ga _{k+1} \cdots \Ga _1, \Ga _1]=1$$ 
But:\\
$$\Ga ' \Ga _1 = \Ga ' (\Ga _2 ^{-1} \cdots \Ga _{k+1} ^{-1} \Ga ') \stackrel{(*)\ + \ ab=ba \Rightarrow ab^{-1}= b^{-1} a}{ =}   (\Ga _2 ^{-1} \cdots \Ga _{k+1} ^{-1} \Ga ') \Ga ' = \Ga _1 \Ga '$$

Now, if we assume $(**)$, all the relations in $(*)$ are equal except of 
$\Ga ' \Ga ' = \Ga ' \Ga '$ which is trivial. \hfill $\qed$

\medskip

Hence we got the following presentation for the fundamental group of the
affine complement of $\cL$:
$$\pi _1 (\C ^2 -\cL) = < \Ga ', \Ga _2, \cdots , \Ga _{k+1} \ | \ [\Ga _i, \Ga '] =1, 2 \leq i \leq k+1 >$$

\medskip

Now, when we are going to the projective case, we add one additional relation, 
according to theorem  \ref{VKP}:
$$\Ga _{k+1} \cdots \Ga _1 = 1$$
In terms of the new generator $\Ga'$, this relation gets 
the following form:
$$\Ga ' =1$$  

Therefore, we can copmute the structure of the fundamental group in the
projective case with $t=1$: \\
$\pi _1(\C\PP ^2 -\cL) = <\Ga ' , \Ga _2 , \cdots, \Ga_{k+1}  \ |\  [\Ga _i, \Ga ']=1, 2 \leq i \leq k+1 ; \Ga '=1> \ \cong \\ \cong \ <\Ga _2, \cdots, \Ga _{k+1}> \oplus <\Ga' \ | \ \Ga '=1> \cong \F ^k$. 

\medskip

Now we continue to the general case ($t>1$). 
For simplicity of the proof, 
we will prove it for two multiple points and the proof for $t$ 
multiple points uses exactly the same arguments.

From the last theorem, we get for a line arrangement $\cL$ 
with two multiple points:
$$\pi _1 (\C ^2 -\cL) \cong \F ^{k_1} \oplus \F ^{k_2} \oplus \Z ^2$$
Let $l_1, \cdots, l_{k_1+1}$ be $k_1+1$ lines which pass through $P_1$ and 
let $l_{k_1+2}, \cdots , l_{k_1+k_2+2}$ be $k_2+1$ lines which  pass 
through $P_2$. We choose  $\{ \Ga _1, \cdots, \Ga _{k_1+k_2+2} \}$, a g-base 
of $\pi _1(\C _u-\cL)$ (see section \ref{outline}) where $\Ga _i$ corresponds 
 to the line $l_i$. 

Similarly to the first part of the proof, we can write the following
presentation for $\pi _1 (\C ^2 -\cL)$:\\
Generators: $g = \{ \Ga _1, \cdots ,\Ga_{k_1}, \Ga ', \Ga_{k_1+2} , \cdots, \Ga_{k_1+k_2+1}, \Ga '' \}$.\\
Relations: $\cR = \{ \Ga _i \Ga_j =\Ga _j \Ga _i, 1 \leq i \leq k_1, k_1+2 \leq j \leq k_1+k_2+1 ;$\\
$[\Ga' , \Ga ] = 1 , \forall \Ga \in g ;  \  [\Ga'' , \Ga ] = 1 , \forall \Ga \in g \}$,
where:
$$\Ga '= \Ga _{k_1+1} \cdots \Ga _1;\ \  \Ga ''= \Ga _{k_1+k_2+2} \cdots \Ga _{k_1+2}$$

Now, when we are going to the projective case, we add one additional relation, 
according to theorem  \ref{VKP}:
$$\Ga _{k_1+k_2+2} \cdots \Ga _1 = 1$$
In terms of the new generators $\Ga',\Ga ''$, this relation gets 
the following form:
$$\Ga '' \Ga ' =1$$  

Now, we can finish to compute the structure the fundamental group 
in the projective case: \\
$\pi _1(\C\PP ^2 -\cL) = < g \ | \ \cR, \Ga '' \Ga '=1> \ \cong \\ \cong \ <\Ga _1, \cdots, \Ga _{k_1}> \oplus <\Ga_{k_1+2}, \cdots, \Ga _{k_1+k_2+1}> \oplus <\Ga', \Ga '' \ | \ \Ga '' \Ga '=1>  \cong \F ^{k_1} \oplus \F ^{k_2} \oplus \Z$. \hfill $\qed$

\bigskip
 
As a consequence of the last theorem, we get:

\bco
$$\pi _1 (\C ^2 -\cL) \cong \pi _1 (\C\PP ^2 -\cL)\oplus \Z $$
\eco

Therefore, the short exact sequence which was proved by Oka (theorem \ref{oka}):
$$ 1 \to \Z \to \pi _1 (\C ^2 -\cL) \to \pi _1 (\C\PP ^2 -\cL) \to 1 $$  
splits.

\section{Arrangements with $t$ collinear multiple points}

In this section, we are going to calculate the fundamental group
of the complement of line arrangements which consist of $t$ subsets of lines
where all the lines in each subset intersect at a single (multiple) point, 
all the $t$ multiple intersection points lie on a single line which belongs
to all the subsets and any two subsets of lines intersect in that line
and in simple points out of that line. We define:

\bde {\bf An arrangement with $t$ collinear multiple points}\\
An {\bf  arrangement with $t$ collinear multiple points} is a line arrangement which 
contains a line where all the $t$ multiple points lie on it.
\ede

\subsection{The affine case}

\bthm\label{thm_t_cnntd}
Let $\cL$ be a real line arrangement in $\C\PP ^2$
with $t$ collinear multiple points $P_1, \cdots, P_t$ with multiplicities 
$k_1+1, \cdots, k_t+1$, respectively. Then:
$$\pi _1 (\C ^2 -\cL) \cong {\bigoplus _{i=1} ^t} \F ^{k_i} \oplus \Z $$
\ethm

\medskip
 
It has to be noted that this theorem has a similar result to what we have got
in the previous section in the non-collinear case. In both cases, the multiple
points induced the free groups. The difference between the cases is that the
connected line of the collinear case degenerates all the infinite cyclic groups 
of the non-collinear case into one infinite cyclic group.  

\bigskip

Let $L$ be the line on which all the multiple points lie. 
We choose  $\{ \Ga _1, \cdots, \Ga _n \}$ $(n=\# \{ l \in \cL\})$, a g-base 
of $\pi _1(\C _{u_0}-\cL)$ (see section \ref{outline}), where $\Ga _i$  
corresponds to the line $l_i$ in $\cL$.
The proof of the theorem is based on the following two lemmas:

\ble\label{lem1}
In the situation of the theorem,
let $\cL _i$ be the subset of lines meet in $P_i$ apart from $L$.
Then: $[\Ga _i,\Ga _j]=1$ where $l_i \in \cL _i, l_j \in \cL _j$ and
$1 \leq i < j \leq t$.
\ele

\ble\label{lem2}
Let $\cL _i \cup L = \{ l_{p_1}, \cdots, l_{p_{k_i+1}} \}$ be the $k_i+1$ lines that meet
in the multiple point $P_i$. Then, the relations that
are induced from this multiple point are:
$$[ \Ga _{p_{k_i+1}} \cdots \Ga _{p_1}, \Ga _{p_j}]=1, \quad 1 \leq j \leq k_i+1$$
\ele

The proof of lemma \ref{lem1} is in section \ref{lem1-section}. 
The proof of lemma \ref{lem2} is in section \ref{lem2-section}. 
The proof of the theorem (\ref{thm_t_cnntd}) is in section \ref{lem3-section}.

\subsection{Proof of lemma \ref{lem1}}\label{lem1-section}

For simplicity, we prove the lemma only for two multiple points, 
and the proof for $t$ multiple points uses exactly the same arguments.

\medskip

We will split the proof of this lemma into two cases: with the restriction 
that all the simple intersection points are to the right of the multiple 
points, and without this restriction. This restriction simplifies the proof 
significantly, and help to understand the proof of the general case.

\subsubsection{First case - with the restriction}
In this case, all the simple points are to the right of the multiple points.

\medskip

Let $N = \{ x \in \C \ | \ (x,y)\ {\rm is\ an\ intersection\ point} \}$, and
let $u_0 \in \R$ such that $x \ll u_0$ for all $x \in N$.
Let $\C _{u_0} = \{ (u_0,y) \ | \ y \in \C \}$. 
We numerate the lines according to their intersection with $\C _{u_0}$. 
By a proper choosing of the line in infinity and homotopic movements 
of the lines, we can assume that the line arrangement has the 
following property: for $1 \leq i < j \leq k_1$,   
$$x(l_i \cap l_t) < x (l_j \cap l_s), \ \ k_1 +1 \leq t,s \leq k_1+k_2$$
Therefore, we get the following line arrangement:

\medskip

\begin{center}
\epsfxsize=12cm
\epsfbox{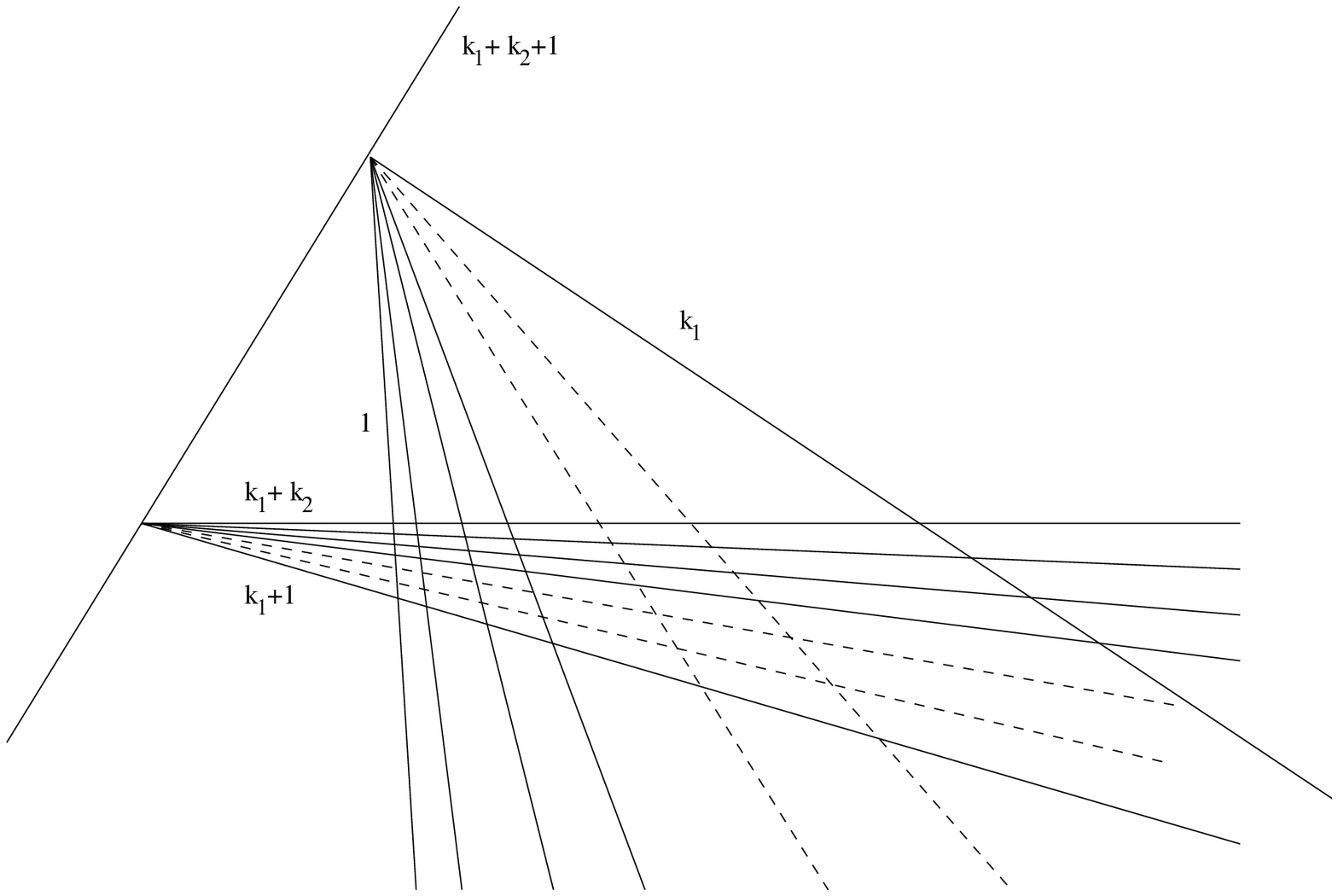}
\end{center}

\medskip

Let $g=\{ \Ga _1, \cdots, \Ga _{k_1+k_2+1} \}$ be a g-base of $\pi _1 (\C _{u_0} - \cL)$. 
By abuse of notations, let us 
denote the images of $\Ga _i$ in $\pi _1 (\C ^2 - \cL)$ by the same notation.
   
Now, we prove this lemma using the braid monodromy techniques (\ref{MT})
and the Van-Kampen theorem (\ref{VK2}). First, let us calculate 
the skeletons representing the $\cL$.V.C.s of the braid monodromy.

\medskip

According to this line arrangement, we have the following set of 
Lefschetz pairs:
 
\begin{center}
\begin{tabular}{|c|c|}
\hline
 $j$ & $\la _{x_j}$ \\
\hline
$1$ & $(k_1,k_1+1)$ \\
$2$ & $(k_1+1,k_1+2)$ \\
$3$ & $(k_1+2,k_1+3)$ \\
\vdots & \vdots \\
$k_2$ & $(k_1+k_2-1,k_1+k_2)$ \\
$k_2+1$ & $(k_1-1,k_1)$ \\
$k_2+2$ & $(k_1,k_1+1)$ \\
\vdots & \vdots \\
$2 k_2$ & $(k_1+k_2-2,k_1+k_2-1)$ \\
\vdots & \vdots \\
$(k_1-1)k_2+1$ & $(1,2)$ \\
$(k_1-1)k_2+2$ & $(2,3)$ \\
\vdots & \vdots \\
 $k_1 k_2$ & $(k_2,k_2 +1)$ \\
 $k_1 k_2 +1$ & $(k_2+1,k_1 +k_2+1)$ \\
 $k_1 k_2 +2$ & $(1, k_2+1)$ \\
\hline
\end{tabular}
\end{center}

Let $\{ \de _i \ | \ 1 \leq i \leq k_1 k_2 +2 \}$ be a g-base for $\pi _1(\C ^X-N,u_0)$ (where
$\C ^X$ is the $x$-axis). Let $\varphi$ be the braid monodromy of $\cL$ w.r.t. 
$\pi _1,u_0$.

Now, using the table of Lefschetz pairs, we can calculate the skeletons 
representing the $\cL$.V.C.s for the braids $\varphi (\de _i)$ 
(according to Moishezon-Teicher's algorithm (\ref{MT})). 
Here, we will calculate  the $\cL$.V.C.s of the two general cases.

\medskip

{\bf Skeleton representing the $\cL$.V.C. of $\varphi (\de _{l k_2+1}), 0 \leq l \leq k_1-1$}: 
The Lefschetz pair is 
 $(k_1 -l,k_1-l+1)$. So the skeleton representing the local $\cL$.V.C. is:

\medskip
\begin{center}
\epsfxsize=11cm
\epsfbox{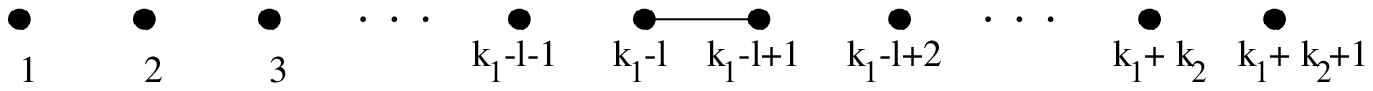}  
\end{center}

\medskip

\noindent
According to the algorithm, we have to apply on the skeleton the composition 
of the following $l$ sequences of braids:  
$$\De <k_1+k_2-l ,k_1+k_2-l+1> \De <k_1 +k_2-l-1,k_1+k_2-l> \cdots $$
$$ \De <k_1-l+2,k_1 -l+3> \De <k_1-l+1,k_1 -l+2>$$
$$\De <k_1+k_2-l+1 ,k_1+k_2-l+2> \De <k_1 +k_2-l,k_1+k_2-l+1> \cdots $$
$$ \De <k_1-l+3,k_1 -l+4> \De <k_1-l+2,k_1 -l+3>$$
$$\vdots$$
$$\De <k_1+k_2-1 ,k_1+k_2> \De <k_1 +k_2-2,k_1+k_2-1> \cdots \De <k_1,k_1 +1>$$
In every sequence, only the last braid of the sequence affects 
the skeleton (because the region of the others has no intersection with 
the region of the skeleton). Therefore, we get the following skeleton:   

\medskip

\begin{center}
\epsfxsize=13cm
\epsfbox{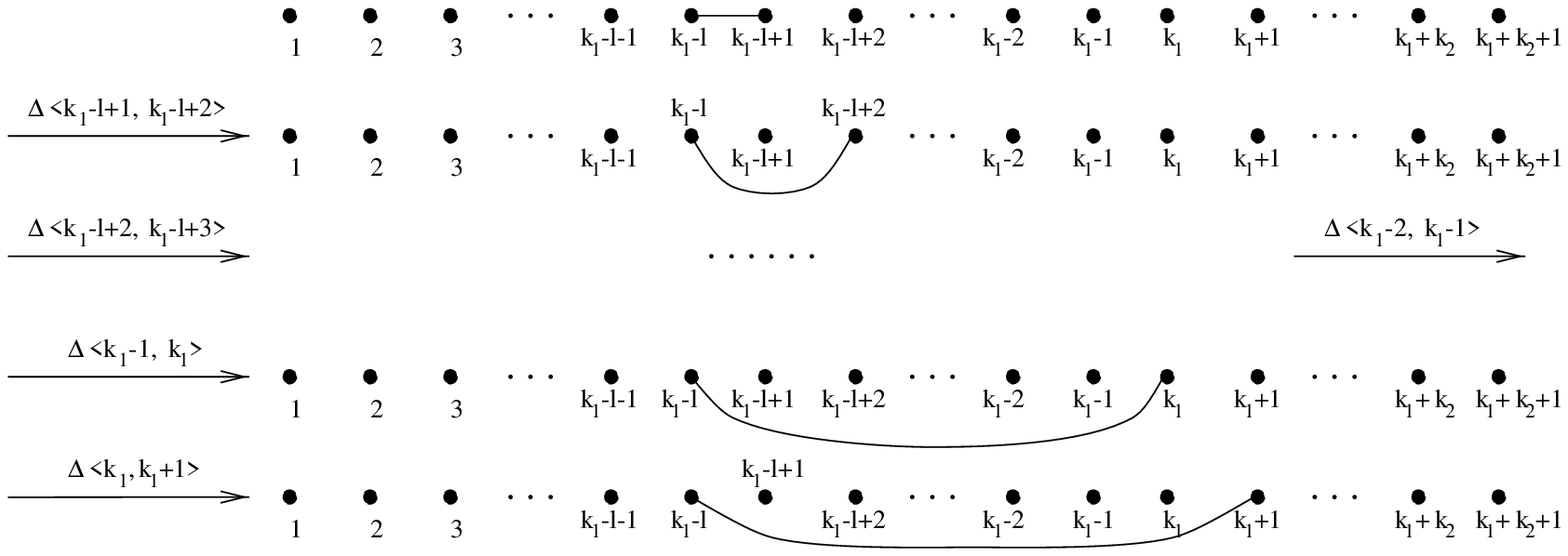}  
\end{center}

\medskip

{\bf Skeleton representing the $\cL$.V.C. of $\varphi (\de _{l k_2+i}), 0 \leq l \leq k_1-1, 2 \leq i \leq k_2$}: 
The Lefschetz pair is $(k_1 -l+i-1,k_1-l+i)$. So the skeleton representing local 
$\cL$.V.C. is:

\medskip

\begin{center}
\epsfxsize=11cm
\epsfbox{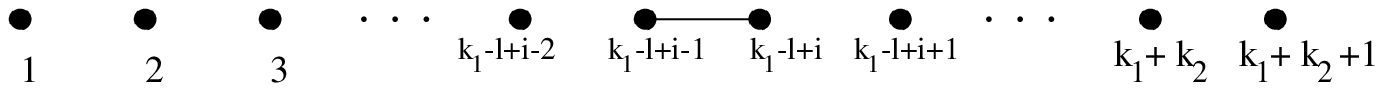}  
\end{center}

\medskip

\noindent
According to the algorithm, we have to apply on the skeleton the composition 
of the following $l+1$ sequences of braids:
$$\De <k_1 -l+i-2,k_1-l+i-1> \De < k_1 -l+i-3,k_1-l+i-2> \cdots$$
$$\De < k_1 -l+1,k_1-l+2 > \De < k_1 -l,k_1-l+1 > $$  
$$\De <k_1+k_2-l ,k_1+k_2-l+1> \De <k_1 +k_2-l-1,k_1+k_2-l> \cdots$$ 
$$ \De <k_1-l+2,k_1 -l+3> \De <k_1-l+1,k_1 -l+2>$$
$$\De <k_1+k_2-l+1 ,k_1+k_2-l+2> \De <k_1 +k_2-l,k_1+k_2-l+1> \cdots$$
$$\De <k_1-l+3,k_1 -l+4> \De <k_1-l+2,k_1 -l+3>$$
$$\vdots$$
$$\De <k_1+k_2-1 ,k_1+k_2> \De <k_1 +k_2-2,k_1+k_2-1> \cdots \De <k_1,k_1 +1>$$
The first sequence causes the following effect to the skeleton:
 
\medskip

\begin{center}
\epsfxsize=14cm
\epsfbox{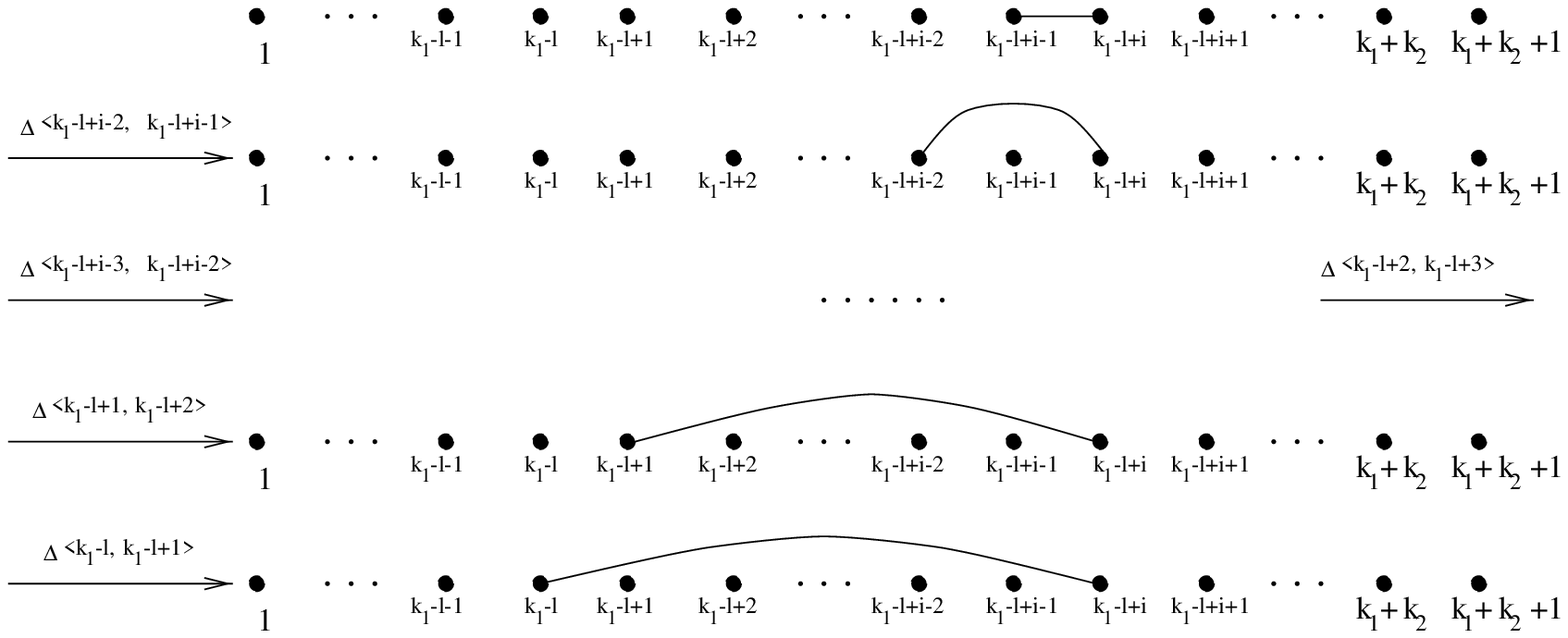}  
\end{center}

\medskip

\noindent
Only the last part of the second sequence affects the skeleton as follows:

\medskip

\begin{center}
\epsfxsize=14cm
\epsfbox{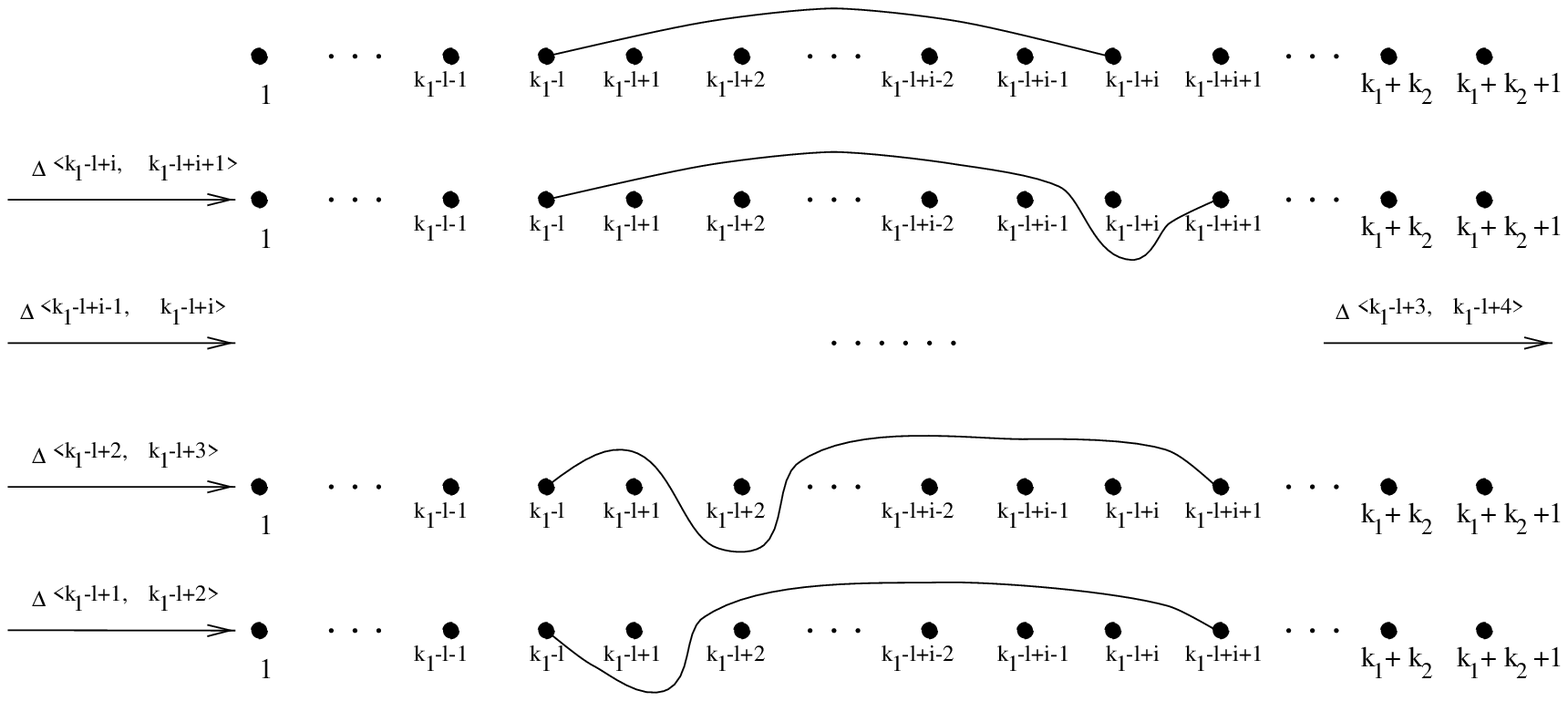}  
\end{center}

\medskip

\noindent
In the other $l-1$ sequences of braids, only the second part of the sequence
 affects, i.e. only  
the braids whose region intersects the region of the skeleton.
Therefore, we get the following skeleton representing the $\cL$.V.C.:
  
\medskip

\begin{center}
\epsfxsize=14cm
\epsfbox{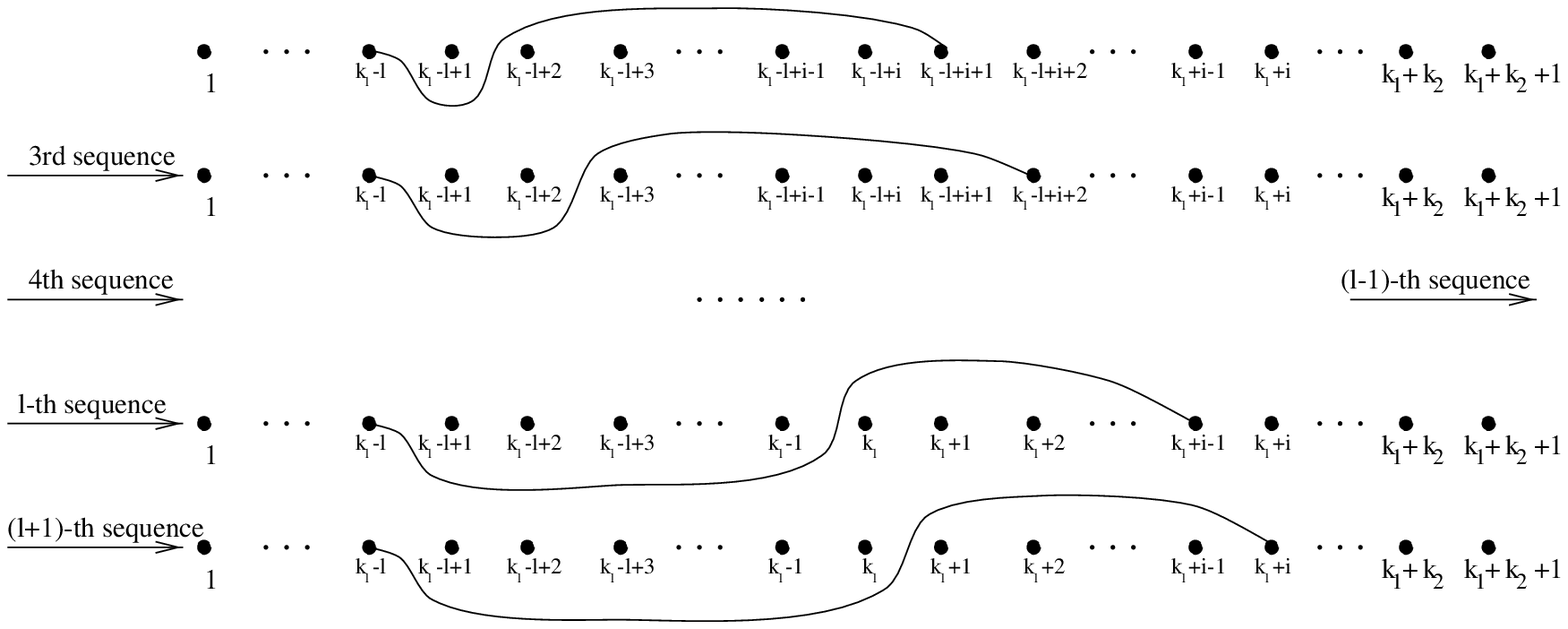}  
\end{center}

\medskip

After we have calculated the skeletons representing $\cL$.V.C.s for the 
braid monodromy, we can calculate the relations that they induced. 
As we have introduced in the previous section, according to 
Van-Kampen theorem (\ref{VK2}), every $\cL$.V.C. induces a relation. Now, we will
calculate the general relations which are induced from the general $\cL$.V.C.s .

\medskip

The relation which is induced from $\varphi (\de _{lk_2+1}),\ 0 \leq l \leq k_1-1$:

\medskip

\begin{center}
\epsfxsize=13cm
\epsfbox{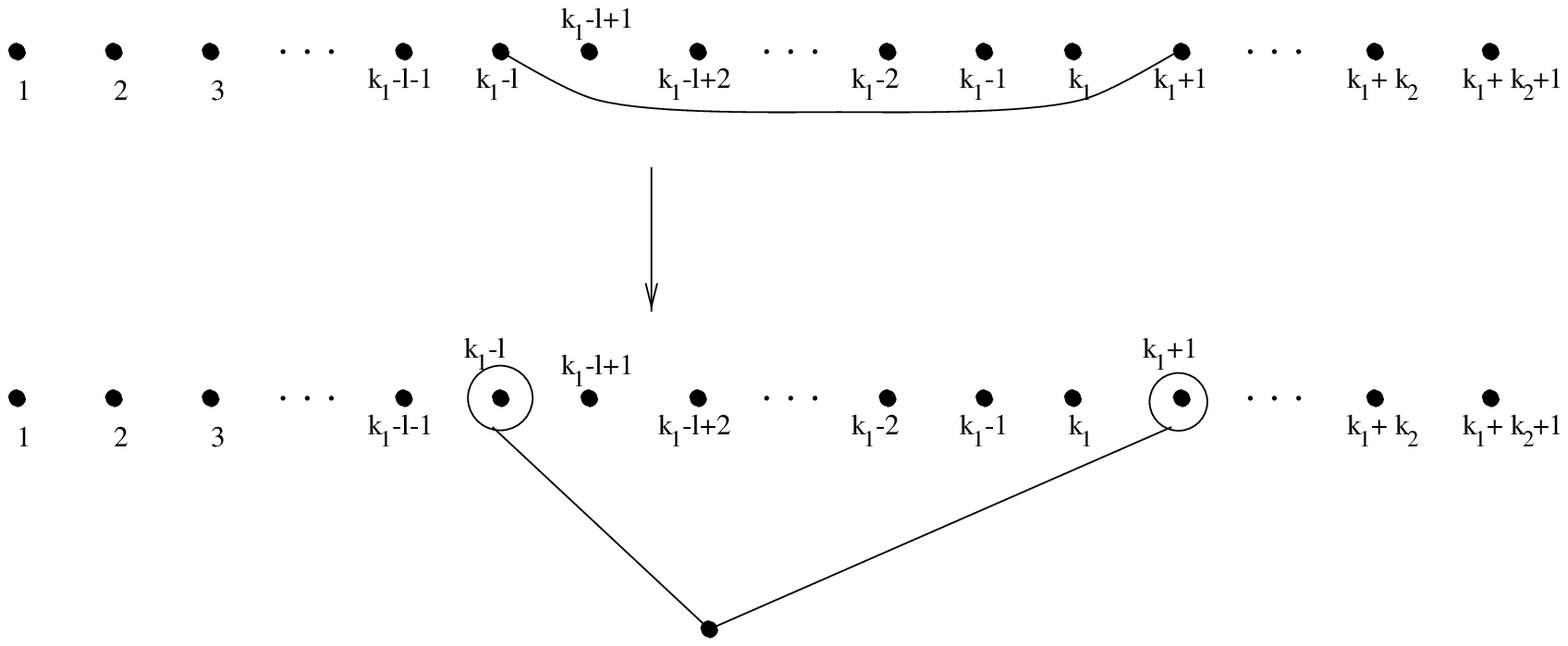}  
\end{center}

\medskip

\noindent
Therefore, the relation is: 
$$\Ga _{k_1-l} \Ga _{k_1+1} = \Ga _{k_1+1} \Ga _{k_1-l}$$

The relation which is induced from $\varphi (\de _{lk_2+i}),\ 0 \leq l \leq k_1-1,\ 2 \leq i \leq k_2$:

\medskip

\begin{center}
\epsfxsize=13cm
\epsfbox{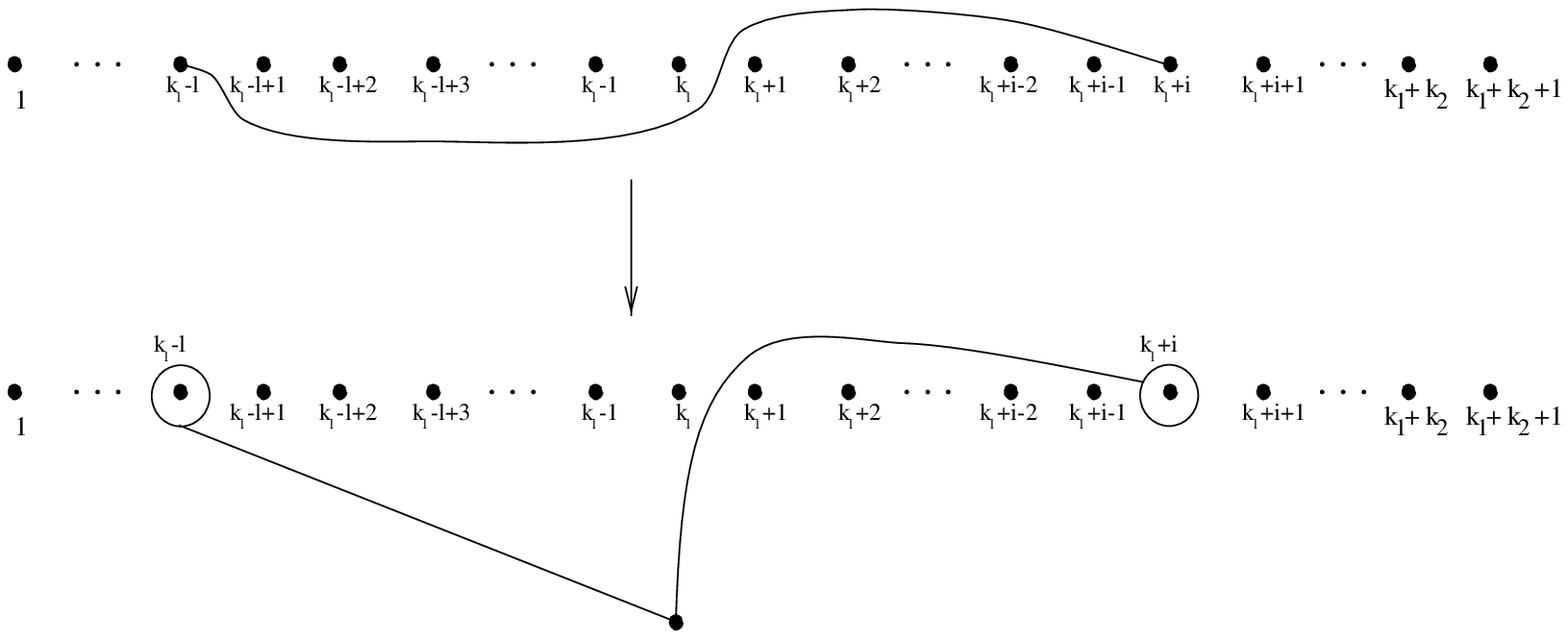}  
\end{center}

\medskip

\noindent
Therefore, the relation is:
$$\Ga _{k_1-l} \Ga ^{-1} _{k_1+1} \cdots \Ga ^{-1} _{k_1+i-1} \Ga _{k_1+i} \Ga _{k_1+i-1} \cdots \Ga _{k_1+1} =$$
$$\Ga ^{-1} _{k_1+1} \cdots \Ga ^{-1} _{k_1+i-1} \Ga _{k_1+i} \Ga _{k_1+i-1} \cdots \Ga _{k_1+1} \Ga _{k_1-l}$$

\medskip

Therefore, we got the following set of relations:\\
 for all $\ 0 \leq l \leq k_1-1,\ 1 \leq i \leq k_2,$
$$\Ga _{k_1-l} \Ga ^{-1} _{k_1+1} \cdots \Ga ^{-1} _{k_1+i-1} \Ga _{k_1+i} \Ga _{k_1+i-1} \cdots \Ga _{k_1+1} =$$
$$\Ga ^{-1} _{k_1+1} \cdots \Ga ^{-1} _{k_1+i-1} \Ga _{k_1+i} \Ga _{k_1+i-1} \cdots \Ga _{k_1+1} \Ga _{k_1-l}$$

Now, it is easy to see that this set of relations is equivalent to 
the following set of relations (see [Ga]):
$$ \Ga _i \Ga _j = \Ga _j \Ga _i; \ 1 \leq i \leq k_1, \ k_1+1 \leq j \leq k_1+k_2$$
and this finished the proof of the first case of the first lemma 
(\ref{lem1}). \hfill $\qed$

\subsubsection{Second case - without the restriction}

Let $N = \{ x \in \C \ | \ (x,y)\ {\rm is\ an\ intersection\ point} \}$, and
let $u_0 \in \R$ such that $x \ll u_0$ for all $x \in N$.
Let $\C _{u_0} = \{ (u_0,y) \ | \ y \in \C \}$. 
We numerate the lines according to their intersection with $\C _{u_0}$.
We organized this line arrangement in such a way that the following
property holds:\\ 
for $1 \leq i < j \leq l$ and $k_1+l+1 \leq i < j \leq k_1+k_2$, 
$$x(L_i \cap L_t) < x(L_j \cap L_s), \ \  l+1 \leq s,t \leq k_1+l$$
It is easy to see that this is the general case, i.e. every line arrangement
is homotopic  to this situation by rotations and a proper choosing of the
line at infinity.

\medskip
 
Therefore, we get the following line arrangement:

\medskip

\begin{center}
\epsfxsize=12cm
\epsfbox{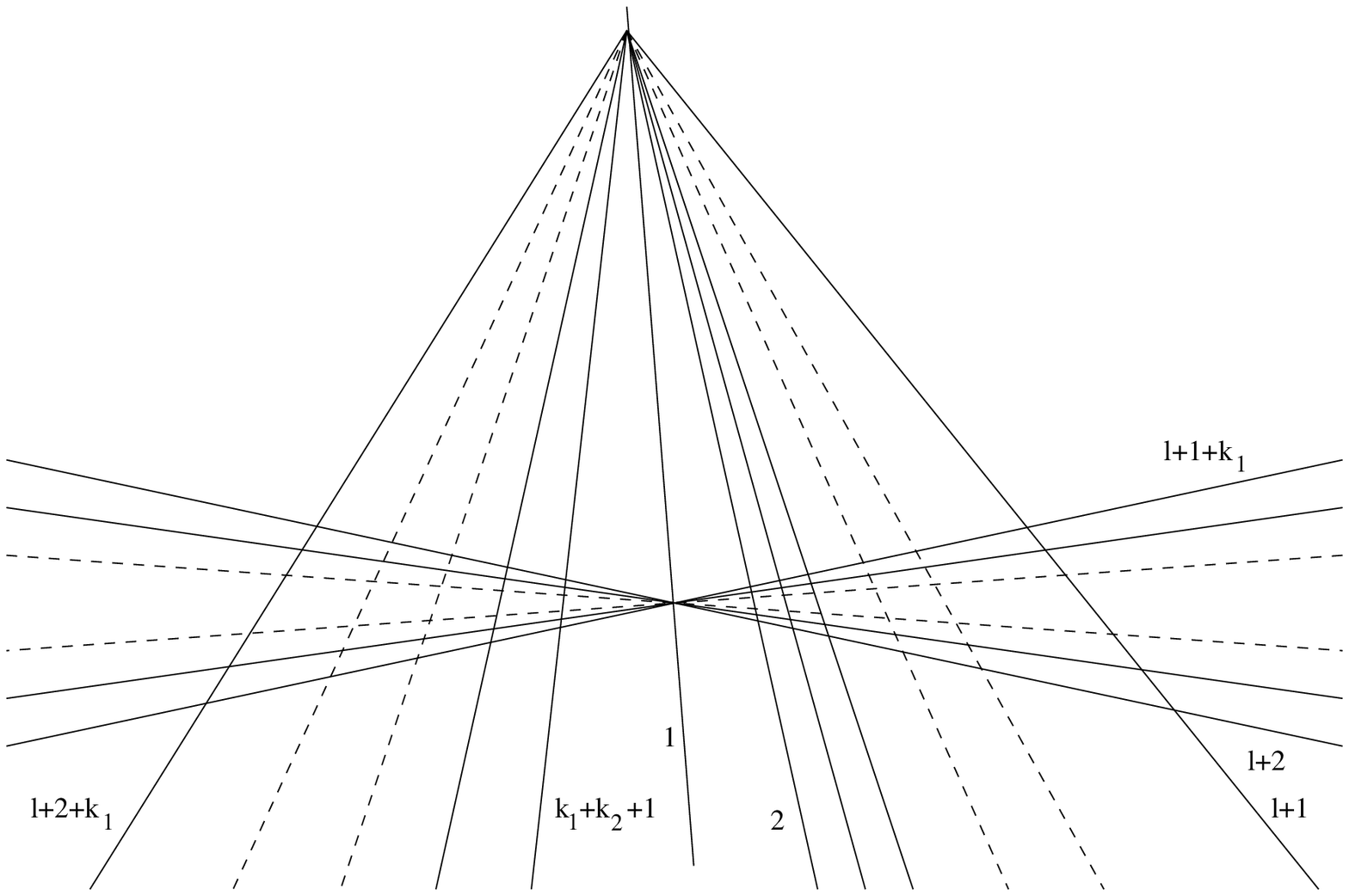}
\end{center}

\medskip

Let $g= \{ \Ga _1, \cdots , \Ga _{k_1+k_2+1} \}$ be a g-base of
 $\pi _1 (\C _{u_0} - \cL)$. 
By abuse of notations, let us
denote the images of $\Ga _i$ in $\pi _1 (\C ^2 - \cL)$ by the same notation.
   
Now, we prove this lemma using the braid monodromy techniques (\ref{MT})
and the Van-Kampen theorem (\ref{VK2}). First, let us calculate 
the skeletons representing the $\cL$.V.C.s of the braid monodromy.

\medskip

According to this line arrangement, we have the following set of 
Lefschetz pairs:

\begin{center}
\begin{tabular}{|c|c|}
\hline
 $j$ & $\la _{x_j}$ \\
\hline
 $1$ & $(l+1,l+2)$ \\
 $2$ & $(l+2,l+3)$ \\
\vdots & \vdots \\
 $k_1$ & $(k_1+l,k_1+l+1)$ \\
 $k_1+1$ & $(l,l+1)$ \\
 $k_1+2$ & $(l+1,l+2)$ \\
 \vdots & \vdots \\
 $2 k_1$ & $(k_1+l-1,k_1+l)$ \\
 \vdots & \vdots \\ 
 $(l-1) k_1 +1$ & $(2,3)$ \\
 $(l-1) k_1+2$ & $(3,4)$ \\
 \vdots & \vdots \\
 $l k_1$ & $(k_1+1,k_1+2)$ \\
 $l k_1+1$ & $(1,k_1+1)$ \\
 $l k_1+2$ & $(k_1+1,k_1+k_2+1)$ \\
 $(l k_1+2)+1$ & $(k_1,k_1+1)$ \\
 $(l k_1+2)+2$ & $(k_1-1,k_1)$ \\
 \vdots & \vdots \\
 $(l k_1+2)+k_1$ & $(1,2)$ \\
 $(l k_1+2)+k_1+1$ & $(k_1+1,k_1+2)$ \\
 $(l k_1+2)+k_1+2$ & $(k_1,k_1+1)$ \\
 \vdots & \vdots \\
 $(l k_1+2)+2 k_1$ & $(2,3)$ \\
 \vdots & \vdots \\
 $(l k_1+2)+(k_2-l-1)k_1+1$ & $(k_1+k_2-l-1,k_1+k_2-l)$ \\
 $(l k_1+2)+(k_2-l-1)k_1+2$ & $(k_1+k_2-l-2,k_1+k_2-l-1)$ \\
\vdots & \vdots \\
 $(l k_1+2)+(k_2-l)k_1 \ [=k_1 k_2 +2]$ & $(k_2-l, k_2-l+1)$ \\
\hline
\end{tabular}
\end{center}

Let $\{ \de _i \ | \ 1 \leq i \leq k_1 k_2 +2 \}$ be a g-base for $\pi _1(\C ^X-N,u_0)$ (where
$\C ^X$ is the $x$-axis). Let $\varphi$ be the braid monodromy of $\cL$ w.r.t. 
$\pi _1,u_0$.

Now, using the table of the Lefschetz pairs, we can calculate 
the skeletons representing $\cL$.V.C.s for
the braids $\varphi (\de _i)$ 
(according to the Moishezon-Teicher algorithm (\ref{MT})).

\medskip

Until singular point number $lk_1$ we have almost 
the same configuration 
as in the first case of the lemma, hence the general skeleton, 
which represents the $\cL$.V.C., which we have found there is identical 
(but its center is shifted one point left) 
to the general skeleton in this case of the lemma until point number $lk_1$. 
Therefore:

\medskip

{\bf Skeleton representing the $\cL$.V.C. of $\varphi (\de _{i k_1+1}), 0 \leq i \leq l-1$}:

\medskip

\begin{center}
\epsfxsize=12cm
\epsfbox{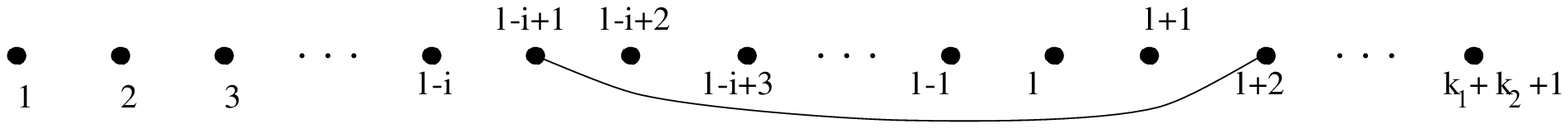}  
\end{center}

\medskip

{\bf Skeleton representing the $\cL$.V.C. of $\varphi (\de _{i k_1+j}), \ 0 \leq i \leq l-1, \ 2 \leq j \leq k_1$}: 
  
\medskip

\begin{center}
\epsfxsize=12cm
\epsfbox{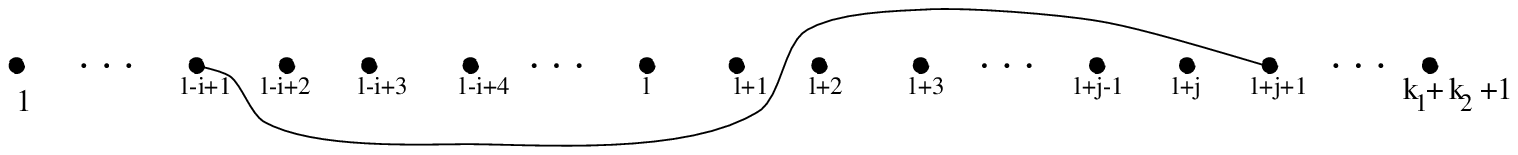}  
\end{center}

\medskip

We skip the calculations of the braid monodromy of the two multiple points
(which will be done in the proof of the next lemma (\ref{lem2})),
and we continue with the rest of the simple points and we pass directly 
to the general case:

\medskip

{\bf Skeleton representing the $\cL$.V.C. of $\varphi (\de _{(l k_1+2)+i k_1+1}), 0 \leq i \leq (k_2-l-1)$}: 
The Lefschetz pair is
 $$(k_1+i,k_1+i+1),$$
 therefore the skeleton representing the local $\cL$.V.C. is:

\medskip

\begin{center}
\epsfxsize=12cm
\epsfbox{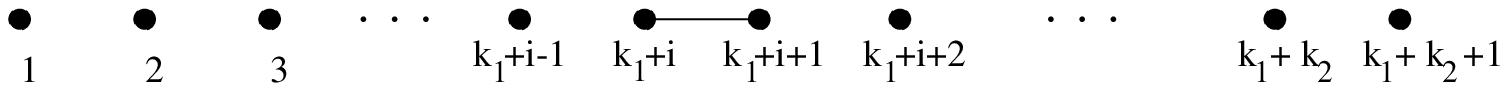}  
\end{center}

\medskip

\noindent
We have to apply on this skeleton the following sequences of braids:
$$\De <i,i+1> \De <i+1,i+2> \cdots \De <k_1+i-1,k_1+i>$$
$$\vdots$$
$$\De <1,2> \De <2,3> \cdots \De <k_1,k_1+1>$$
$$\De <k_1+1,k_1+k_2+1> \De <1,k_1+1>$$
$$\De <k_1+1,k_1+2> \De <k_1,k_1+1> \cdots \De <2,3>$$
$$\vdots$$
$$\De <k_1+l,k_1+l+1> \De <k_1+l-1,k_1+l> \cdots \De <l+1,l+2>$$

\noindent
In the first $i-1$ sequences, only the last braid in each sequence affects
 the skeleton, hence we get:

\medskip

\begin{center}
\epsfxsize=14cm
\epsfbox{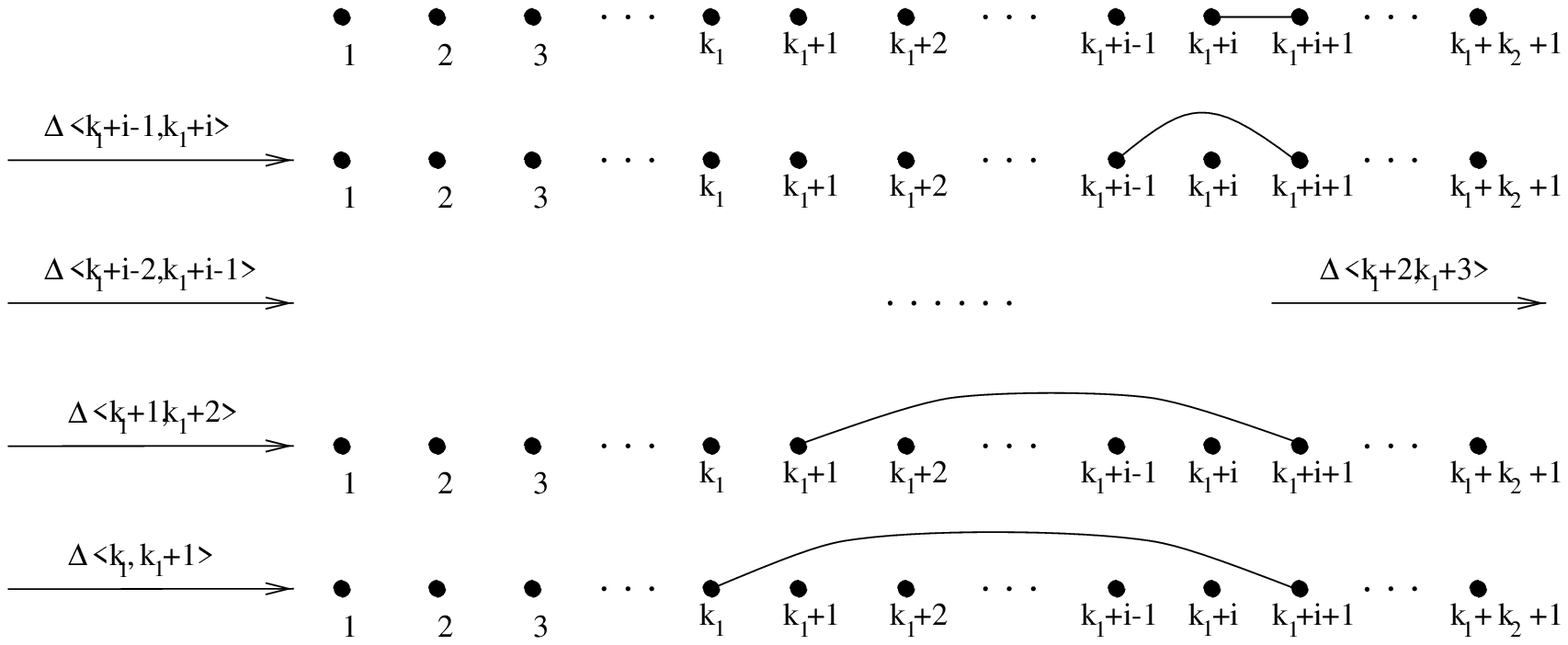}  
\end{center}

\medskip

\noindent 
Next, the action of the braids $\De <k_1+1,k_1+k_2+1>$ and $\De <1,k_1+1>$ 
is as follows:

\medskip

\begin{center}
\epsfxsize=14cm
\epsfbox{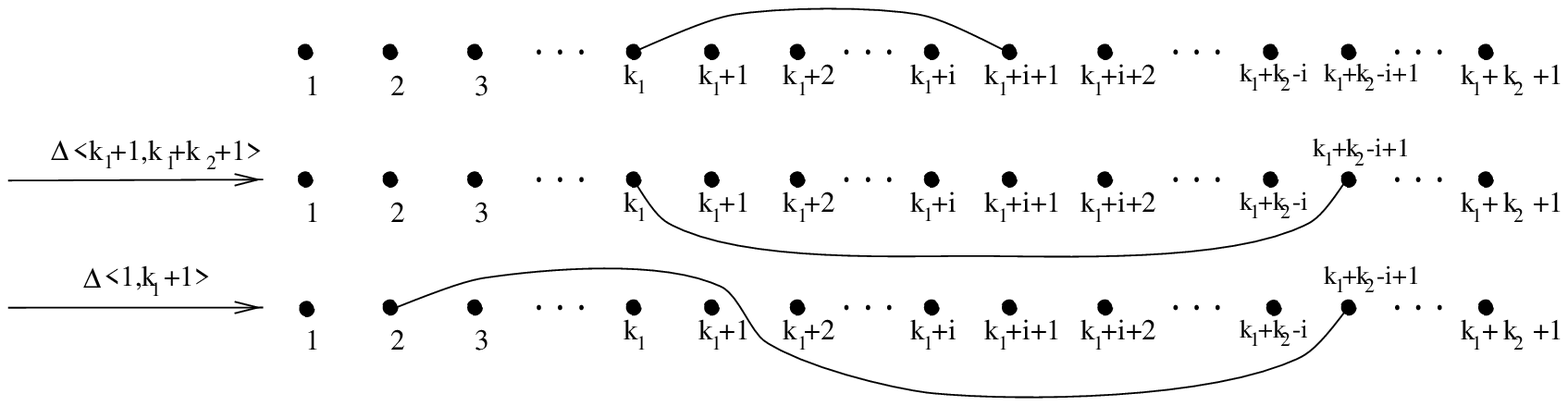}  
\end{center}

\medskip

\noindent
Then, the $l$ sequences of braids move the leftest side of the skeleton
$l$ points right:

\medskip

\begin{center}
\epsfxsize=14cm
\epsfbox{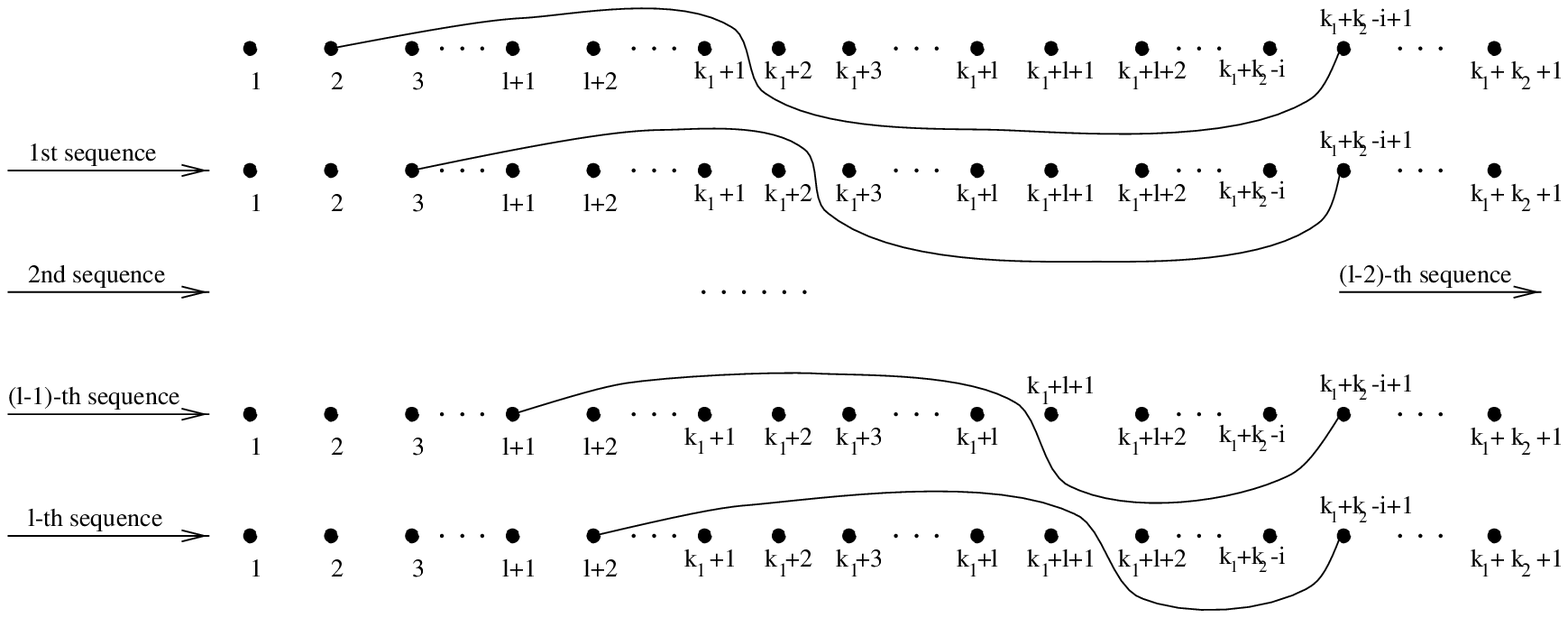}  
\end{center}

\medskip

{\bf Skeleton representing the $\cL$.V.C. of $\varphi (\de _{(l k_1+2)+i k_1+j});\ 0 \leq i \leq k_2-l-1,\\ 2 \leq j \leq k_1$}: 
The Lefschetz pair is 
 $$(k_1+i-j+1,k_1+i-j+2)$$
 therefore the skeleton representing the local $\cL$.V.C. is:

\medskip

\begin{center}
\epsfxsize=10.5cm
\epsfbox{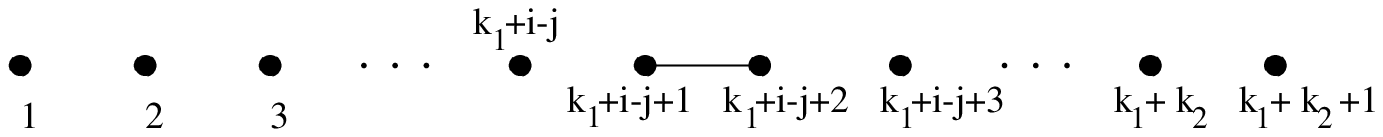}  
\end{center}

\medskip

\noindent
We have to apply on this skeleton the following sequences of braids:
$$\De <k_1+i-j+2,k_1+i-j+3> \De <k_1+i-j+3,k_1+i-j+4> \cdots$$
$$\De <k_1+i-1,k_1+i> \De <k_1+i,k_1+i+1>$$
$$\De <i,i+1> \De <i+1,i+2> \cdots \De <k_1+i-1,k_1+i>$$
$$\vdots$$
$$\De <1,2> \De <2,3> \cdots \De <k_1,k_1+1>$$
$$\De <k_1+1,k_1+k_2+1> \De <1,k_1+1>$$
$$\De <k_1+1,k_1+2> \De <k_1,k_1+1> \cdots \De <2,3>$$
$$\vdots$$
$$\De <k_1+l,k_1+l+1> \De <k_1+l-1,k_1+l> \cdots \De <l+1,l+2>$$
 
\noindent
The first sequence acts as follows:

\medskip

\begin{center}
\epsfxsize=14cm
\epsfbox{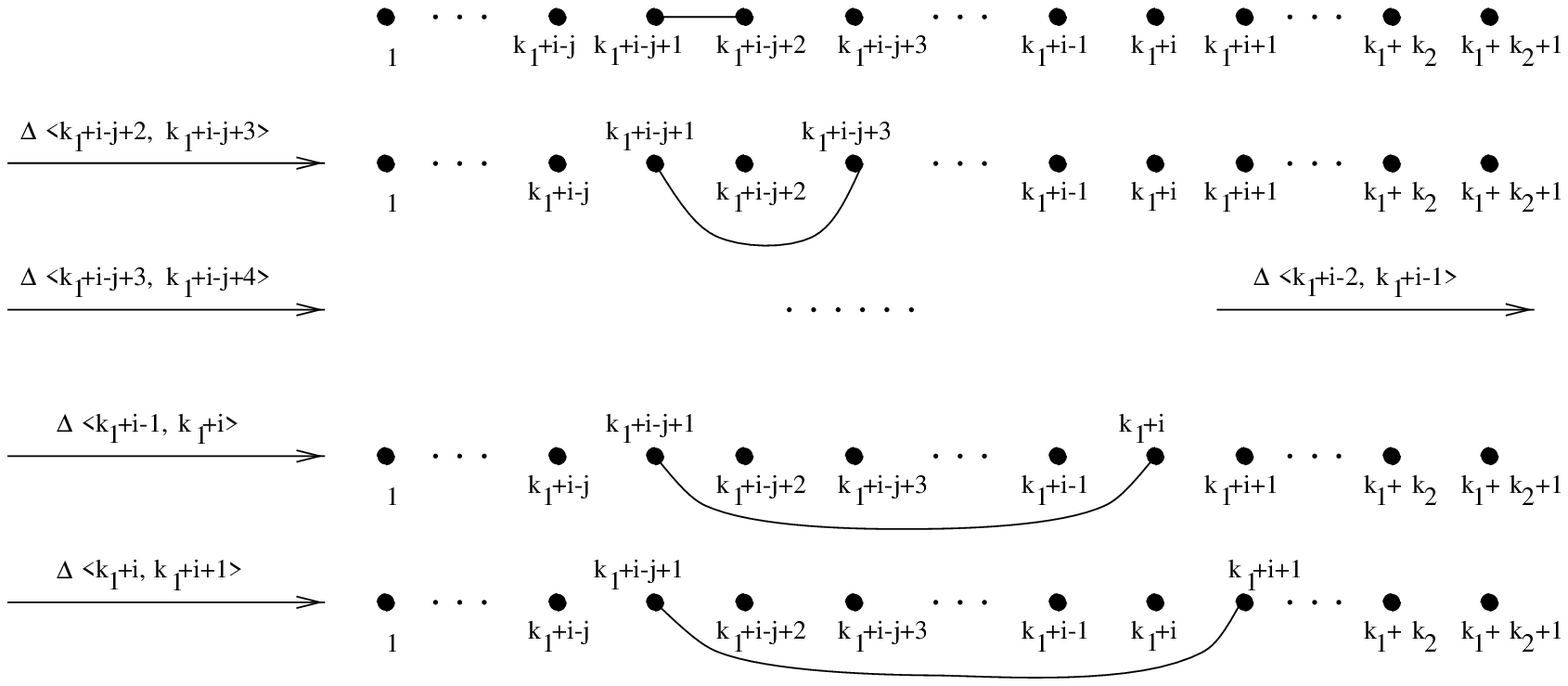}  
\end{center}

\medskip

\noindent
The second sequence moves the left side of the skeleton one point left
(the first part of the sequence does not affect the skeleton):

\medskip

\begin{center}
\epsfxsize=14cm
\epsfbox{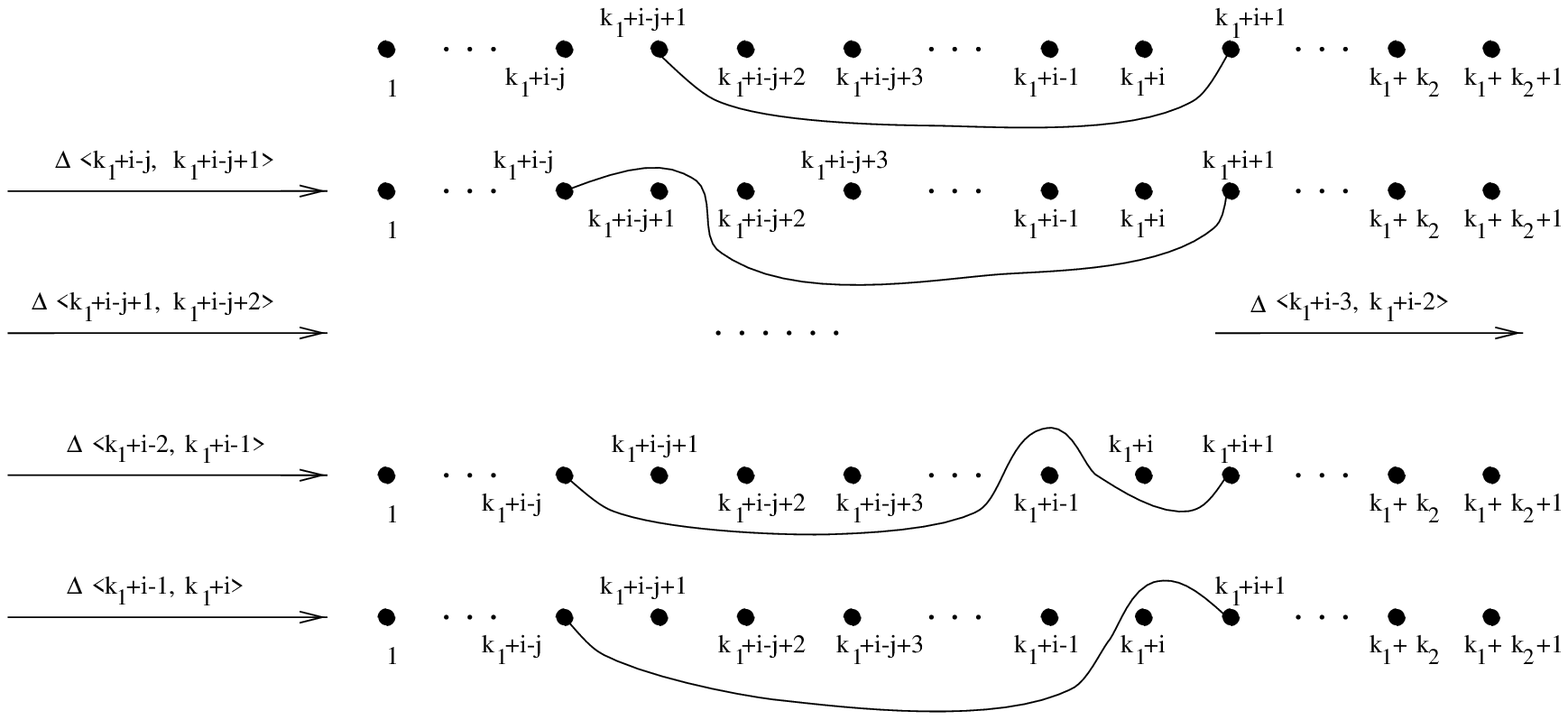}  
\end{center}

\medskip

\noindent
Each of the next $i-1$ sequences moves the left side of the skeleton another
step left, so we get the following: 

\medskip

\begin{center}
\epsfxsize=14cm
\epsfbox{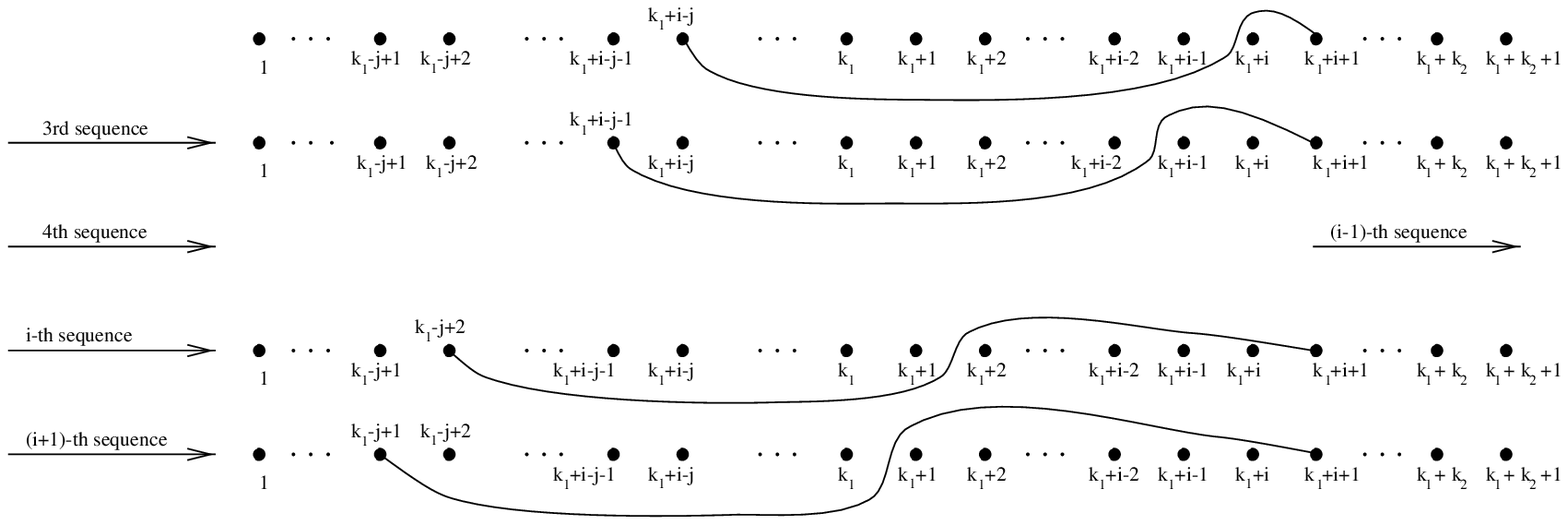}  
\end{center}

\medskip
 
\noindent
Next, the action of the braids $\De <k_1+1,k_1+k_2+1>$ and $\De <1,k_1+1>$ 
is as follows:

\medskip

\begin{center}
\epsfxsize=14cm
\epsfbox{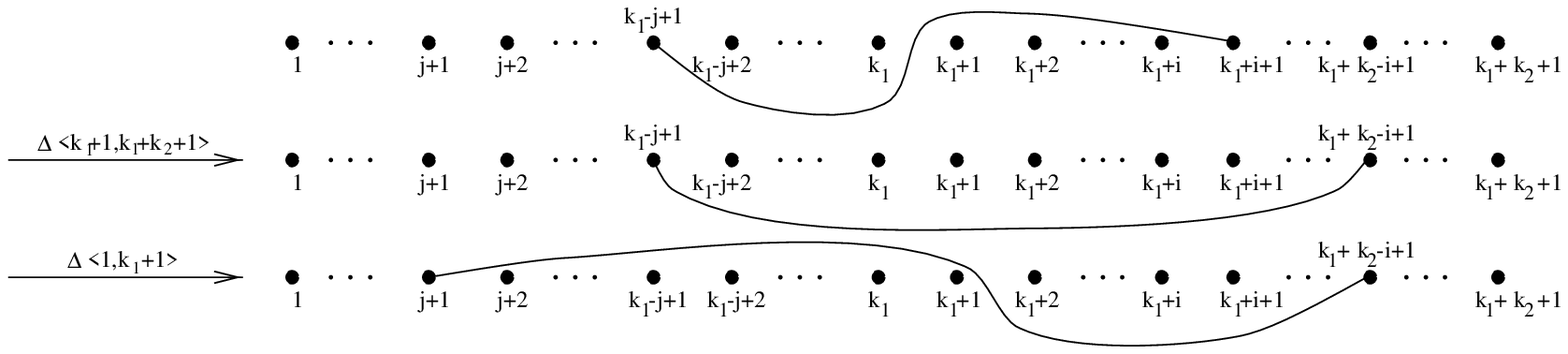}  
\end{center}

\medskip

\noindent
Then, the $l$ sequences of braids move the leftest side of the skeleton
$l$ points right:
  
\medskip

\begin{center}
\epsfxsize=14cm
\epsfbox{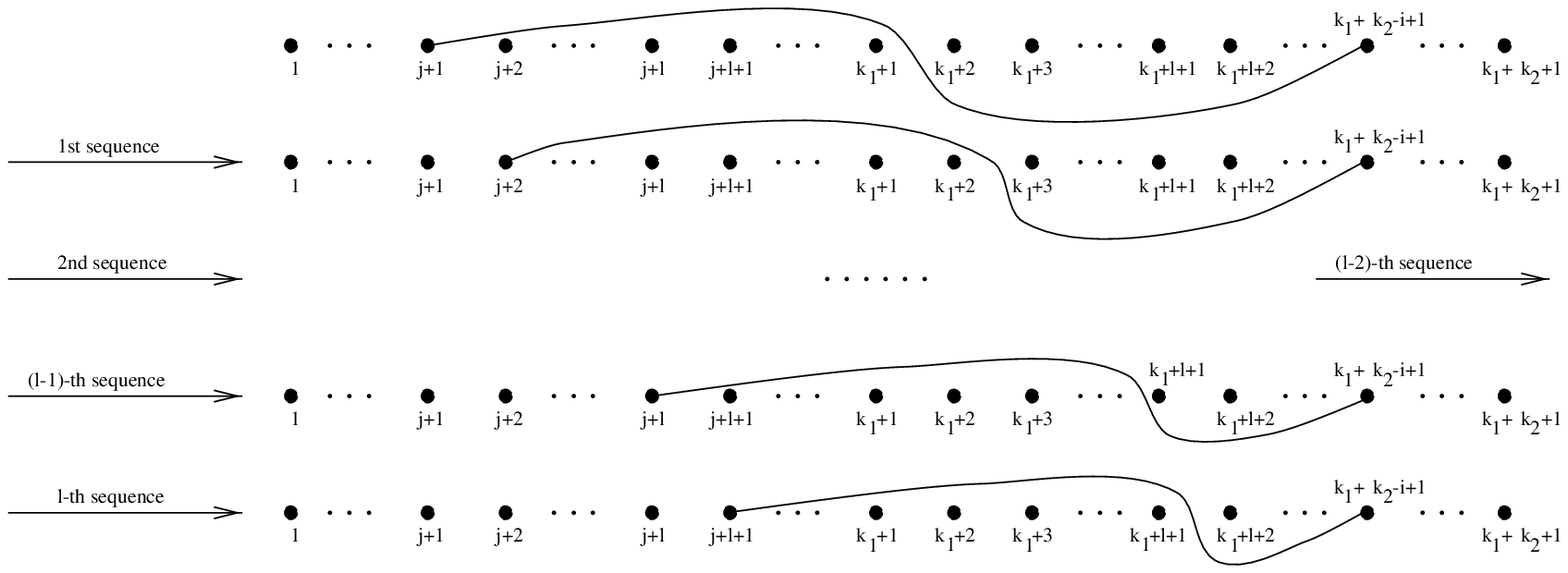}  
\end{center}

\medskip

After we have calculated the skeletons representing $\cL$.V.C.s for the 
braid monodromy, we can calculate the relations that they induced. 
As we have introduced in the previous section, according to 
Van-Kampen's theorem (\ref{VK2}), every $\cL$.V.C. induces a relation. 
Now, we will
calculate the general relations which are induced from the general $\cL$.V.C.s .

\medskip

The relation which is induced from $\varphi(\de _{i k_1+1}), 0 \leq i \leq l-1$:

\medskip

\begin{center}
\epsfxsize=11cm
\epsfbox{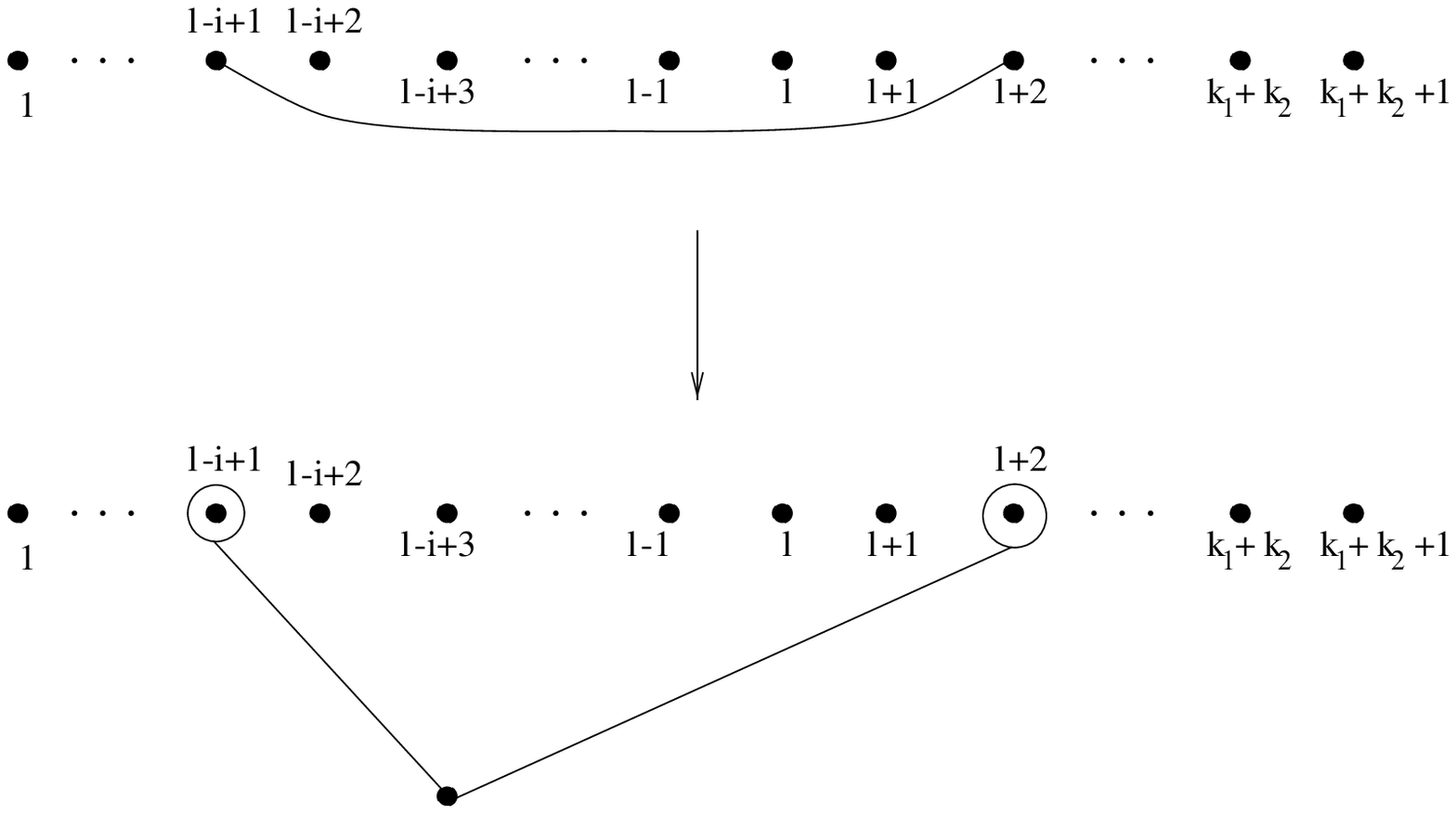}
\end{center}

\medskip

\noindent
Therefore, the relation is: 
$$\Ga _{l-i+1} \Ga _{l+2} = \Ga _{l+2} \Ga _{l-i+1}$$

The relation which is induced from 
$\varphi(\de _{i k_1+j}), \ 0 \leq i \leq l-1, \ 2 \leq j \leq k_1$:

\medskip

\begin{center}
\epsfxsize=13cm
\epsfbox{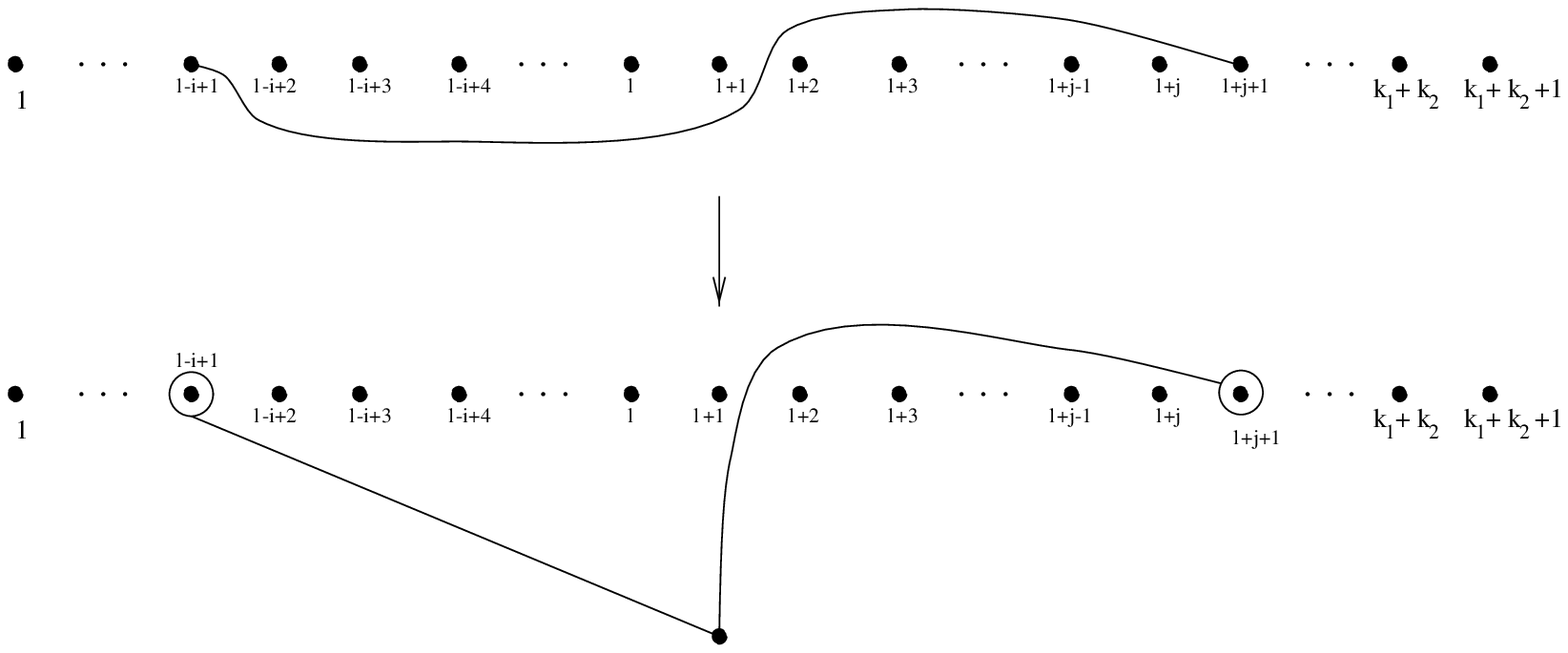}  
\end{center}

\medskip

\noindent
Therefore, the relation is: 
$$\Ga _{l-i+1} \Ga _{l+2} ^{-1} \cdots \Ga _{l+j} ^{-1} \Ga _{l+j+1} \Ga _{l+j} \cdots \Ga _{l+2} = \Ga _{l+2}^{-1} \cdots \Ga _{l+j} ^{-1} \Ga _{l+j+1} \Ga _{l+j} \cdots \Ga _{l+2} \Ga _{l-i+1}$$

The relation which is induced from 
$\varphi(\de _{(l k_1+2)+ik_1+1}), \ 0 \leq i \leq (k_2-l-1)$:

\medskip

\begin{center}
\epsfxsize=13cm
\epsfbox{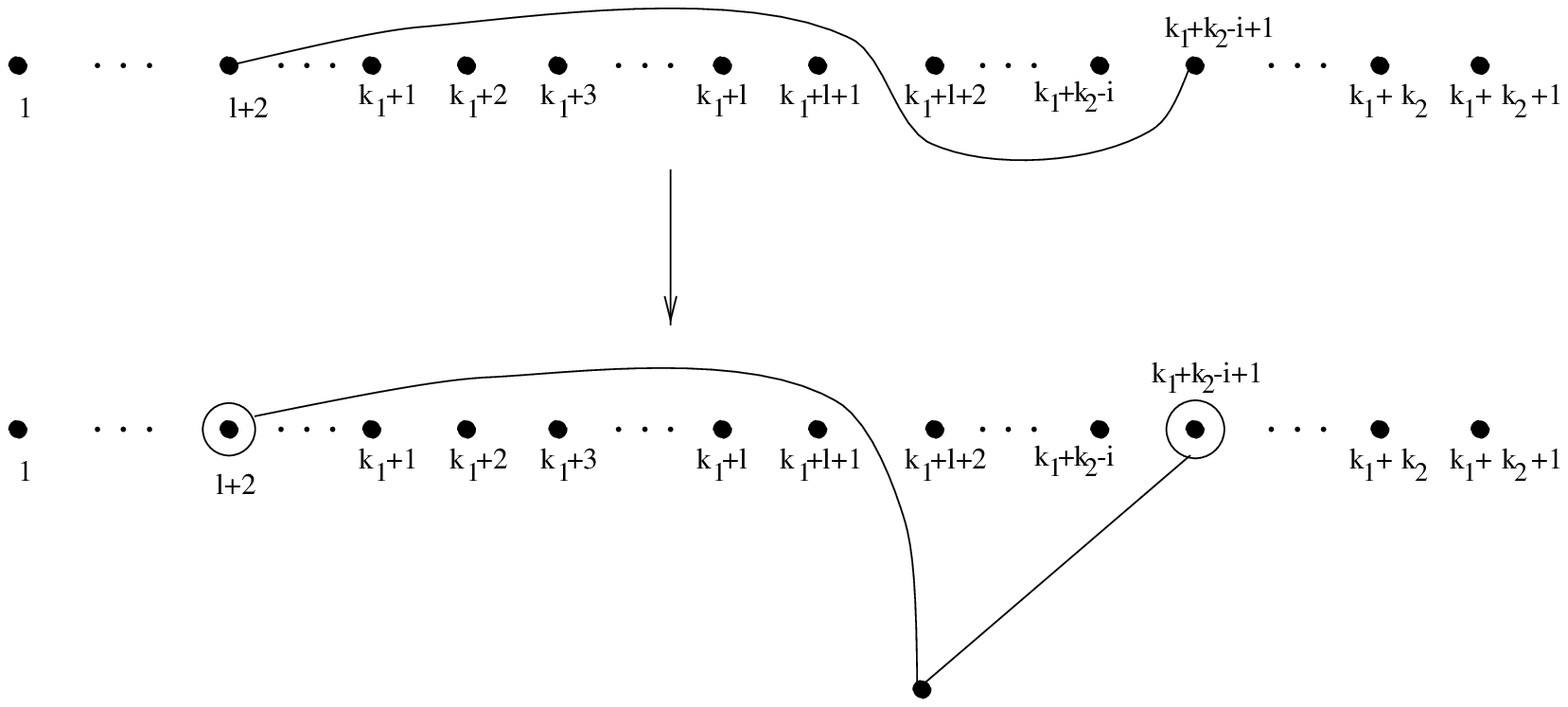}  
\end{center}

\medskip

\noindent
Therefore, the relation is:
$$\Ga _{k_1+l+1} \cdots \Ga _{l+3} \Ga _{l+2} \Ga _{l+3} ^{-1} \cdots \Ga _{k_1+l+1} ^{-1} \Ga _{k_1+k_2-i+1} =$$
$$ \Ga _{k_1+k_2-i+1} \Ga _{k_1+l+1} \cdots \Ga _{l+3} \Ga _{l+2} \Ga _{l+3} ^{-1} \cdots \Ga _{k_1+l+1} ^{-1}$$   

The relation which is induced from $\varphi(\de _{(l k_1+2)+ik_1+j}), 0 \leq i \leq (k_2-l-1), \\ 2 \leq j \leq k_1$:

\medskip

\begin{center}
\epsfxsize=11.5cm
\epsfbox{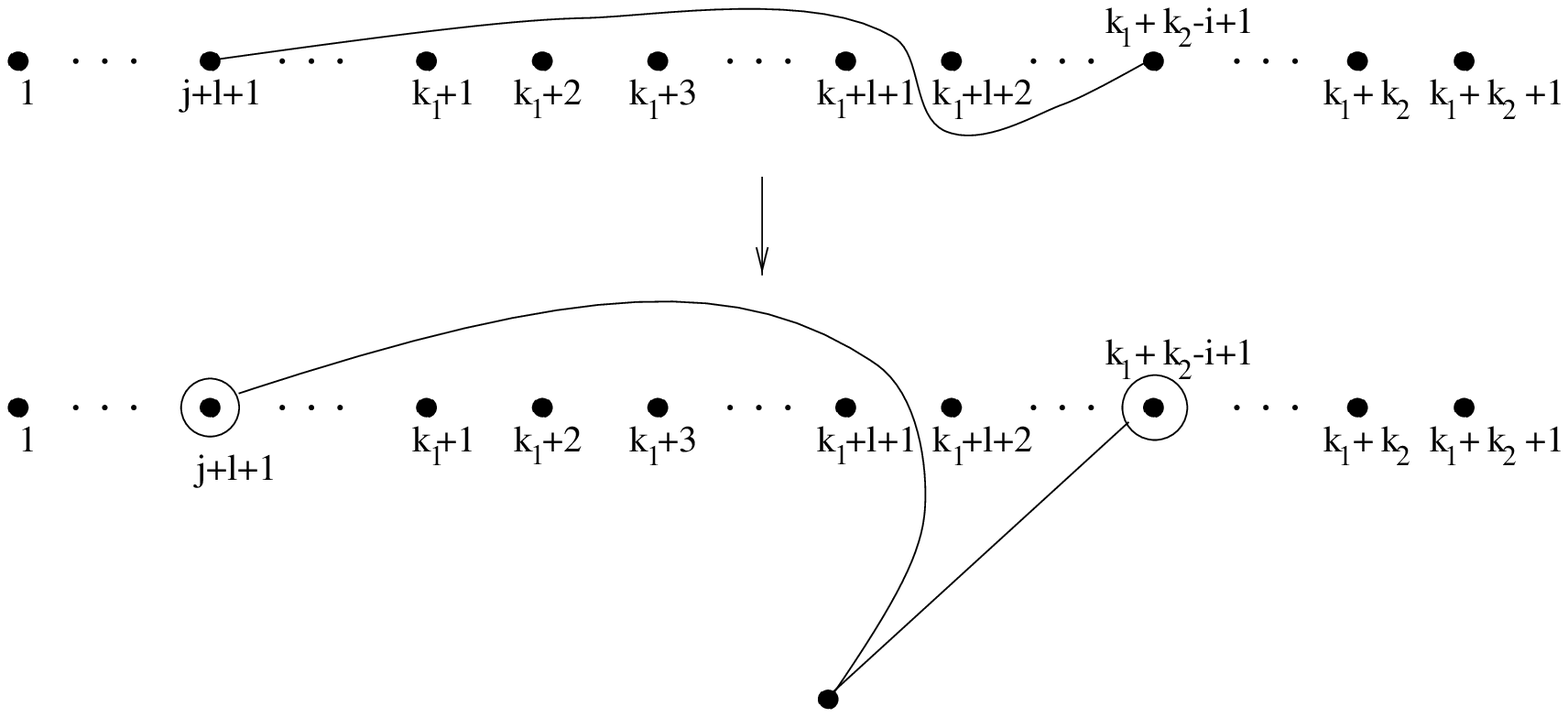}  
\end{center}

\medskip

\noindent
Therefore, the relation is:
$$\Ga _{k_1+l+1} \cdots \Ga _{l+j+2} \Ga _{l+j+1} \Ga _{l+j+2} ^{-1} \cdots \Ga _{k_1+l+1} ^{-1} \Ga _{k_1+k_2-i+1} =$$
$$ \Ga _{k_1+k_2-i+1} \Ga _{k_1+l+1} \cdots \Ga _{l+j+2} \Ga _{l+j+1} \Ga _{l+j+2} ^{-1} \cdots \Ga _{k_1+l+1}^{-1}$$   

\medskip

Therefore, we got the following two sets of relations: \\
for all $0 \leq i \leq l-1 ,\ 1 \leq j \leq k_1$:
$$\Ga _{l-i+1} \Ga _{l+2} ^{-1} \cdots \Ga _{l+j} ^{-1} \Ga _{l+j+1} \Ga _{l+j} \cdots \Ga _{l+2} =$$
$$\Ga _{l+2}^{-1} \cdots \Ga _{l+j} ^{-1} \Ga _{l+j+1} \Ga _{l+j} \cdots \Ga _{l+2} \Ga _{l-i+1}$$
and for all $0 \leq i \leq k_2 -l-1, \ 1 \leq j \leq k_1$:
$$\Ga _{k_1+l+1} \cdots \Ga _{l+j+2} \Ga _{l+j+1} \Ga _{l+j+2} ^{-1} \cdots \Ga _{k_1+l+1} ^{-1} \Ga _{k_1+k_2-i+1} =$$
$$ \Ga _{k_1+k_2-i+1} \Ga _{k_1+l+1} \cdots \Ga _{l+j+2} \Ga _{l+j+1} \Ga _{l+j+2} ^{-1} \cdots \Ga _{k_1+l+1}^{-1}$$

Now, it is easy to see that these two sets of relations are equivalent to 
the following two sets  of relations:
$$ \Ga _i \Ga _j = \Ga _j \Ga _i; \ 2 \leq i \leq l+1, \ l+2 \leq j \leq l+k_1+1$$
and
$$ \Ga _i \Ga _j = \Ga _j \Ga _i; \ l+2 \leq i \leq l+k_1+1, \ l+k_1+2 \leq j \leq k_1+k_2+1$$
and this finished the proof of the second case of the first lemma 
(\ref{lem1}). \hfill $\qed$

\subsection{Proof of lemma \ref{lem2}}\label{lem2-section}

As in the first lemma,  we prove this lemma only for two multiple points, 
and the proof for $t$ multiple points  uses exactly the same arguments.

We will prove it directly in the general case. By homotopic rotations and
movements and a proper choosing of the line at infinity, we can get 
the following line arrangement from any line arrangement with 
two multiple points:

\medskip

\begin{center}
\epsfxsize=8cm
\epsfbox{prop4_line_arr.ps}
\end{center}

\medskip

In the first lemma (\ref{lem1}), 
we already wrote down the set of Lefschetz pairs of this 
line arrangement. In order to calculate the induced
relations of the multiple points, we have to compute their braid monodromy
 according to the Moishezon-Teicher algorithm (\ref{MT})
and then we have to use the Van-Kampen theorem (\ref{VK2})
to get their induced relations.

\medskip

{\bf Skeleton representing the $\cL$.V.C. of $\varphi (\de _{l k_1+1})$}: 
The Lefschetz pair is 
$$(1,k_1+1),$$ 
then the skeleton representing the local $\cL$.V.C. is:

\medskip

\begin{center}
\epsfxsize=10cm
\epsfbox{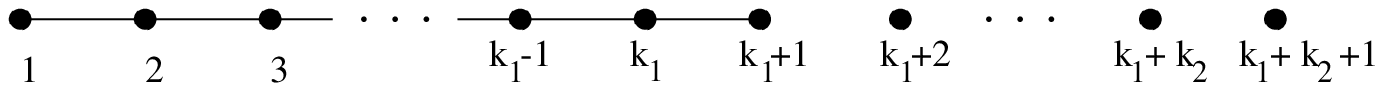}  
\end{center}

\medskip

\noindent
According to the algorithm, we have to apply on the skeleton 
the following sequence of braids:
$$\De <k_1+1,k_1+2> \De <k_1,k_1+1> \cdots \De <2,3>$$
$$\vdots$$
$$\De <k_1+l,k_1+l+1> \De <k_1+l-1,k_1+l> \cdots \De <l+1,l+2>$$

\noindent
The first sequence acts as follows:

\medskip

\begin{center}
\epsfxsize=11cm
\epsfbox{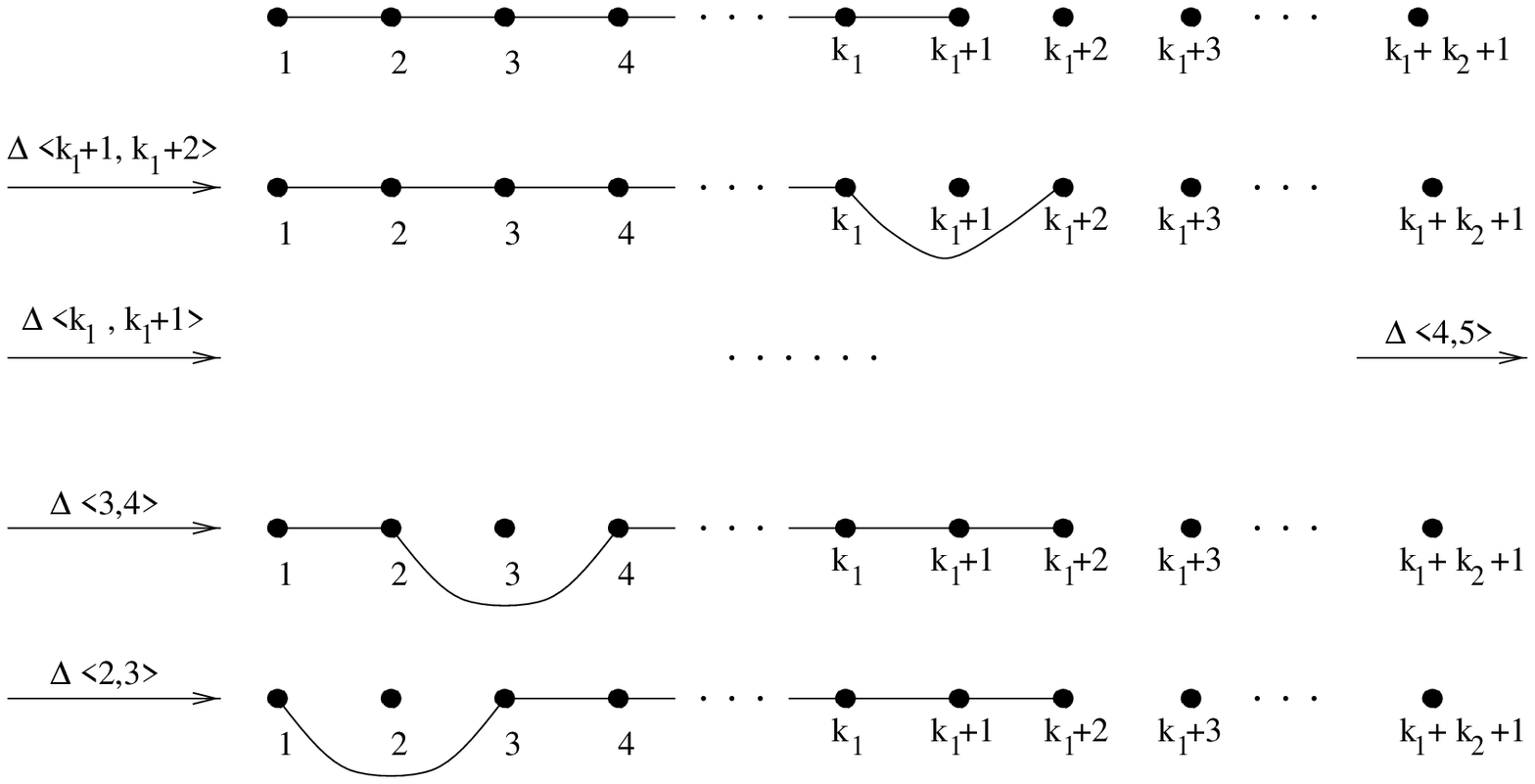}  
\end{center}

\medskip

Each of the next $l-1$ sequences moves the right side of the skeleton 
one step right, so we get the following: 

\medskip

\begin{center}
\epsfxsize=12.5cm
\epsfbox{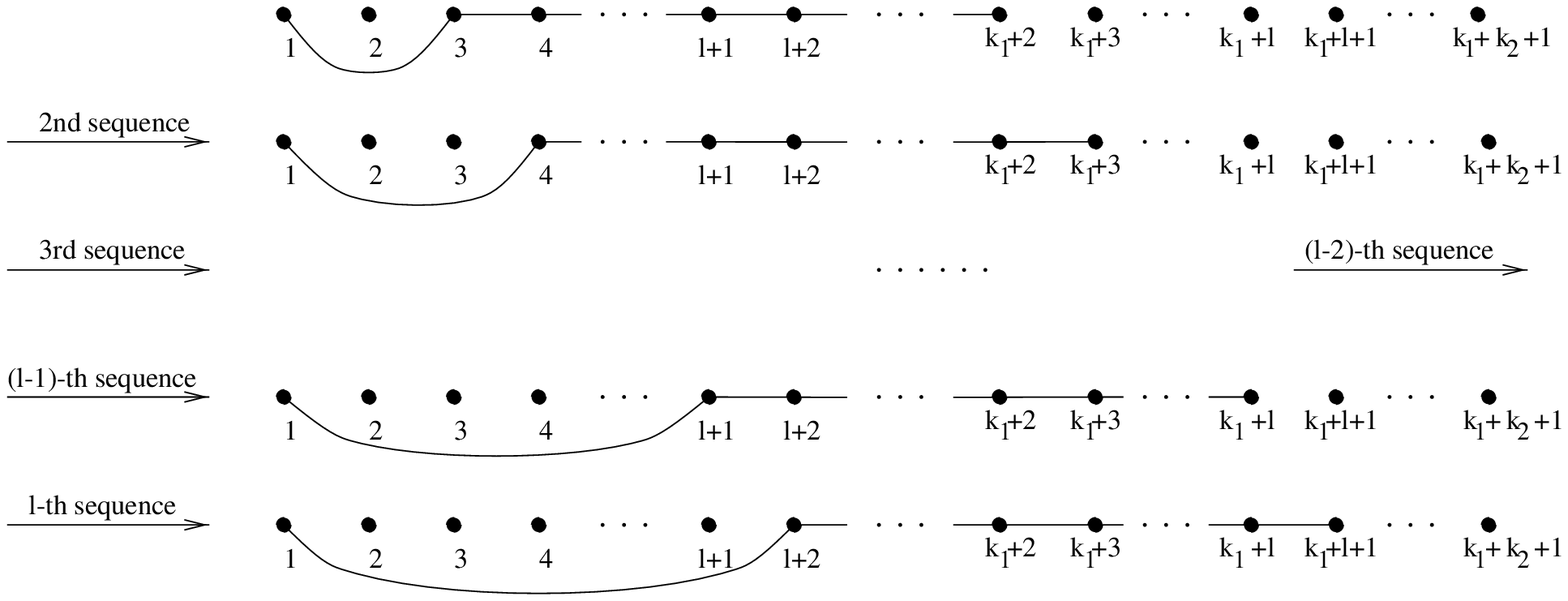}
\end{center}

\medskip

{\bf Skeleton representing the $\cL$.V.C. of $\varphi (\de _{l k_1+2})$}: 
The Lefschetz pair is 
$$(k_1+1,k_1+k_2+1),$$
therefore the skeleton representing the local $\cL$.V.C. is:

\medskip

\begin{center}
\epsfxsize=9.5cm
\epsfbox{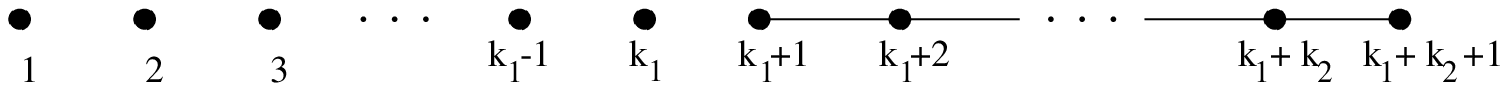}
\end{center}

\medskip

\noindent
According to the algorithm, we have to apply on the skeleton 
the following sequence of braids:
$$\De <1, k_1+1>$$
$$\De <k_1+1,k_1+2> \De <k_1,k_1+1> \cdots \De <2,3>$$
$$\vdots$$
$$\De <k_1+l,k_1+l+1> \De <k_1+l-1,k_1+l> \cdots \De <l+1,l+2>$$

The effect of the braid $\De <1,k_1+1>$ is:

\medskip

\begin{center}
\epsfxsize=10cm
\epsfbox{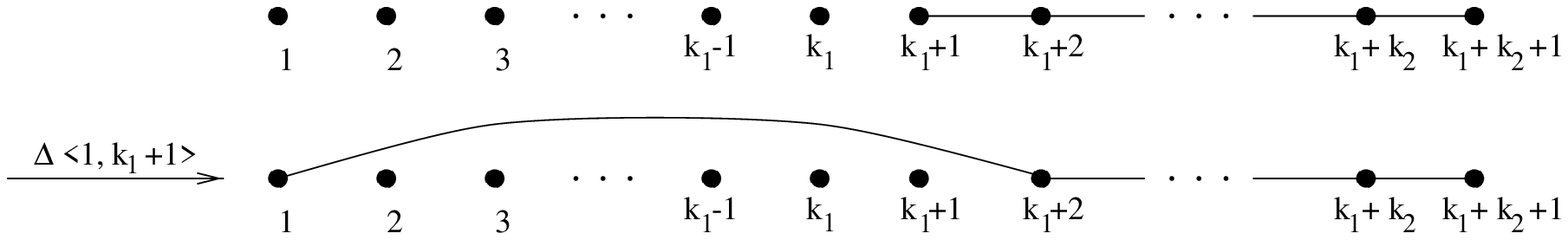}  
\end{center}

\medskip

\noindent
The first sequence acts as follows:

\medskip

\begin{center}
\epsfxsize=10cm
\epsfbox{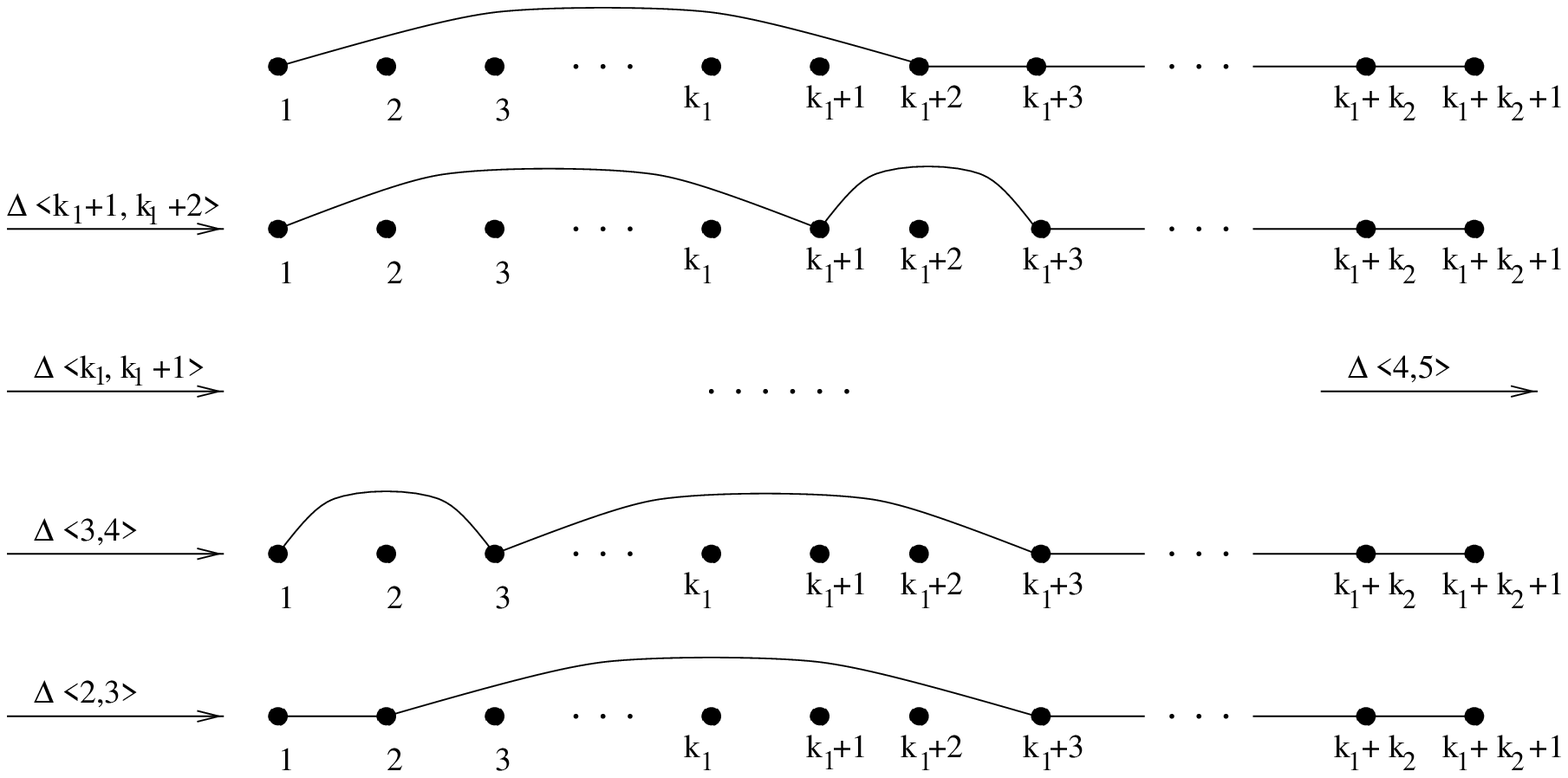}  
\end{center}

\medskip

Each of the next $l-1$ sequences moves the left side of the skeleton 
one step left, so we get the following: 

\medskip

\begin{center}
\epsfxsize=14cm
\epsfbox{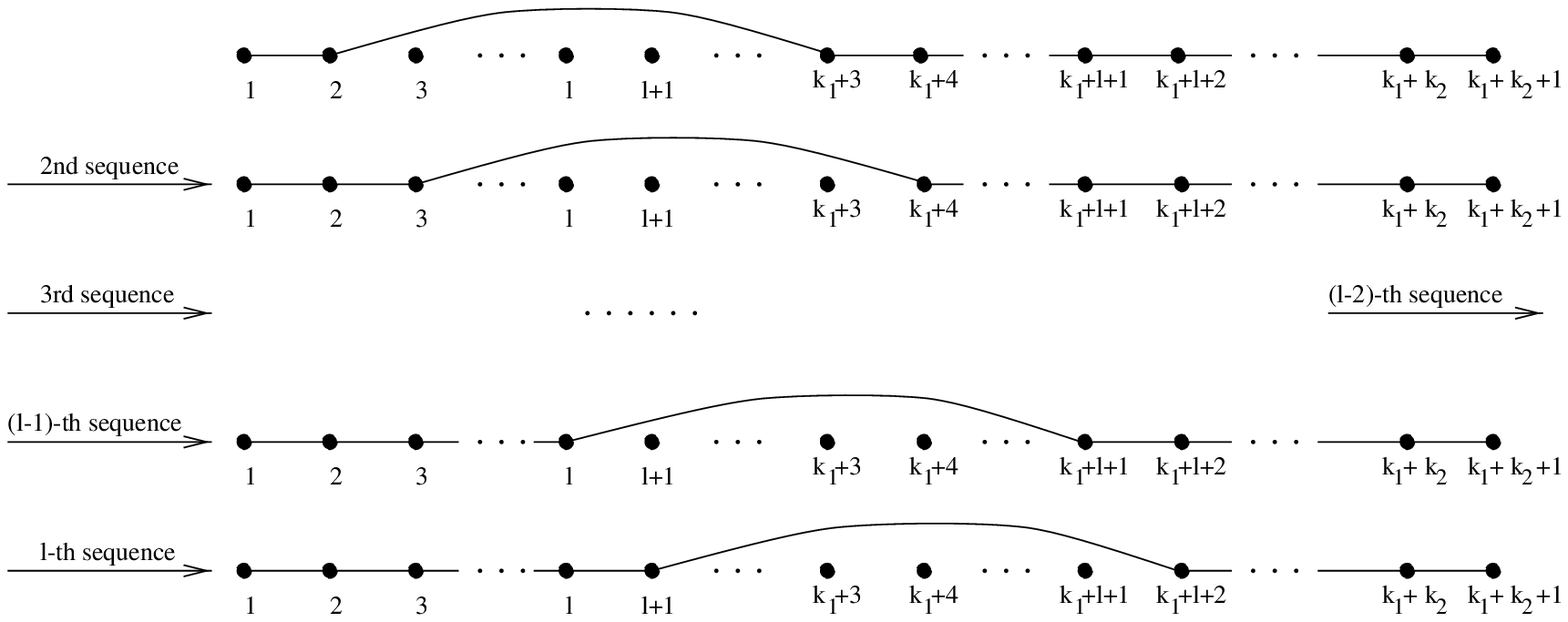}
\end{center}

\medskip

After we have calculated the skeletons representing $\cL$.V.C.s for the 
braid monodromy, we can calculate the relations which they induced. 

\medskip

The relations which are induced from $\varphi (\de _{l k_1+1})$:

\medskip

\begin{center}
\epsfxsize=13cm
\epsfbox{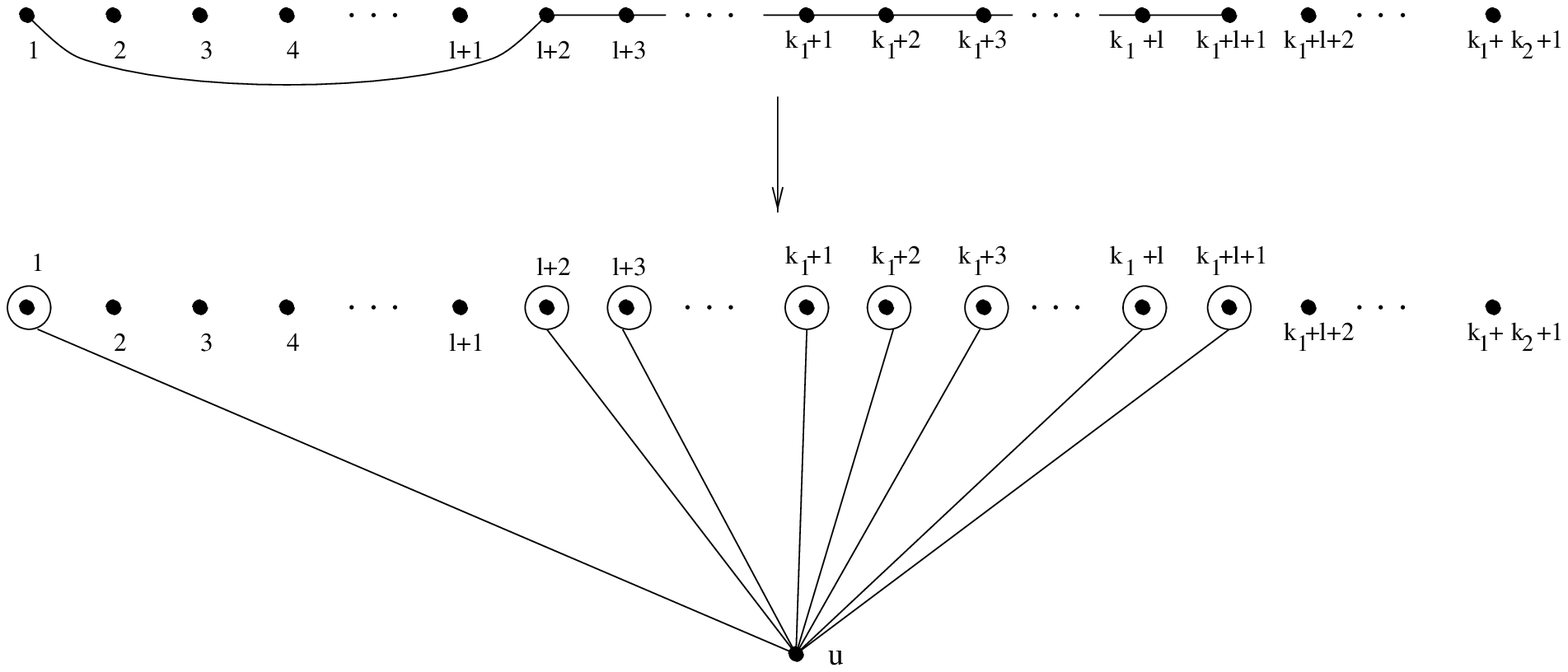}  
\end{center}

\medskip

\noindent
Therefore, according to lemma \ref{vk-cyc}, the relations are: 
$$\Ga _{k_1+l+1} \Ga _{k_1+l} \cdots \Ga _{l+2} \Ga _1 = \Ga _{k_1+l} \cdots \Ga _{l+2} \Ga _1 \Ga _{k_1+l+1} = \cdots = \Ga _1 \Ga _{k_1+l+1} \cdots \Ga _{l+2}$$
\medskip

The relations which are induced from $\varphi (\de _{l k_1+2})$:

\medskip

\begin{center}
\epsfxsize=13cm
\epsfbox{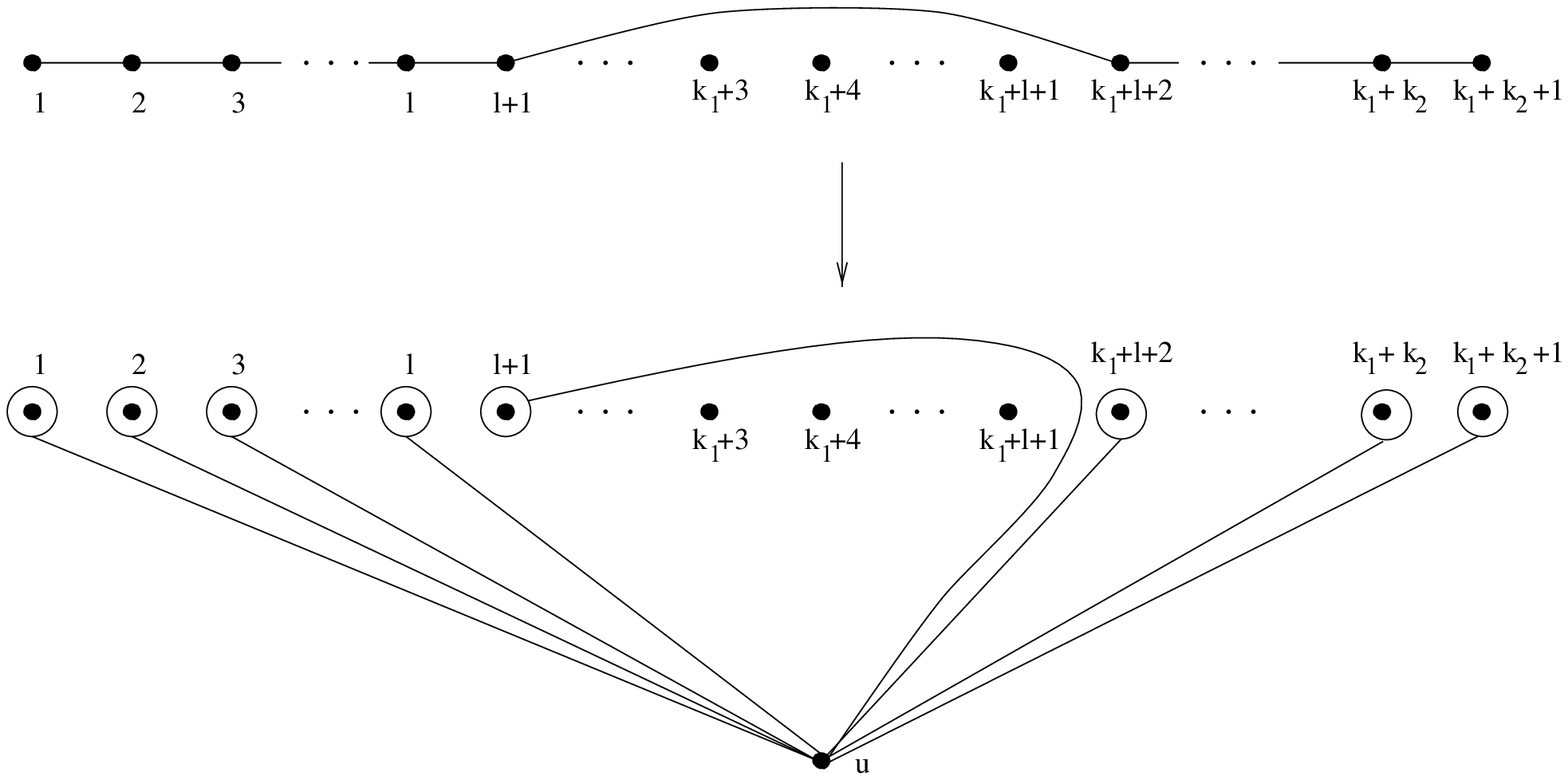}  
\end{center}

\medskip

\noindent
Therefore, according to lemma \ref{vk-cyc}, the relations are: 
$$\Ga _{k_1+k_2+1} \Ga _{k_1+k_2} \cdots \Ga _{k_1+l+2} (\Ga _{k_1+l+1} \cdots \Ga _{l+2} \Ga _{l+1} \Ga _{l+2} ^{-1} \cdots \Ga _{k_1+l+1} ^{-1}) \Ga _l \cdots \Ga _1 = $$
$$\Ga _{k_1+k_2} \cdots \Ga _{k_1+l+2} (\Ga _{k_1+l+1} \cdots \Ga _{l+2} \Ga _{l+1} \Ga _{l+2} ^{-1} \cdots \Ga _{k_1+l+1} ^{-1}) \Ga _l \cdots \Ga _1 \Ga _{k_1+k_2+1} = \cdots = $$
$$\Ga _1 \Ga _{k_1+k_2+1} \Ga _{k_1+k_2} \cdots \Ga _{k_1+l+2} (\Ga _{k_1+l+1} \cdots \Ga _{l+2} \Ga _{l+1} \Ga _{l+2} ^{-1} \cdots \Ga _{k_1+l+1} ^{-1}) \Ga _l \cdots \Ga _2$$

Now, according to the first lemma (\ref{lem1}), second case, 
$\Ga _{l+1}$ commutes with all $\Ga _j, \ l+2 \leq j \leq k_1+l+1$, therefore:
$$\Ga _{k_1+l+1} \cdots \Ga _{l+2} \Ga _{l+1} \Ga _{l+2} ^{-1} \cdots \Ga _{k_1+l+1} ^{-1} = \Ga _{l+1}.$$

Hence, the last set of relations comes to the following simplified form:
$$\Ga _{k_1+k_2+1} \Ga _{k_1+k_2} \cdots \Ga _{k_1+l+2} \Ga _{l+1} \Ga _l \cdots \Ga _1 = \Ga _{k_1+k_2} \cdots \Ga _{k_1+l+2} \Ga _{l+1} \Ga _l \cdots \Ga _1 \Ga _{k_1+k_2+1} = $$
$$ = \cdots = \Ga _1 \Ga _{k_1+k_2+1} \Ga _{k_1+k_2} \cdots \Ga _{k_1+l+2} \Ga _{l+1} \Ga _l \cdots \Ga _2$$
According to the proof of proposition \ref{lml}, these two sets of relations 
(of the two multiple points) are equivalent 
to the following two sets, respectively:
$$[\Ga _{k_1+l+1} \Ga _{k_1+l} \cdots \Ga _{l+2} \Ga _1, \Ga _i]=1, \ \forall i \in \{ 1, l+2, \cdots, k_1+l+1 \} $$
$$[\Ga _{k_1+k_2+1} \Ga _{k_1+k_2} \cdots \Ga _{k_1+l+2} \Ga _{l+1} \Ga _l \cdots \Ga _1, \Ga_i] =1, \ \forall i \in \{ 1 , \cdots, l+1, k_1+l+2, \cdots , k_1+k_2+1 \}$$
And the second lemma (\ref{lem2}) is proved. \hfill $\qed$

\subsection{Proof of  theorem \ref{thm_t_cnntd}}\label{lem3-section}
For simplicity, we prove the theorem only for two multiple points, 
and the proof for $t$ multiple points  uses exactly the same arguments.

Till now, we got the following set of generators:
$$g= \{ \Ga _1, \Ga _2, \cdots , \Ga _{k_1+k_2+1} \}$$
and the following sets of relations:
\begin{itemize}
\item[(1)] $ \Ga _i \Ga _j = \Ga _j \Ga _i; \ 2 \leq i \leq l+1, \ l+2 \leq j \leq k_1+l+1$
\item[(2)] $ \Ga _i \Ga _j = \Ga _j \Ga _i; \ l+2 \leq i \leq k_1+l+1, \ k_1+l+2 \leq j \leq k_1+k_2+1$
\item[(3)] $[\Ga _{k_1+l+1} \Ga _{k_1+l} \cdots \Ga _{l+2} \Ga _1, \Ga _i]=1, \ \forall i \in \{ 1, l+2, \cdots, k_1+l+1 \} $
\item[(4)] $[\Ga _{k_1+k_2+1} \Ga _{k_1+k_2} \cdots \Ga _{k_1+l+2} \Ga _{l+1} \Ga _l \cdots \Ga _1, \Ga_i] =1, \\ \forall i \in \{ 1 , \cdots, l+1, k_1+l+2, \cdots , k_1+k_2+1 \}$
\end{itemize}

We have to show that this finitely presented group is isomorphic to
$$\F ^{k_1} \oplus \F ^{k_2} \oplus \Z $$

Let us modify the set of generators by replacing the generator $\Ga _1$
by the generator 
$$\Ga ' = \Ga _{k_1+k_2+1} \Ga _{k_1+k_2}  \cdots \Ga _{k_1+1} \Ga _{k_1} \cdots  \Ga _2 \Ga _1$$

Now, we have to check that after the modifications 
we get an equivalent set of generators, and then we have to calculate 
the new set of relations. 

\medskip

\bcl
 After replacing $\Ga _1$ by $\Ga '$ 
 (which was defined above) in $g$, 
we again get  a set of generators. We denote 
this set of generators by $\tilde g$.
\ecl

\noindent
{\it Proof}: We have to show that $\Ga _1 \in <\tilde g>$. But this is
obvious, because: 
$$\Ga _1 = \Ga _2 ^{-1} \Ga _3 ^{-1} \cdots \Ga _{k_1+k_2+1}^{-1} \Ga '$$
\hfill $\qed$
   
\medskip

The next step is the calculation of the new set of relations for  $\tilde g$.
The sets (1) and (2) of the old sets of relations have not been changed 
(because these generators in the relations have not been replaced). 
We have to deal with the sets (3) and (4).

\medskip

\bcl  $[\Ga ',\Ga ] = 1, \forall \Ga \in \tilde g$. 
\ecl

\noindent
{\it Proof}: Obviously, $\Ga ' \Ga' = \Ga ' \Ga' $.
We will split the rest of the proof into two cases: \\
(a) $\Ga \in \{  \Ga _2, \cdots ,\Ga _{l+1}, \Ga _{k_1+l+2}, \cdots, \Ga _{k_1+k_2+1} \}$: \\ 
    $\Ga ' \Ga \stackrel{\rm Def}{=} \Ga _{k_1+k_2+1} \cdots \Ga _{k_1+l+2} \Ga _{k_1+l+1} \cdots \Ga _{l+2} \Ga _{l+1} \cdots \Ga _1 \Ga \stackrel{(2)}{=} \\  \Ga _{k_1+l+1} \cdots \Ga _{l+2} \Ga _{k_1+k_2+1} \cdots \Ga _{k_1+l+2} \Ga _{l+1} \cdots \Ga _1 \Ga \stackrel{(4)}{=} \\ \Ga _{k_1+l+1} \cdots \Ga _{l+2} \Ga \Ga _{k_1+k_2+1} \cdots \Ga _{k_1+l+2} \Ga _{l+1} \cdots \Ga _1 \stackrel{(1)(2)}{=} \\ \Ga \Ga _{k_1+l+1} \cdots \Ga _{l+2} \Ga _{k_1+k_2+1} \cdots \Ga _{k_1+l+2} \Ga _{l+1} \cdots \Ga _1 \stackrel{(2)+{\rm Def}}{=} \Ga \Ga ' $ 

\medskip

\noindent
(b) $\Ga \in \{  \Ga _{l+2}, \cdots ,\Ga _{k_1+l+1} \}$: \\ 
$\Ga ' \Ga \stackrel{\rm Def}{=} \Ga _{k_1+k_2+1} \cdots \Ga _{k_1+l+2} \Ga _{k_1+l+1} \cdots \Ga _{l+2} \Ga _{l+1} \cdots \Ga _1 \Ga \stackrel{(1)}{=} \\ \Ga _{k_1+k_2+1} \cdots \Ga _{k_1+l+2} \Ga _{l+1} \cdots \Ga _2 \Ga _{k_1+l+1} \cdots \Ga _{l+2} \Ga _1 \Ga \stackrel{(3)}{=} \\ \Ga _{k_1+k_2+1} \cdots \Ga _{k_1+l+2} \Ga _{l+1} \cdots \Ga _2 \Ga  \Ga _{k_1+l+1} \cdots \Ga _{l+2} \Ga _1 \stackrel{(1)+(2)}{=} \\ \Ga  \Ga _{k_1+k_2+1} \cdots \Ga _{k_1+l+2} \Ga _{l+1} \cdots \Ga _2 \Ga _{k_1+l+1} \cdots \Ga _{l+2} \Ga _1 \stackrel{(1)+{\rm Def}}{=} \Ga  \Ga '$ \hfill $\qed$

\medskip

Now, we can claim:

\bcl\label{clm_lem3} The following set is 
a complete set of relations for $\tilde g$ (we denote it by $\cR'$):
\begin{itemize}
\item[{\rm (1')}] $ \Ga _i \Ga _j = \Ga _j \Ga _i; \ 2 \leq i \leq l+1, \ l+2 \leq j \leq k_1+l+1$
\item[{\rm (2')}] $ \Ga _i \Ga _j = \Ga _j \Ga _i; \ l+2 \leq i \leq k_1+l+1, \ k_1+l+2 \leq j \leq k_1+k_2+1$
\item[{\rm (3')}] $[\Ga ',\Ga] = 1; \ \forall \Ga \in \tilde g$.
\end{itemize} 
\ecl

\noindent
{\it Proof}: 
We have to show that 
 $\{$(1),(2),(3),(4)$\}$ is equivalent to $\{$(1'),(2'),(3')$\}$ 
(with respect to the required replacements).
In the previous claim, we proved that $\{$(1),(2),(3),(4)$\}$ $\Rightarrow$ 
 $\{$(1'),(2'),(3')$\}$. We have to prove the opposite direction. 
Assume the set of relations $\{$(1'),(2'),(3')$\}$,
and prove the relations $\{$(1),(2),(3),(4)$\}$: \\
(1) and (2): this is the same as (1') and (2'), respectively.\\
(3): We have to prove that 
$$[\Ga _{k_1+l+1} \Ga _{k_1+l} \cdots \Ga _{l+2} \Ga _1, \Ga _i]=1, \ \forall i \in \{ 1, l+2, \cdots, k_1+l+1 \} $$
From (3') we know that $[\Ga ', \Ga _i]=1$. Therefore, we have: \\
$\Ga _{k_1+k_2+1} \Ga _{k_1+k_2} \cdots  \Ga _1 \Ga _i = \Ga _i \Ga _{k_1+k_2+1} \Ga _{k_1+k_2} \cdots  \Ga _1$ $\stackrel{(1')}{\Rightarrow}$ \\
$\Ga _{k_1+k_2+1} \cdots \Ga _{k_1+l+2} \Ga _{l+1} \cdots \Ga _2 \Ga _{k_1+l+1} \cdots \Ga _{l+2} \Ga _i = \\ = \Ga _i \Ga _{k_1+k_2+1} \cdots \Ga _{k_1+l+2} \Ga _{l+1} \cdots \Ga _2 \Ga _{k_1+l+1} \cdots \Ga _{l+2}$ $\stackrel{(1')+(2')+ (l+2 \leq i \leq k_1+l+1)}{\Longrightarrow}$ \\
$\Ga _{k_1+k_2+1} \cdots \Ga _{k_1+l+2} \Ga _{l+1} \cdots \Ga _2 \Ga _{k_1+l+1} \cdots \Ga _{l+2} \Ga _i = \\ = \Ga _{k_1+k_2+1} \cdots \Ga _{k_1+l+2} \Ga _{l+1} \cdots \Ga _2 \Ga _i \Ga _{k_1+l+1} \cdots \Ga _{l+2}$ $\stackrel{\Ga_2 ^{-1} \cdots \Ga _{l+1}^{-1} \Ga _{k_1+l+2} ^{-1} \cdots \Ga _{k_1+k_2+1} ^{-1} \cdot}{\Longrightarrow}$ \\
$[\Ga _{k_1+l+1} \Ga _{k_1+l} \cdots \Ga _{l+2} \Ga _1, \Ga _i]=1, \ \forall i \in \{ 1, l+2, \cdots, k_1+l+1 \} $

\medskip

Now, it remains to prove that:
$$[\Ga _{k_1+l+1} \Ga _{k_1+l} \cdots \Ga _{l+2} \Ga _1, \Ga _1]=1$$
$\Ga _{k_1+l+1} \cdots \Ga _{l+2} \Ga _1 \Ga _1 = \\
\Ga _{k_1+l+1}  \cdots \Ga _{l+2} \Ga _1 (\Ga _{l+2} ^{-1} \cdots \Ga _{k_1+l+1}^{-1} \Ga _{k_1+l+1}  \cdots \Ga _{l+2}) \Ga _1 \stackrel{(3)\  +\  ab=ba \Rightarrow ab ^{-1} =b^{-1} a}{=}\\
(\Ga _{l+2} ^{-1} \cdots \Ga _{k_1+l+1}^{-1} ) \Ga _{k_1+l+1}  \cdots \Ga _{l+2} \Ga _1 (\Ga _{k_1+l+1}  \cdots \Ga _{l+2}) \Ga _1 =  \Ga _1 \Ga _{k_1+l+1}  \cdots \Ga _{l+2} \Ga _1$ \\
(4) Same arguments as (3). \hfill $\qed$

\bigskip

We return to the proof of the theorem. Using the above claim, 
we can find the structure of the calculated group:

\medskip

\noindent
 $G = \langle g | \cR \rangle \cong \langle \tilde g | \cR' \rangle = \langle  \Ga _2, \cdots ,\Ga _{k_1+k_2+1}, \Ga ' \  \ | \ \cR'  \rangle \cong $\\
 $\cong \langle \Ga ' \rangle \oplus \langle  \Ga _2, \cdots ,\Ga _{k_1+k_2+1}  \ | \ \cR ' \rangle \cong $ \\
 $\cong \langle \Ga ' \rangle \oplus \langle \Ga _2, \cdots ,\Ga _{l+1}, \Ga _{k_1 +l+2}, \cdots, \Ga _{k_1+k_2+1}  \rangle \oplus\langle \Ga _{l+2}, \cdots ,\Ga _{k_1+l+1}  \rangle \cong$ \\
 $\cong \Z \oplus \F ^{k_2} \oplus \F ^{k_1}$ 

\noindent
Hence, we finished the proof of theorem \ref{thm_t_cnntd}. \hfill $\qed$

\subsection{The projective case}

Now, we will investigate the projective case.

\bthm\label{proj_t_cnntd} 
Let $\cL$ be a real line arrangement in $\C\PP ^2$
where all the $t$ multiple points are on the same line $L \in \cL$.
Let $k_i +1$ be the multiplicity of the
multiple point $P_i$, $1 \leq i \leq t$. Then:
$$\pi _1 (\C\PP ^2 -\cL) \cong {\bigoplus _{i=1} ^t} \F ^{k_i}$$
\ethm

\noindent
{\it Proof:} For simplicity, we will prove it for two multiple points 
and the proof for $t$  multiple points uses exactly the same arguments.

From the last theorem, we get:
$$\pi _1 (\C ^2 -\cL) \cong \F ^{k_1} \oplus \F ^{k_2} \oplus \Z$$
According to claim \ref{clm_lem3}, we get the following presentation 
for this group: \\
Generators: $g=\{ \Ga ', \Ga_2, \cdots, \Ga_{k_1+k_2+1} \}$. \\
Relations: $\cR = \{ \Ga _i \Ga _j = \Ga _j \Ga _i, \ 2 \leq i \leq l+1, \ l+2 \leq j \leq k_1+l+1; \\ \Ga _i \Ga _j = \Ga _j \Ga _i, \ l+2 \leq i \leq k_1+l+1, \ k_1+l+2 \leq j \leq k_1+k_2+1;[\Ga ',\Ga] = 1, \ \forall \Ga \in g \}$\\
where $\Ga ' = \Ga_{k_1+k_2+1} \cdots \Ga _1$ .

Now, when we are going to the projective case, we add one additional relation,
according to theorem \ref{VKP}:
$$\Ga_{k_1+k_2+1} \cdots \Ga _1 =1 $$
In terms of $\Ga '$, this relation gets the following form:
$$\Ga ' =1$$
Now,

\medskip

\noindent
$\pi_1 (\C\PP ^2 -\cL) = \langle g \ | \ \cR, \Ga ' =1 \rangle \cong \\  
\langle \Ga ' \ | \ \Ga ' =1 \rangle \oplus \langle \Ga _2, \cdots ,\Ga _{l+1}, \Ga _{k_1 +l+2}, \cdots, \Ga _{k_1+k_2+1}  \rangle \oplus\langle \Ga _{l+2}, \cdots ,\Ga _{k_1+l+1}  \rangle \cong$ \\
 $\cong   \F ^{k_2} \oplus \F ^{k_1}$ \hfill $\qed$
 
\bigskip

As a consequence of the last theorem, we get:

\bco
$$\pi _1 (\C ^2 -\cL) \cong \pi _1 (\C\PP ^2 -\cL)\oplus \Z $$
\eco

Therefore, the short exact sequence which was proved by Oka (theorem \ref{oka}):
$$ 1 \to \Z \to \pi _1 (\C ^2 -\cL) \to \pi _1 (\C\PP ^2 -\cL) \to 1 $$  
splits.

\section{Arrangements with more than one equivalence class}
\subsection{The definition of the equivalence relation}

The above results can be generalized more. Let us define the following relation
on the set of multiple intersection points: 

\bde
Let $p_1,p_2$ be two
multiple intersection points. We  say that $p_1 \sim p_2$ if $p_1$
is connected to $p_2$ by a ``path'' which its vertices are multiple
intersection points.
\ede

\bcl
$\sim$ is an equivalence relation on the set of multiple 
intersection points.
\ecl

\noindent
{\it Proof:} {\bf Reflexive:} each point is connected to itself by
the empty path.\\ 
{\bf Symmetry:} if $p_1$ is connected to $p_2$ by a path $P$,
$p_2$ is connected to $p_1$ by $P^{-1}$ - the opposite path of $P$ (which
is also a path of multiple points).\\
{\bf Transitive:} if $p_1$ is connected to $p_2$ by $P$, and $p_2$ is connected 
to $p_3$ by $Q$, $p_1$ is connected to $p_3$ by $P \cdot Q$, which is the 
concatenation of $P$ and $Q$ and therefore it is a path of multiple points,
because $p_2$ itself is a multiple point too. \hfill $\qed$

\medskip

This equivalence relation induces equivalence classes on the set of  multiple 
intersection points.
We also want to show that this equivalence relation induces a partition on the
lines of the arrangement:
\bcl
Let $C_1= \{ p_1, \cdots, p_k \}$ be the multiple points of one equivalence class 
and $C_2 = \{ q_1, \cdots,q_l \}$ be the multiple points of another equivalence class.
Let $\cL_i$ be the set of lines which  pass through one of 
the multiple points in $C_i$.\\
 Then: $\cL _1 \cap \cL _2 = \emptyset$.
\ecl

\noindent
{\it Proof:} Assume, on the contrary, that there exists a line $L$, 
such that $L \in \cL _1 \cap \cL _2$. Therefore, $L \in \cL_1$ and $L \in \cL _2$. 
From the definitions of $\cL_1$ and $\cL _2$, 
there exist points $p \in C_1$ and $q \in C_2$ 
such that $L$ passes through $p$ and $q$. Therefore, $p \sim q$, and hence
$C_1=C_2$, a contradiction to the assumption that $C_1$ and $C_2$ are distinct
equivalence classes. \hfill $\qed$

\subsection{The affine case}

Now, we can claim the following:

\bthm\label{general_thm}
Let $\cL$ be a real line arrangement in $\C\PP ^2$ consists of $n$ lines. We
choose the line at infinity such that all the lines are intersected in $\C ^2$.
Assume that there are $k$ multiple intersection points $p_1, \cdots ,p_k$
with multiplicities $m_1, \cdots, m_k$ respectively. 
Assume also that all the multiple 
intersection points in every equivalence class are collinear, i.e. every
equivalence class contains a unique line which connects all the multiple
points of that class. 
Then: 
$$\pi _1 (\C ^2 -\cL, u_0 ) \cong {\bigoplus_{i=1}^k} \F ^{m_i -1} 
\oplus \Z ^{n-({\sum_{i=1}^k} (m_i -1))}$$     
The number of infinite cyclic groups is a sum of two numbers: the number of 
equivalence classes and the number of lines which have only simple intersection
points.
\ethm

\noindent
{\it Proof:}  Let $C_i, 1 \leq i \leq t$ be the different equivalence classes of 
multiple points. According to the last claim, we define $\cL_i$ to be
the lines which  pass through points in $C_i$. Let $l_1, \cdots, l_r$ be lines
which are not in any $\cL_i$ (which means that they do not pass through any
multiple point, or equivalently, they intersect all the other lines 
at  simple points only). 

In every $\cL_i$, we have a line $L_i$ which connects all the multiple points
in $C_i$. Therefore, according to theorem \ref{thm_t_cnntd}, we have:
$$\pi_1 (\C ^2 - \cL _i) = ({\bigoplus _{j=1} ^{n_i}} \F ^{m_{P_{i,j}}-1}) \oplus \Z$$
where $n_i = \# C_i$ and $m_{P_{i,j}}$ is the multiplicity of the j-th point 
in $C_i$, $1 \leq j \leq n_i$.

For $l_i$, we know:
$$\pi_1 (\C ^2 - l _i) = \Z$$

Now, we use the Oka-Sakamoto theorem (see section 2.1) to get:
$$\pi_1 (\C ^2 - \cL) = \pi_1 (\C ^2 - ({\bigcup _{i=1} ^t} \cL _i \cup {\bigcup _{i=1} ^r} l_i)) \cong$$
$$\cong  ({\bigoplus _{i=1} ^t} \pi_1 (\C ^2 - \cL _i)) \oplus ({\bigoplus _{i=1} ^r} \pi_1 (\C ^2 - l _i)) \cong $$
$$\cong ({\bigoplus _{i=1} ^t} ({\bigoplus _{j=1} ^{n_i}} \F ^{m_{P_{i,j}}-1}) \oplus \Z) \oplus ({\bigoplus _{i=1} ^r} \Z) \cong ({\bigoplus _{i=1} ^t}{\bigoplus _{j=1} ^{n_i}} \F ^{m_{P_{i,j}}-1})\oplus \Z ^{t+r}$$ 
 
It remains to show that this group is equal to  the group 
mentioned in the formulation of the theorem.
First, in the double sum, every multiple point appears exactly once, 
because it appears in only one equivalence class. Therefore:
$$({\bigoplus _{i=1} ^t}{\bigoplus _{j=1} ^{n_i}} \F ^{m_{P_{i,j}}-1}) =  {\bigoplus_{i=1}^k} \F ^{m_i -1}$$
Now we have to show that:
$$t+r = n-({\sum_{i=1}^k} (m_i -1))$$ 
Let $o_i$ be the number of lines in $\cL_i$. We know that 
$$({\sum _{i=1} ^t } o_i) + r = n$$ 
It is easy to see that:
$$o_i = ({\sum _{j=1} ^{n_i}} (m_{P_{i,j}}-1))+1,$$
because there is a unique line which connects all the multiple points in every
equivalence class.  

When we combine the last two equations, we get:
$${\sum _{i=1} ^t }{\sum _{j=1} ^{n_i}} (m_{P_{i,j}}-1) +t+r=n$$
As before, due to the fact that every multiple point appears exactly in one 
equivalence class, we get:
$${\sum _{i=1} ^t }{\sum _{j=1} ^{n_i}} (m_{P_{i,j}}-1) = {\sum_{i=1}^k} (m_i -1)$$
and therefore, we get:
$$t+r = n-({\sum_{i=1}^k} (m_i -1))$$ 
\hfill $\qed$

\subsection{The projective case}
Now, we will investigate the projective case.

\bthm\label{proj_thm} 
Let $\cL$ be a real line arrangement in $\C\PP ^2$ consists of $n$ lines. We
choose the line at infinity such that all the lines are intersected in $\C ^2$.
Assume that there are $k$ multiple intersection points $p_1, \cdots ,p_k$
with multiplicities $m_1, \cdots, m_k$ respectively. 
Assume also that all the multiple 
intersection points in every equivalence class are collinear, i.e.  every 
equivalence class contains a unique line which connects all the multiple
points of that class. 
Then: 
$$\pi _1 (\C\PP ^2 -\cL, u_0 ) \cong {\bigoplus_{i=1}^k} \F ^{m_i -1} 
\oplus \Z ^{n-1-({\sum_{i=1}^k} (m_i -1))}$$     
The number of infinite cyclic groups is a sum of two numbers: the number of 
equivalence classes minus 1 and the number of lines 
which have only simple intersection points.
\ethm

\noindent
{\it Proof:} This is the projective analogue of theorem \ref{general_thm}. We 
induce it using the same techniques as we induced  
theorem \ref{proj_non_cnntd} from theorem \ref{thm_not_cnntd}. \hfill $\qed$

\medskip

As a consequence of the last theorem, we get:

\bco
$$\pi _1 (\C ^2 -\cL) \cong \pi _1 (\C\PP ^2 -\cL)\oplus \Z $$
\eco

Therefore, the short exact sequence which was proved by Oka (theorem \ref{oka}):
$$ 1 \to \Z \to \pi _1 (\C ^2 -\cL) \to \pi _1 (\C\PP ^2 -\cL) \to 1 $$  
splits. 

\medskip

\noindent
{\bf Remark:} Simultaneously and independently, Fan [Fa2] got similar 
results (see section 2.1), with entirely different methods, in even 
more general case, when there is no equivalence class which has a cycle 
of multiple points in it.

\section{Results concerning the bigness of the  fundamental group}
\stepcounter{subsection}

\bde
A group $G$ is called {\bf big} if $\F ^2 \sbs G$. 
\ede

As a result from the general theorems (\ref{general_thm},\ref{proj_thm}), 
we can say the following:

\bco\label{bigness}
Let $\cL$ be a real line arrangement in $\C\PP ^2$ consisting of 
$n$ lines which satisfies the conditions of theorem \ref{general_thm}.  
Then, the fundamental groups of its complement, $\pi _1 (\C ^2 -\cL,u_0)$ and
$\pi _1 (\C\PP ^2 -\cL,u_0)$, are big. 
\eco

\noindent
{\it Proof}: According to theorem \ref{general_thm}, the fundamental group of its
affine complement is of the form:
$$\pi _1 (\C ^2 -\cL, u_0 ) \cong {\bigoplus_{i=1}^k} \F ^{m_i -1} 
\oplus \Z ^{n-({\sum_{i=1}^k} (m_i -1))}$$     
Now, $m_i \geq 3$ in every multiple point, and hence  $\F^2$ is contained 
in  this group. Therefore, the fundamental group of its affine complement 
is big. The proof for the projective case is the same. \hfill $\qed$

\medskip

In fact, this result has been recently proven [DOZ] for any arrangement
which has at least one multiple intersection point:
  
\bthm\label{thm_doz} {\bf (Dethloff, Orevkov, Zaidenberg)} \\
Let $\cL$ be a real line arrangement in $\C\PP ^2$ consisting of $n$ lines. We
choose the line at infinity such that all the lines are intersected in $\C ^2$.
Assume that there exists in $\cL$ at least one multiple intersection point.\\
Then, $\pi _1 (\C\PP ^2 - \cL, u_0)$ is big.
\ethm

\medskip

\noindent
{\bf Remark:} It seems that this phenomena is not happen for branch curves of
surfaces, unlike previous expectations which followed earlier results of
Zariski and Moishezon. Most fundamental groups of complements of branch 
curves are ``almost solvable'', i.e. they contain a solvable subgroup of
finite index and they are not ``big'' (see [Te2]).   

\bigskip\bigskip

\noindent
{\it Acknowledgments.} We  thank Prof. Leonid Makar-Limanov
 for suggestions which led to the crucial 
part of the proof of proposition \ref{lml}.

\bigskip\bigskip

\begin{\bib}{10}
\bibitem[CS]{CS} Cohen, D.C. and Suciu, A.I., {\it The braid monodromy of 
   plane algebraic curves and hyperplane arrangements}, Comment. Math. 
   Helv. {\bf 72(2)}, 285-315 (1997).
\bibitem[DOZ]{DOZ} Dethloff, G., Orevkov, S. and Zaidenberg, M., 
  {\it Plane curves with a big fundamental group of the complement}, 
{\it in:} Voronezh Winter Mathematical Schools: Dedicated
to Selim Krein (P. Kuchment, V. Lin, eds.), 
American Mathematical Society Translations--Series 2 {\bf 184} (1998).
\bibitem[Fa1]{Fa1} Fan, K.M., {\it Position of singularities and fundamental group of the
   complement of a union of lines}, Proc. Amer. Math. Soc. {\bf 124(11)}, 
   3299-3303 (1996).
\bibitem[Fa2]{Fa2} Fan, K.M., {\it Direct product of free groups as the fundamental 
   group of the complement of a union of lines}, Michigan Math. J. 
   {\bf 44(2)}, 283-291 (1997).  
\bibitem[Ga]{Ga} Garber, D., {\it On the fundamental group of complement of real line 
   arrangements}, M.Sc. thesis, Bar-Ilan University (1997).
\bibitem[Mo]{Mo} Moishezon, B., {\it Stable branch curves and braid monodromies}, 
  Lect. Notes in Math. {\bf 862}, 107-192 (1981). 
\bibitem[MoTe1]{MoTe1} Moishezon, B. and  Teicher, M., {\it Braid group 
   techniques in complex geometry I, Line arrangements in $\C\PP ^2$}, 
   Contemporary Math. {\bf 78}, 425-555 (1988).  
\bibitem[MoTe2]{MoTe2} Moishezon, B. and  Teicher, M., {\it Braid group 
   techniques in complex geometry II, From arrangements of lines and conics 
   to cuspidal curves}, Algebraic Geometry, Lect. Notes in Math. {\bf 1479} 
   (1990).
\bibitem[MoTe3]{MoTe3} Moishezon, B. and  Teicher, M., {\it Braid group
   techniques in complex geometry V: The fundamental group of a complement of
   a branched curve of a Veronese generic projection}, Com. in Analysis and
   Geometry {\bf 4(1)}, 1-120 (1996).  
\bibitem[MoTe4]{MoTe4} Moishezon, B. and  Teicher, M., {\it Braid groups, 
   singularities, and algebraic surfaces}, Academic Press, to appear.
\bibitem[O]{O} Oka, M., {\it On the fundamental group of a reducible curve in $\PP ^2$},
   J. London Math. Soc. (2) {\bf 12}, 239-252 (1976).
\bibitem[OS]{OS1} Oka, M. and Sakamoto, K., {\it Product theorem of the 
   fundamental group of a reducible curve}, J. Math. Soc. Japan {\bf 30(4)}, 
   599-602 (1978).
\bibitem[OT]{OT} Orlik, P. and Terao, H., {\it Arrangements of hyperplanes},
   Grundlehren {\bf 300}, Springer-Verlag (1992).
\bibitem[Ra]{Ra} Randell, R., {\it The fundamental group of the complement of a union of
  complex hyperplanes}, Invent. Math. {\bf 69}, 103-108 (1982). {\it Correction},
  Invent. Math. {\bf 80}, 467-468 (1985).
\bibitem[Sa]{Sa} Salvetti, M., {\it Topology of the complement of real hyperplanes in 
  $\C  ^N$}, Invent. Math. {\bf 88}, 603-618 (1987). 
\bibitem[Te1]{Te1} Teicher, M., {\it Braid groups, algebraic surfaces and fundamental
   groups of complement of branch curves}, Proc. Symp. Pure Math. 
   {\bf 62(1)}, 127-150 (1997).
\bibitem[Te2]{Te2} Teicher, M., {\it New invariants of surfaces}, Contemp.
Math., to appear. 
\bibitem[VK]{VK} Van Kampen, E.R., {\it On the fundamental group of 
   an algebraic curve}, Amer. J. Math. {\bf 55}, 255-260 (1933). 
\bibitem[Z]{Z} Zariski, O., {\it On the problem of existence of algebraic functions
   of two variables possessing a given branch curve}, Amer. J. Math. {\bf 51},
   305-328 (1929).
\end{\bib}

\end{document}